\documentclass[10pt,UTF8]{article}
\usepackage{amsmath,mathrsfs}
\usepackage{amsfonts}
\usepackage{amssymb,amsthm,mathtools}
\usepackage{graphicx}
\usepackage{xcolor}
\usepackage{subfigure,enumitem}
\usepackage{cite}

\usepackage{caption}

\usepackage{indentfirst}
\usepackage[numbers,sort&compress]{natbib}

\usepackage{hyperref}
\hypersetup{
colorlinks=true,
linkcolor=red,
urlcolor=blue,
citecolor=blue}

\usepackage[capitalise, nameinlink]{cleveref}

\usepackage{ulem}
\usepackage{cases}

\newtheorem{theorem}{Theorem}[section]
\newtheorem{corollary}{Corollary}[section]
\newtheorem{lemma}{Lemma}[section]
\newtheorem{example}{Example}[section]
\newtheorem{definition}{Definition}[section]
\newtheorem{remark}{Remark}[section]
\newtheorem{proposition}{Proposition}[section]
\numberwithin{equation}{section}

\def\to{{\rightarrow}}
\allowdisplaybreaks

\begin{document}

\title{
    Evolution Dynamics and Cauchy Problems with Quaternionic Variables and Lifting to 
    Banach Bilateral $\mathbb{A}$-Algebras
}
\author{
    Chao Wang$^a$\thanks{Email: chaowang@ynu.edu.cn.},\,
 	Guangzhou Qin$^{a}$, \,Jibin Li$^{b,c}$\\ \\
 	$^a$Department of Mathematics, Yunnan University\\
 	Kunming, Yunnan 650091, China \\ \\
 	$^b$ School of Mathematical Sciences, Zhejiang Normal University\\
Jinhua, Zhejiang 321004, China\\ \\
$^c$School of Mathematical Science, Huaqiao University\\
Quanzhou, Fujian 362021, China\\
}
\date{}
\maketitle{}

\begin{abstract}
    In this paper, we develop a theory of evolution dynamics in a Banach bilateral
    $\mathbb{H}$-algebra and establish a well-posedness theory for the Cauchy problem of
    evolution equations with quaternionic variables. First, we study the geometric
    properties of slice structures of Banach bilateral $\mathbb{H}$-modules.
    By employing different slice structures, we construct some noncommutative Banach
    bilateral $\mathbb{H}$-subalgebras which are compatible with the slice regularity of
    the functions involved. Moreover, through the quaternionic sliced geometric path
    integral over Banach bilateral $\mathbb{A}$-module, an integral representation of
    the solution to the homogeneous Cauchy problem with quaternionic variables is
    obtained in a quaternionic $\sigma$-ball, and for the first time, a local existence
    and uniqueness theorem is established for $\sigma$-analytic solutions.
    Then the local $\sigma$-analytic solution is extended to the entire quaternionic
    spherical sectorial region, thereby proving a global existence and uniqueness
    theorem for the Cauchy problem. Based on it, the notion of an evolution system
    in a quaternionic $\mathbb{H}$-algebra is introduced and its core properties
    are established. In addition, the global existence and uniqueness theorems
    are obtained under hyperbolic-like and parabolic-like conditions.
    Finally, the structure of Banach bilateral $\mathbb{H}$-modules and
    $\mathbb{H}$-algebras is lifted to Banach bilateral $\mathbb{A}$-modules
    and $\mathbb{A}$-algebras, and a well-posedness theory for evolution equations in
    this setting is developed. This reveals, from a new perspective, the universality
    of well-posedness for the Cauchy problem under the generic noncommutative
    $\mathbb{A}$-algebraic structures.
\end{abstract}


\indent \textbf{2020MSC:} 47A10; 47D06; 47S05; 34G10; 30G35.

\tableofcontents

\section{Introduction}

The concept of quaternions was first introduced by the Irish mathematician Hamilton in
1843 (see \cite{ha1}), and has since found extensive application in both pure and
applied mathematics (see \cite{kexi-basp0,kexi-xdao4,kexi-qpzz2}),
as well as in physics, engineering, and other interdisciplinary domains (see \cite{cwgu1}).
In 1936, Birkhoff and Neumann demonstrated that quantum mechanics can be formulated over the real,
complex, or quaternionic number systems (see \cite{Bi1}),
a result that subsequently established quaternionic theory as a powerful tool widely adopted
across numerous important research areas (see \cite{chen1,qin1,liz1,ggg1}).
During the 1980s and 1990s, McIntosh, Jefferies, and their colleagues carried
out pioneering work on functional calculus within the framework of Clifford
algebras (see \cite{70,71,73,69}). In 1995, Adler provided a comprehensive exposition
of the theory underlying quaternionic quantum mechanics and quantum fields (see \cite{adj1}).
However, the precise proof of the spectral theorem for quaternionic linear operators
remained an unresolved issue even after the formulation of quaternionic quantum mechanics,
primarily due to the lack of a clear definition for the quaternionic spectrum and the
associated resolvent operator. A breakthrough did not occur until 2006, when Colombo and
Sabadini introduced the notion of the $S$-spectrum, which laid the foundation for a
comprehensive quaternionic operator theory (see Subsection 1.2.1 in \cite{kexi-co3}
for details on the discovery of the $S$-spectrum and how this novel spectral notion emerged).
The development of a spectral theorem based on the $S$-spectrum required several additional years;
finally, in 2016, Alpay et al. (see \cite{Al1}) provided a complete proof of this fundamental
theorem for both bounded and unbounded operators. The $S$-functional calculus, along with the
broader spectral theory on the $S$-spectrum, began to take shape only after 2006,
and the principal references for this theory and its associated function theory
are the monographs \cite{Al5,kexi-co5,kexi-co3,kexi-x8} and the references therein.
In \cite{colo1}, the authors established the universality property of the $S$-functional
calculus for noncommuting matrix variables and Clifford operators. Nonetheless,
the noncommutative nature of quaternionic theory distinguishes it from complex
analysis (see \cite{26em}) and poses substantial challenges for ongoing research in this field,
one may refer to Alpay et al.\cite{Al2,kexi-xxx9}, Cerejeiras et al.\cite{cer1,cer2},
Colombo et al.\cite{jia6,kexi-co1}, Gantner \cite{kexi-basp2} and references therein
for more detailed progress of theory on quaternionic operator and $S$-spectrum.

Building upon the development of quaternionic operator theory and $S$-spectrum,
it is therefore natural to turn to the theory of slice regular functions,
which has emerged as the fundamental analytic tool underpinning these advances
and enabling a coherent extension of complex function theory to the quaternionic setting.
The subtlety of the geometric structure of quaternionic slices lies in their ability
to generate more intricate and refined noncommutative subspaces,
on each of which slice regular functions can be defined (see Alpay et al.\cite{kexi-srs2}).
Some important fundamental theories of slice regular functions,
such as the Cauchy kernel (see \cite{kexi-cxh0}), the extension theorem (see \cite{jia4}),
the Cauchy transform (see \cite{28}), and the duality theorem (see \cite{17}),
have greatly enriched the theoretical framework of quaternionic slice hyperholomorphic functions.
Moreover, in 2023, Dou et al. (see \cite{mb11}) extended the extension theorem
and representation formulas for regular functions to quaternionic non-axially-symmetric domains,
achieving a series of significant advances,  and simultaneously obtained important progress
for a class of slice monogenic functions (see \cite{jia5,jia7}).
As early as 2006, Gentili and Struppa (see \cite{kexi-qpzz0})
announced the basic elements of a new theory of regular functions of one quaternionic variable,
following an idea of Cullen but adopting a more geometric approach to demonstrate that
a rather complete theory can be built. They then developed the fundamental results in full,
showing that this theory allows one to extend important results for quaternionic polynomials
to the case of power series (see \cite{kexi-qpzz1}).
This theory of quaternionic regular functions was subsequently refined 
(see \cite{kexi-x8,kexi-qpzz3,kexi-daoshu0})
and into an important theoretical framework, encompassing
key achievements such as local representation formulas (see \cite{gen4}),
quaternionic power series and analyticity theory (see \cite{kexi-jishu}).
In particular, it also played an important role in the study of orthogonal complex
structures on slice domain (see \cite{gee1}). In addition, in 2012, Stoppato
(see \cite{stopo,18c}) established the theory of singularities and series
extension for slice regular functions.

The structure of quaternionic slice regular functions can be extended to broader
algebraic frameworks while retaining favorable generic properties, a fact profoundly
demonstrated by the significant work of Ghiloni et al. In 2011, Ghiloni and Perotti
(see \cite{jia11}) developed a theory of slice regular functions on real alternative algebras,
and within the context of real alternative $^*$-algebras, they further established theories
of power series and spherical series (see \cite{ghi3}), as well as a singularity theory for
slice regular functions (see \cite{gh1}). Moreover, they extended slice regularity to more
general noncommutative algebraic structures by developing the theory of slice regular
semigroups and that of quaternionic spherical sector operators on such algebras (see \cite{ghp1}).
This framework adapts slice regularity to generalized Banach bilateral modules whose
scalar sets can be any real associative $^*$-algebra or Clifford algebra
(see \cite{kexi-srs0,kexi-srs1}).

On the other hand, the Cauchy problem (see \cite{hada1,kexi-8qq}),
a term introduced by Hadamard in his 1921 Yale University Press monograph,
stands as a central theme in the theory of partial differential equations.
Its significance spans both the classical era dating from the late nineteenth
to early twentieth centuries and the modern period, which has witnessed remarkable
progress since recent decades. In its classical formulation, the Cauchy problem
emerged in two distinct contexts: first, as the foundation for the most essential
results in the theory of analytic partial differential equations; and second,
in connection with classical boundary value problems, particularly those arising
in the study of the wave equation and the broader class of second-order hyperbolic equations.
With the progressive deepening and advancement of research on the Cauchy problem,
evolution equations in abstract spaces, along with semigroup theory, operator theory,
spectral theory, and related disciplines, have consequently evolved into
a systematic and integrated theoretical framework
(see \cite{hi3,hil1,kat1,kat2,kat3,kexi-8q,phi1}).
Regarding quaternionic evolution equations, Colombo and Sabadini (see \cite{kexi-co1})
proposed the following Cauchy problem as an open issue:
\begin{equation}\label{99p}
    dq/dt=Tq(t),\,\,q(0)=q_{0},
\end{equation}
where $T$ is a given a closed unbounded quaternionic operator and $q(t)$ is a quaternionic
function defined for $t\geqslant0$ such that $q(t)$ belongs to the domain of
$T$, $q_0$ is a given initial datum.
In \cite{kexi-wq10}, the theory of quaternionic exponentially dichotomous operators
was developed via $S$-spectral splitting and applied to the Cauchy problem.
Furthermore, in \cite{kexi-xty}, the ergodic theory of continuous and discrete quaternionic
semigroups, as well as Tauberian theory, were developed within the quaternionic setting,
together with the ergodicity of solutions to the Cauchy problem.
However, the problem \eqref{99p} triggers an extremely important Cauchy problem for a
quaternionic evolution equations with a distinguished direction of time flow:
\begin{equation}\label{99p1}
    dq/dt=T(t)q(t),\,\,q(0)=q_{0},
\end{equation}
where $T(t)$ is a given a closed unbounded quaternionic operator which
is dependent on time $t$. More generally, one may directly embed \eqref{99p1} into
a quaternionic space with multi-dimensional evolution directions, thereby elevating
it to a Cauchy problem dependent on quaternionic variables:
\begin{equation}\label{99p2}
    \partial_{S}^{R} u (q) = A(q)u(q),\,\,
    u(s_{0})=u_{0},\qquad s_{0},q\in \Sigma.
\end{equation}
where $\partial_{S}^{R}$ is a right slice derivative operator, $A(q)$ is a given
quaternionic operator which is dependent on $q$ and $\Sigma$ is a spherical sector
in quaternionic space. Although quaternionic variables formally include a real
component that may be interpreted as time, the quaternionic Cauchy problem does not
intrinsically encode a temporal evolution structure. Instead, it provides a higher-dimensional
and noncommutative analytical framework within which classical time-dependent Cauchy problems
can be embedded once a distinguished real direction (or slice) is specified.
Therefore, rather than directly generalizing time-dependent well-posedness,
the quaternionic setting offers a broader structural environment
that requires additional geometric or algebraic constraints to recover the notion of
dynamical evolution.

From this perspective, the analysis of Cauchy problems inherently relies on the
identification of an appropriate evolution mechanism that governs the propagation of states.
In classical real or complex Banach space, such mechanisms are typically realized through
strongly continuous semigroups or, more generally, evolution families, which furnish a
unified framework for establishing existence, uniqueness, and stability results.
However, in the quaternionic context, the noncommutative algebraic structure prevents a
direct transfer of these classical theories, as the underlying assumptions of linear
ordering and commutativity no longer hold. Consequently, the construction of a suitable
notion of  quaternionic evolution system is not merely a technical extension,
but rather a fundamental prerequisite for formulating and analyzing
Cauchy problems in this setting. Due to the complex geometric structure of the
quaternionic space, these unknown challenges give rise to a series of core questions:
what conditions guarantee the existence and uniqueness of a local $\sigma$-analytic solution
to the evolution equation \eqref{99p2}? In what manner can the standard notion
of a quaternionic evolution system be defined so that such a system exists and is unique?
Is the local $\sigma$-analytic solution of the quaternionic Cauchy problem \eqref{99p2}
uniquely extendable from a closed $\sigma$-ball to the entire spherical sectorial region?
How do hyperbolic-like and parabolic-like constraints independently ensure the existence
and uniqueness of such a global analytic solution? Furthermore, after global extension,
is the $\sigma$-analyticity of the solution to the inhomogeneous evolution equation
(via the standard notion of a quaternionic evolution system)
\begin{equation*}
    \left\{\begin{aligned}
        &\partial_{S}^{R} u(q)+A(q)u(q)=f(q), \\
        &u(s_{0})=u_{0},
    \end{aligned}\right. \quad s_{0},q\in\Sigma
\end{equation*}
preserved, or does any nontrivial evolution inevitably destroy the cross-slice
analytic structure, thereby restricting the evolution to a single isolated slice?
Finally,  can  the entire well-posedness theory be lifted from $\mathbb{H}$-algebras
to $\mathbb{A}$-algebras, thereby demonstrating the universality of the well-posedness
theory for evolution equations thus established? All these successive questions will
be answered in this paper.

In Section \ref{sec2}, we establish the geometric properties of slice structures in
Banach bilateral $\mathbb{H}$-modules and construct noncommutative Banach bilateral
$\mathbb{H}$-subalgebras compatible with the slice regularity of the functions involved.
Section \ref{sec3} is devoted to establishing essential properties of slice regular
functions valued in such modules and of the slice structure of Banach bilateral
$\mathbb{H}$-algebras. In Section \ref{sec4}, a quaternionic sliced geometric path
integral over a Banach bilateral $\mathbb{A}$-module is employed to obtain, for the first time,
an integral representation of the solution to the homogeneous Cauchy problem within a
quaternionic $\sigma$-ball, thereby proving a local existence and uniqueness theorem
for $\sigma$-analytic solutions; this local solution is then extended to the entire
spherical sectorial region, yielding a global existence and uniqueness theorem,
from which the notion of an evolution system in a quaternionic $\mathbb{H}$-algebra
is introduced. Moreover, global existence and uniqueness is also established under
hyperbolic-like and parabolic-like conditions. Section \ref{sec5} is devoted to
lifting the well-posedness theory for evolution equations to  Banach bilateral
$\mathbb{A}$-modules and $\mathbb{A}$-algebras and revealing the universality of
well-posedness for Cauchy problems under generic algebraic structures.
In Section \ref{sec6}, we make a conclusion to end the paper.


\section{The Structural Properties of Banach Bilateral  $\mathbb{H}$-Modules}
\label{sec2}

Quaternions, an extension of complex numbers, are used to represent rotations
in three-dimensional space. A quaternion $s$ is typically expressed as
$ s = \hat{s}_{0} + \hat{s}_{1} i  + \hat{s}_{2} j  + \hat{s}_{3} k$,
where $\hat{s}_{\ell} $ $(\ell=0,1,2,3)$ are real numbers,
and $  i ,  j ,  k  $ are the imaginary units,
satisfying
\[
    i ^{2} =  j ^{2} =  k ^{2}= -1,\quad
    i  j =  k =- j  i ,\quad
    j  k = i =- k  j ,\quad
    k  i = j =- i  k .
\]
The real and imaginary parts of the quaternion $ s $
are defined as
\[
    \mathrm{Re}(s):=\hat{s}_{0}, \quad\mathrm{Im}(s):= \hat{s}_{1} i  + \hat{s}_{2} j  +
    \hat{s}_{3} k ,
\]
while the conjugate is defined as
$ \bar{s} = \hat{s}_{0} - \hat{s}_{1} i  - \hat{s}_{2} j  - \hat{s}_{3} k  $.
Additionally, the modulus and inverse of $ s$ are given by
\[
    \vert s\vert=\sqrt{s\bar{s}}=\sqrt{\bar{s}s},\quad s^{-1}=\frac{\bar{s}}{\vert s\vert^{2}}.
\]

In the following discussion, we will use the symbols $\mathbb{R}, \mathbb{C}, \mathbb{H}$
to denote the sets of real numbers, complex numbers, and quaternions, respectively.
Let us denote by $ \mathbb{S} $ the unit sphere consisting of purely imaginary quaternions,
i.e.,
\[
    \mathbb{S}=\{\hat{s}_{1}i+\hat{s}_{2}j+\hat{s}_{3}k\in\mathbb{H}:
    \hat{s}_{1}^{2}+\hat{s}_{2}^{2}+\hat{s}_{3}^{2}=1\}.
\]
It can be concluded that for any element $I$ in $\mathbb{S}$,
the property $I^{2}=-1$ holds true. Therefore,
for a given $ I \in \mathbb{S} $, we denote by $ \mathbb{C}_I $ the complex plane
$ \mathbb{R} + I\mathbb{R} $, which consists of elements expressed in the form
$ a + Ib $, where $ a, b \in \mathbb{R} $.
Moreover, for any $s\in\mathbb{H}$, the sphere $\mathbb{S}$ contains an element linked to it,
i.e.,
\begin{equation*}
     I_{s} :=
    \left\{\begin{aligned}
        &\frac{\mathrm{Im}(s)}{\vert \mathrm{Im}(s)\vert}&&\text{if }\mathrm{Im}(s)\neq 0,\\
        &\text{any element of }\mathbb{S}\quad &&\text{otherwise}.
    \end{aligned}\right.
\end{equation*}
For further details regarding the two-dimensional sphere
of purely imaginary unit quaternions
$\mathbb{S}$, one may refer to \cite{kexi-x8}.

The non-commutativity of quaternions yields vector space structures fundamentally different from those in the complex case, necessitating the consideration of distinct cases.

\begin{definition}[\cite{kexi-basp2}]\label{kexi-def2.1}
    A right $\mathbb{H}$-module is defined as an additive group $(X,+)$
    equipped with a right scalar multiplication
    $X\times \mathbb{H}\rightarrow X$, denoted by $(\mathbf{v},s)\mapsto \mathbf{v}s$, satisfying
    the following properties for all $\mathbf{u},\mathbf{v}\in X$
    and $q,s\in\mathbb{H}:$
    \begin{align*}
        (\mathbf{u}+\mathbf{v})s=\mathbf{u}s+\mathbf{v}s,\quad \mathbf{u}\cdot 1=\mathbf{u},\quad
        \mathbf{u}(s+q)= \mathbf{u}s+ \mathbf{u}q,\quad  ( \mathbf{u}q)s=\mathbf{u}(qs ).
    \end{align*}
    Similarly,
    a  left
    $\mathbb{H}$-module is an additive group $(X,+)$ endowed with a left
    scalar multiplication
    $\mathbb{H}\times X\rightarrow X$, denoted by $(s,\mathbf{v})\mapsto s\mathbf{v}$, satisfying
    the following properties for all
    $\mathbf{u},\mathbf{v}\in X$ and $q,s\in\mathbb{H}:$
    \begin{align*}
        s(\mathbf{u}+\mathbf{v})=s\mathbf{u}+s\mathbf{v},\quad 1\cdot\mathbf{u}=\mathbf{u},\quad
        (s+q)\mathbf{u}=s\mathbf{u} +q\mathbf{u} ,\quad s(q\mathbf{u})=(sq)\mathbf{u}.
    \end{align*}
    A bilateral $\mathbb{H}$-module is an additive group $(X,+)$ equipped with
    both a left and a
    right scalar multiplication, such that $(X,+)$ with the left scalar multiplication forms
    a left
    $\mathbb{H}$-module, and $(X,+)$ with the right scalar multiplication forms a
    right $\mathbb{H}$-module. Additionally, it satisfies the condition
    $s\mathbf{v}=\mathbf{v}s$ for all $\mathbf{v}\in X$ and $s\in\mathbb{R}$.
\end{definition}

The following result from Lemma 1.5 in \cite{kexi-basp0}
provides a crucial conclusion concerning bilateral $\mathbb{H}$-modules.

\begin{lemma}[\cite{kexi-basp0}]\label{kexi-basp1}
    Let $X$ be a bilateral $\mathbb{H}$-module and let
    $X_{\mathbb{R}} = \{\phi \in X : a\phi = \phi a, \forall a \in\mathbb{H} \}$.
    Then $X_{\mathbb{R}}$ is a real vector space, called the real part of $X$, and
    \begin{equation*}
        \Re(\varphi)=\frac{1}{4}(\varphi-i\varphi i- j \varphi j - k\varphi k)
        \quad  (\varphi\in X)
    \end{equation*}
    defines an $\mathbb{R}$-linear mapping from $X$ onto $X_{\mathbb{R}}$ that satisfies
    $\Re\circ\Re=\Re$. Moreover, any $\varphi\in X$ satisfies
    the polarization identity
    \begin{align*}
        \varphi
        =&~\Re(\varphi)+\Re(-i\varphi)i+\Re(-j\varphi)j +\Re(-k\varphi)k\\
        =&~\Re(\varphi)+i\Re(-i\varphi)+j\Re(-j\varphi) +k\Re(-k\varphi).
    \end{align*}
\end{lemma}

\begin{remark}\label{kexi-rem2}
    It should be noted that  if $X$ is a nontrivial bilateral $\mathbb{H}$-module,
    then the real vector space $X_{\mathbb{R}}$ is also nontrivial. In fact,
    for any $\mathbf{0}\neq\varphi\in X$,
    one has $\mathbf{0}\neq\Re(\varphi)\in X$ and
    \begin{align*}
        i \Re(\varphi)
       =&~ i
       \frac{1}{4}(\varphi- i \varphi  i
       - j \varphi  j - k \varphi  k )
        =\frac{1}{4}( i \varphi+\varphi  i
        - k \varphi  j  + j \varphi  k )\\
        =&~\frac{1}{4}(- i \varphi i +\varphi
        - k \varphi  k  - j \varphi  j ) i
        =\Re(\varphi) i \in X,
    \end{align*}
    the same result can be derived for $ j $ and $ k $,
    which implies that for any $a\in \mathbb{H}$ and $\varphi\in X$,
    $a\Re(\varphi)=\Re(\varphi)a$, i.e., $\mathbf{0}\neq \Re(\varphi)\in X_{\mathbb{R}}\neq \emptyset$.
    Therefore, without loss of generality, in the following discussion,
    we shall assume that $X$ is nontrivial.
\end{remark}

\begin{remark}\label{ree0}
For Lemma \ref{kexi-basp1}, by substituting the orthonormal basis
$\{ i , j , k \}$
with an arbitrary orthogonal basis selected from $\mathbb{S}$, say $\{I,J,K\}$,
we can derive the following extension
    \begin{equation*}
    \varphi = \Re' (\varphi) + \Re' (\overline{I}\varphi)I
    + \Re' (\overline{J}\varphi)J + \Re'(\overline{K}\varphi)K,
    \end{equation*}
    where $I, J, K$ are the quaternionic units satisfying $I^2 = J^2 = K^2 = IJK = -1$ and
        \begin{equation*}
        \Re'(\varphi)=\frac{1}{4}(\varphi-I\varphi I- J\varphi J- K\varphi K)
        \quad  (\varphi\in X).
    \end{equation*}
    Thus, the expression of $\varphi$ can be further rearranged as
    \begin{equation*}\label{exc0}
    \varphi = \left(\Re '(\varphi)
    + \Re' (\overline{I}\varphi)I\right)
    + \left(\Re'(\overline{J}\varphi)
    + \Re'(\overline{K}\varphi)I\right)J.
    \end{equation*}
    Since
    $\Re'(\varphi)\in X_{\mathbb{R}}$
    and $\overline{I}\varphi,\overline{J}\varphi,\overline{K}\varphi\in X$,
    then
    $\Re'(\overline{I}\varphi),\Re'(\overline{J}\varphi),\Re'(\overline{K}\varphi)\in X_{\mathbb{R}}$. Therefore,
    \begin{equation}\label{hh6}
    \varphi\in (X_{\mathbb{R}}+X_{\mathbb{R}}I)+(X_{\mathbb{R}}+X_{\mathbb{R}}I)J.
    \end{equation}
\end{remark}

\begin{example}
    Let
    \[
        X := C([0,1], \mathbb{H})
    =
        \{ f:[0,1] \to \mathbb{H} \mid f \text{ is continuous} \},
    \]
    equipped with the supremum norm
    \[
        \|f\|_\infty := \sup_{t\in[0,1]} |f(t)|.
    \]
    Then $(X, \|\cdot\|_\infty)$ is a natural Banach bilateral $\mathbb{R}$-module.

    For $a,b\in\mathbb{H}$, define the left and right scalar multiplications by
    \[
        (a \cdot f \cdot b)(t) := a\,f(t)\,b,\quad t\in[0,1],
    \]
    which makes $X$ a Banach bilateral $\mathbb{H}$-module.
    Now the real part of $X$ is
    \[
        X_{\mathbb{R}}
        := \{ f \in X : a f = f a,\ \forall a \in \mathbb{H} \}.
    \]
    Since a quaternion commutes with all elements of $\mathbb{H}$ if and only if it is real,
    it follows that
    $X_{\mathbb{R}} = C([0,1],\mathbb{R})$,
    the space of all real-valued continuous functions on $[0,1]$, naturally embedded in $X$.
\end{example}

According to Remark \ref{ree0}, for any $\varphi\in X$, it admits a splitting representation 
as given in \eqref{hh6}, which allows the bilateral $\mathbb{H}$-module to be defined in a 
split form as follows.

\begin{lemma}
    Let $X$ be a bilateral $\mathbb{H}$-module and
    \[
        X_{I}:=\{x+yI:x,y \in X_{\mathbb{R}}\},\quad I\in \mathbb{S}.
    \]
    Then $X=X_{I}+X_{I}J$, where $I,J\in\mathbb{S}$ and  $I\perp J$.
\end{lemma}

Although we can construct a quaternionic space isomorphic to $\mathbb{C}_{I}$
via a bilateral $\mathbb{H}$-module, it is essentially different from the split form of 
the bilateral $\mathbb{H}$-module. In the following, 
we will use the concept of a positive cone (see \cite{kexi-cone}) to construct a 
quaternionic space isomorphic to $\mathbb{C}_{I}$.

\begin{definition}
    Let $X$ be a Banach bilateral $\mathbb{H}$-module,
    and let $X_{\mathbb{R}}$ be the real part of $X$.
    A subset $P \subset X_{\mathbb{R}}$ is called a
    positive cone of $X$ if
    it satisfies the following conditions:
    \begin{itemize}
        \item[\rm(i)] For all $x,y \in P$, then $x+y\in P$.
        \item[\rm(ii)] For all $x \in P$ and all real numbers
        $t \geq 0$, then $t x \in P$.
        \item[\rm(iii)] If $x \in P$ and $-x \in P$, then $x = 0$.
         \item[\rm(iv)] $P\neq \{0\}$.
    \end{itemize}
\end{definition}

\begin{remark}
    Note that for a Banach bilateral $\mathbb{H}$-module $X$, by Remark \ref{kexi-rem2},
    $X_{\mathbb{R}}\neq \emptyset$.
    If $\dim (X_{\mathbb{R}}) \geq 1$,
    then there exists a nontrivial positive cone $P \subset X_{\mathbb{R}}$.
    Indeed, choose a nonzero linear functional $f \in X_{\mathbb{R}}^{*}$.
    Define
    \[
    P := \{\, x \in X_{\mathbb{R}} : f(x) \ge 0 \,\}.
    \]
    Then $P$ is a
    positive cone of $X$, since $P+P \subset P$ and $\lambda P \subset P$
    for all $\lambda \in \mathbb{R}^{+}$.
    Moreover, $P \cap (-P) = \ker f=\{0\}$.
\end{remark}

Based on the concept of a positive cone,
we introduce the following notion for equivalence classes
to construct a quaternionic Banach space which will be isomorphic to
$\mathbb{C}_{I}$.

\begin{figure}[htbp]
        \centering{
        \includegraphics[scale=0.6]{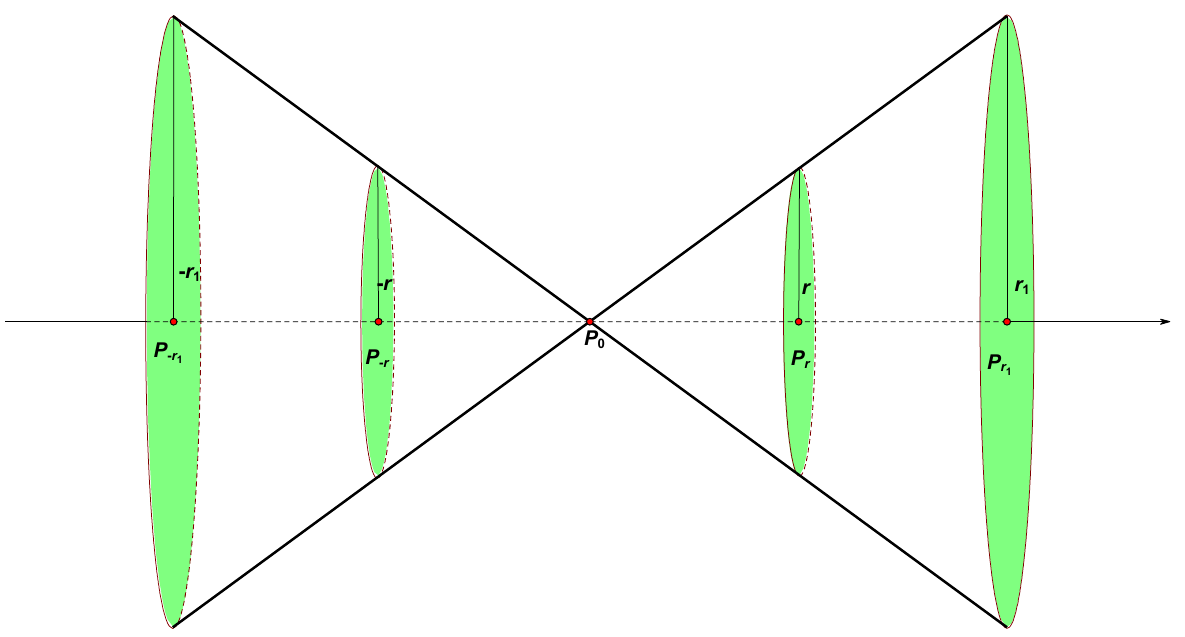}}
        \caption{A cone $\widetilde{P}$ and the equivalence class $P_{r}$}
        \label{ke-xi3}
\end{figure}

\begin{definition}\label{hh7}
Let $(X,\|\cdot\|_{X})$ be a Banach bilateral $\mathbb{H}$-module, $X_{\mathbb{R}}$
be the real part of $X$
and let $P\subset X_{\mathbb{R}}$ be a positive cone of $X$.
Define a relation $\overset{r}{\sim}$ on $X_{\mathbb{R}}$ as follows:
\begin{itemize}
    \item[\rm(i)] If $x,y\in P$ and $r\geq 0$,
    then
    \[
    x \overset{r}{\sim} y \quad \Longleftrightarrow \quad \|x\|_{X}=\|y\|_{X}=r.
    \]
    \item[\rm(ii)] If $x,y\in - P$ and $ r< 0$,
    then
    \[
        x \overset{r}{\sim} y \quad \Longleftrightarrow \quad \|x\|_{X}=\|y\|_{X}=-r.
    \]
\end{itemize}
Then $\overset{r}{\sim}$ is an equivalence relation on $\widetilde{P}:=P\cup(-P)$.
For each $x\in \widetilde{P}$,
its equivalence class $P_{r}$ $($see Figure \ref{ke-xi3}$)$ is defined by
\[
    P_{r}:= \{ y\in \widetilde{P} : y\overset{r}{\sim} x  \},\quad r\in\mathbb{R} .
\]
The quotient set induced by the equivalence relation $\overset{r}{\sim}$ is defined as
\[
    \widetilde{P}/\overset{r}{\sim}~\cong \mathbb{R}.
\]
\end{definition}

Based on Definition \ref{hh7}, the following quaternionic quotient space can be constructed,
and the two-dimensional imaginary sphere of the quotient space can be defined.
\begin{definition}
    Let $X$ be a Banach bilateral $\mathbb{H}$-module.
    On the quotient set $\widetilde{P}/ \overset{r}{\sim}$, define a norm
    \[
    \|\widetilde{r}\|_{\widetilde{P}/ \overset{r}{\sim}} := |r|,
    \qquad \widetilde{r}\in \widetilde{P}/ \overset{r}{\sim},
    r \in \mathbb{R}
    \]
    and define
    \[
        \mathbb{H}_{\widetilde{P}/ \overset{r}{\sim}}:=
        \big\{\widetilde{r}_{0}+\widetilde{r}_{1}i+\widetilde{r}_{2}j+\widetilde{r}_{3}k:
        \widetilde{r}_{\ell}\in \widetilde{P}/ \overset{r}{\sim},\ell=1,2,3\big\}.
    \]
    Furthermore, define
    \[
        \mathbb{S}_{\widetilde{P}/ \overset{r}{\sim}}:=
        \big\{\widetilde{r}_{1}i+\widetilde{r}_{2}j+\widetilde{r}_{3}k:
        \|\widetilde{r}_{1}\|^{2}_{\widetilde{P}/ \overset{r}{\sim}}
        +\|\widetilde{r}_{2}\|^{2}_{\widetilde{P}/ \overset{r}{\sim}}
        +\|\widetilde{r}_{3}\|^{2}_{\widetilde{P}/ \overset{r}{\sim}}=1,
        \widetilde{r}_{\ell}\in \widetilde{P}/ \overset{r}{\sim},\ell=1,2,3\big\}
    \]
    and
    \[
        \mathbb{C}_{\mathfrak{I}}:=\{\widetilde{r}_{0}+\eta\mathfrak{I}:\widetilde{r}_{0}\in \widetilde{P}/ \overset{r}{\sim},
        \eta\in\mathbb{R},\mathfrak{I}\in\mathbb{S}_{\widetilde{P}/ \overset{r}{\sim}}\}.
    \]
\end{definition}

\begin{remark}
    Note that an isomorphism $\psi$ is defined between
    $\mathbb{H}_{\widetilde{P}/ \overset{r}{\sim}}$ and
    $\mathbb{H}$ as follows:
\begin{align*}
\psi :\mathbb{H}_{\widetilde{P}/ \overset{r}{\sim}}\longrightarrow \mathbb{H}  ,
\qquad \widetilde{r}_{0} +\widetilde{r}_{1}i+\widetilde{r}_{2}j+\widetilde{r}_{3}k
\longmapsto r_{0}+r_{1}i+r_{2}j+r_{3}k.
\end{align*}
Thus, for any $I\in\mathbb{S}$, there exists
$\mathfrak{I}\in \mathbb{S}_{\widetilde{P}/ \overset{r}{\sim}}$
such that $\mathbb{C}_{I}\cong \mathbb{C}_{\mathfrak{I}}$.
\end{remark}

Next, we further introduce the concept of the ``unit sphere'' in $X$.
    Let $X$ be a Banach bilateral $\mathbb{H}$-module, we define
    \[
        \mathbb{S}_{X}:=\big\{x_{1} i +x_{2} j +x_{3} k \in X:
        \Vert x_{1}\Vert^{2}_{X}+\Vert x_{2}\Vert^{2}_{X}+\Vert x_{3}\Vert^{2}_{X}=1,x_{l}\in
        X_{\mathbb{R}} ,l=1,2,3\big\}.
    \]
Furthermore,
    \[
        X_{\aleph }:=\{x+\alpha\aleph:x \in X_{\mathbb{R}},\alpha\in\mathbb{R}\},
        \quad  \aleph\in \mathbb{S}_{X}.
    \]

For any element
    $x = x_{0} + x_{1} i  + x_{2} j  + x_{3} k  \in X$
    with a nonzero vector part
    $x_{1} i  + x_{2} j  + x_{3} k  \neq 0$,
    we obtain that the unitization of ``imaginary part'' of $x$ as follows:
    \[
        \mathbb{S}_{X} \ni \aleph = \frac{x_{1} i  + x_{2} j
        + x_{3} k }{\sqrt{\Vert x_{1}\Vert^{2}_{X} + \Vert x_{2}\Vert^{2} _{X}
        +\Vert x_{3}\Vert^{2}_{X}}}.
    \]

\begin{proposition}
If $\mathrm{dim}(X_{\mathbb{R}})=n$, then $\mathbb{S}_{X}\cong \mathbb{S}^{n}$.
\end{proposition}
\begin{proof}
Let $\mathcal{V} = \{e_{1}, e_{2}, \cdots, e_{n}\}$ be a basis of $X_{\mathbb{R}}$.
Then for any vector $x \in X_{\mathbb{R}}$
such that
$x = \sum_{\ell=1}^{n} \xi_{\ell} e_{\ell}$, where $\xi_{1}, \xi_{2}, \cdots, \xi_{n} \in \mathbb{R}$,
we define
a mapping $\widetilde{\phi}: X_{\mathbb{R}} \to \mathbb{R}^{n}$, i.e., $\widetilde{\phi}(x) = (\xi_{1}, \xi_{2}, \cdots, \xi_{n})^T$.

Now, for any $u, v \in X_{\mathbb{R}}$ and $c \in \mathbb{R}$,
let $u = \sum_{\ell=1}^{n} a_{\ell} e_{\ell}$ and $v = \sum_{\ell=1}^{n} b_{\ell} e_{\ell}$.
Then $cu + v = \sum_{\ell=1}^{n} (ca_{\ell} + b_{\ell}) e_{\ell}$.
Thus,
\[
\widetilde{\phi}(cu + v) = (ca_1+b_1, \dots, ca_n+b_n)^T
= c\widetilde{\phi}(u) + \widetilde{\phi}(v),
\]
which implies that $\widetilde{\phi}$ is a linear mapping. On the other hand,
for $w=\sum_{\ell=1}^{n}\eta_{\ell}e_{\ell}$,
if $\widetilde{\phi}(w) =0$,
since $e_{\ell}\neq 0$ $(\ell=1,\cdots,n)$,
then $\eta_{1} = \eta_{2} = \cdots = \eta_{n} = 0$. Thus,
 $w = \sum_{\ell=1}^{n} 0 \cdot e_{\ell} = 0$, and hence
$\ker(\widetilde{\phi}) = \{0\}$, which implies that $\widetilde{\phi}$ is injective.
Furthermore,
for any $m = (m_1, \cdots, m_n)^{T} \in \mathbb{R}^{n}$, there exists
vector $\widetilde{m} = \sum_{\ell=1}^{n} m_{\ell} e_{\ell}=\widetilde{\phi}^{-1}(m)\in X_{\mathbb{R}}$.
Therefore, $\widetilde{\phi}$ is a surjective. Then
$\widetilde{\phi}$ is a linear isomorphism,
i.e., $X_{\mathbb{R}} \cong \mathbb{R}^{n}$. Thus,
    \[
        \mathbb{S}_{X} \subseteq X_{\mathbb{R}}i\times X_{\mathbb{R}}j
        \times X_{\mathbb{R}}k \cong \mathbb{R}^{3n}.
    \]
    Let $\mathbb{S}_{\mathbb{R}^{3}}=\{y=(y_{1},y_{2},y_{3})^{T}\in\mathbb{R}^{3}:y_{1}^{2}
    +y_{2}^{2}+y_{3}^{2}=1\}$.
    Since $\mathbb{S}\cong \mathbb{S}_{\mathbb{R}^{3}}$,
     then
    \[
        \mathbb{S}_{X}\cong \underbrace{\mathbb{S}\times\mathbb{S}\times\cdots\times\mathbb{S}}_{n}
        =:\mathbb{S}^{n}.
    \]
    This completes the proof.
\end{proof}

In what follows, we provide some necessary knowledge of quaternionic operators.
\begin{definition}[\cite{kexi-x8}]\label{kexi-su0}
    Let $X$ be a bilateral $\mathbb{H}$-module. 
    A map $T:X\rightarrow X$
    is called a right linear operator if for all $\lambda\in\mathbb{H}$ and $u,v\in X$,
    the following relationships hold
    \[
        T(u+v)=T(u)+T(v),\quad T(u\lambda) =T(u)\lambda.
    \]
    Note that the set consisting of all quaternionic right linear operators with respect to the
    following operations
    \[
        (\lambda T)(v):=\lambda T(v) ,\quad (T\lambda)(v):= T(\lambda v),\quad\text{ for all }
        \lambda\in\mathbb{H},\text{ and for all } v\in X,
    \]
    is both a left $\mathbb{H}$-module and a right $\mathbb{H}$-module. Similarly, 
    A map $T:X\rightarrow X$
    is called a left linear operator if for any $\lambda\in\mathbb{H}$ and $u,v\in X$, the following
    relationships hold
    \begin{align*}
        T(u+v)=T(u) +T(v), \quad T(\lambda u)=\lambda T(u).
    \end{align*}
    Moreover, the set consisting of all quaternionic left linear operators with respect to the
    following operations
    \begin{align*}
        (T\lambda)(u):=T(u)\lambda,\quad (\lambda T)(u):=T(u\lambda),\quad\text{ for all }
        \lambda\in\mathbb{H},\text{ and for all } u\in X,
    \end{align*}
    is both a left $\mathbb{H}$-module and a right $\mathbb{H}$-module.
\end{definition}

    Let $X$ and $Y$ be two bilateral $\mathbb{H}$-modules.
    A mapping $T:X\to Y$ is said to be bounded if for every $x_{1},x_{2}\in X$ there exists constant
    $M$ such that
    \[
        \Vert T(x_{1})-T(x_{2})\Vert_{Y}\leq M\Vert x_{1}-x_{2}\Vert_{X}.
    \]
    The bilateral $\mathbb{H}$-module
    of all bounded linear quaternionic operators from $X$ to $Y$ is denoted
    by $\mathcal{B}(X,Y)$.
    The norm of $T\in \mathcal{B}(X, Y)$ is defined in usual way:
    \[
       \Vert  T\Vert _{\mathcal{B}(X,Y)}\overset{def}{=}\sup\limits_{0\ne x\in X}
       \frac{\Vert  Tx\Vert _{Y}}{\Vert  x\Vert _{X}}=\sup\limits_{\Vert  x\Vert _{X}=1}
       \Vert  Tx\Vert _{Y}.
    \]
    If $X=Y$, then we write $\mathcal{B}(X)$ instead of $\mathcal{B}(X, Y)$, and we write
    $X \oplus  Y$ for the algebraic
    direct sum of the $\mathbb{H}$-modules $X$ and $Y.$

For convenience of our discussion, we will use the following notations related to the
bilateral $\mathbb{H}$-module
$X$ (see Colombo et al. \cite{kexi-x8}, pp. 136-138).

\begin{itemize}
    \item[\rm (i)] $\mathcal{B}^{R}(X)$ is the bilateral $\mathbb{H}$-module
    of all right linear
    bounded quaternionic operators on $X$;
    \item[\rm (ii)] $\mathcal{B}^{L}(X)$ is the bilateral $\mathbb{H}$-module
    of all left linear
    bounded quaternionic operators on $X$;
    \item[\rm (iii)] when it is not necessary to differentiate between left or right linear
    bounded quaternionic operators on $X$, we use the symbol $\mathcal{B}(X)$ and we call an element in
    $\mathcal{B}(X)$ a ``linear operator''.
\end{itemize}

\begin{remark}
    Note that $\mathbb{H}$ is a finite dimensional associative nontrivial real
    $\ast$-algebra with unite.
    Therefore, based on Section 5.3 of \cite{kexi-srs0} p. 5012, the sum of operators, the scalar
multiplications, the composition and the operator norm make
    $\mathcal{B}^{L}(X)$ and $\mathcal{B}^{R}(X)$
    the Banach bilateral $\mathbb{H}$-algebras, respectively, if $X$ is
     a Banach bilateral $\mathbb{H}$-module.
\end{remark}

The pivotal concept of the $S$-spectrum, essential in quaternionic operator theory,
will now be introduced for our later discussion.

\begin{definition}[The \textit{S}-spectrum and the \textit{S}-resolvent sets for quaternionic
    operators, \cite{kexi-co1}]
    Let $T\in\mathcal{B}(X)$. Define
     the \textit{S}-spectrum $\sigma_{S}(T)$ of $T$ as$:$
    \[
        \sigma_{S}(T):=\{\lambda\in\mathbb{H}:T^{2}-2\mathrm{Re}(\lambda)T+\vert \lambda\vert^{2}\mathcal{I}
        \text{ is not invertible in }\mathcal{B}(X)\}.
    \]
    The \textit{S}-resolvent set $\rho_{S}(T)$ is defined by
    $\rho_{S}(T):=\mathbb{H}\backslash\sigma_{S}(T)$.
\end{definition}

\begin{definition}[The \textit{S}-resolvent operator, \cite{kexi-co1}]
    Let $T\in\mathcal{B}(X)$. Define the left \textit{S}-resolvent operator as
 	\begin{equation*}
 		S^{-1}_{L}(\lambda,T):=-\big(T^2-2\mathrm{Re} (\lambda)T+\vert\lambda
 		\vert^2\mathcal{I}\big)^{-1}(T-\overline{\lambda}\mathcal{I}),
 	\end{equation*}
 	and the right $S$-resolvent operator as
 	\begin{equation*}
 		S^{-1}_{R}(\lambda,T):=-(T-\overline{\lambda}\mathcal{I})(T^2-2\mathrm{Re}(\lambda)T+\vert\lambda
 		\vert  ^2\mathcal{I})^{-1}.
 	\end{equation*}
\end{definition}

\begin{definition}[\cite{kexi-x8}]
        Let $X$ be a Banach bilateral $\mathbb{H}$-module.
        \begin{description}
            \item[{\rm (i)}] The notation $\mathcal{K}^{R}(X)$
            will be used to denote the set of right linear closed operators
        $T:\mathcal{D}(T)\subset X\rightarrow X$, such that
        \begin{description}
            \item[\rm(1)] $\mathcal{D}(T)$ is dense in $X$,
            \item[\rm(2)]  $\mathcal{D}(T^2)\subset\mathcal{D} (T)$ is dense in $X$,
            \item[\rm(3)] $T-\bar{s}\mathcal{I}$ is densely defined in $X$.
        \end{description}
        \item[{\rm (ii)}] The notation $\mathcal{K}^{L}(X )$ 
        will be used to denote the set of left linear closed operators satisfying $(1)$ and
        $(2)$ and such that $T-\mathcal{I}\bar{s}$ is densely defined in $X$.
        \item[{\rm (iii)}] When no distinction is made between $\mathcal{K}^{L}(X)$ 
        and $\mathcal{K}^{R}(X)$, the symbol $\mathcal{K}(X)$ is used directly.
    \end{description}
        Note that $T^2-2 \mathrm{Re}(s)T+\vert s\vert^2I:\mathcal{D}(T^2)\subset X\rightarrow X$ is a closed operator. If
        $T\in\mathcal{K}(X)$, 
        denote by $\rho_{S}(T)$ the $S$-resolvent set of $T$ as
        \[
            \rho_{S}(T)=\big\{s\in\mathbb{H}:
            T^2-2 \mathrm{Re}(s)T+\vert s\vert^2\mathcal{I}
            \text{ is invertible in }\mathcal{B}(X)\big\},
        \]
        and $\sigma_{S}(T)=\mathbb{H}\backslash \rho_{S}(T)$. 
        Now, define $Q_s(T)$ as follows:
        \[
            Q_s(T):=(T^2-2 \mathrm{Re} (s)T+\vert s\vert^2\mathcal{I})^{-1}\quad where ~~~Q_s(T):X\rightarrow \mathcal{D}(T^2).
        \]
\end{definition}

By leveraging the structural properties of Banach bilateral $\mathbb{H}$-module and
$\mathbb{H}$-algebra, together with the
$S$-spectrum theory and
quaternionic operator theory, we will establish a theory of evolution dynamics with
quaternionic variables. Moreover, the Banach bilateral $\mathbb{H}$-algebra will be
lifted to an arbitrary associative real $*$-algebra, enabling us to propose a structural
evolution equation within the framework of a real $*$-algebra and develop its corresponding
dynamical theory.

\section{Slice Regularity in Banach Bilateral \text{$\mathbb{H}$}-Modules and \text{$\mathbb{H}$}-Algebras}\label{sec3}

This section first presents the definition and related properties of slice regular functions in Banach bilateral $\mathbb{H}$-module spaces. Subsequently, the concept of slice Banach left and right $\mathbb{H}$-subalgebras is introduced, their important properties are derived, and the splitting form of operator-valued functions taking values in left and right $\mathbb{H}$-subalgebras is obtained.

In addition, it is proved that, in the quaternionic $\sigma$-ball domain, the regularity of operator-valued functions in $\mathbb{H}$-subalgebras implies the regularity of functions in the $\mathbb{H}$-module space. These fundamental results lay an important foundation for the subsequent study of evolution dynamics with quaternionic variables and the extension of Banach $\mathbb{A}$-algebras.

\subsection{Slice Regular Functions Valued in Banach Bilateral $\mathbb{H}$-Modules}

In quaternionic spaces, the definition of ``holomorphic'' functions in the traditional sense is invalid, or directly leads to the linearity of functions. Therefore, slice derivatives and their geometric transformations have become one of the important tools for studying operators and spectral theory in quaternionic spaces.

\begin{proposition}[\cite{kexi-xdao4}]\label{zsh7}
    Let $\Omega$ be a simply-connected domain in $\mathbb{H}$ and let $f:\Omega\rightarrow \mathbb{H}$.
    If for any $q_{0}\in\Omega$
    \[
        \lim_{q\rightarrow q_{0}}(q-q_{0})^{-1}\bigl(f(q)-f(q_{0})\bigr)
    \]
    exists in $\mathbb{H}$, then necessarily $f(q)=qa+b$ for some $a,b\in\mathbb{H}$. If for any
    $q_{0}\in \Omega$
    \[
        \lim_{q\rightarrow q_{0}}\bigl(f(q)-f(q_{0})\bigr)(q-q_{0})^{-1}
    \]
    exists in $\mathbb{H}$, then necessarily $f(q)=aq+b$ for some $a,b\in\mathbb{H}$.
\end{proposition}

As shown by Proposition \ref{zsh7}, adhering strictly to classical complex differentiability limits quaternionic
differentiable functions to affine functions (linear polynomials). Therefore, the following quaternionic slice derivative is needed.

\begin{definition}[\cite{kexi-srs2,Al5,kexi-co3}]
    Let $f : U \subseteq \mathbb{H} \to \mathbb{H}$ and let $q = u + Iv \in U$ $(I\in\mathbb{S})$.
    If $q$ is not real, then $f$ is said to admit a
    left slice derivative in $q$ if
    \begin{equation}\label{e0}
        \partial_S f(q) := \lim_{p \to q,\, p \in \mathbb{C}_{I}} (p - q)^{-1} (f_{I}(p) - f_{I}(q))
    \end{equation}
    exists and is finite. If $q$ is real, then $f$ is said to admit a
    left slice derivative in $q$ if \eqref{e0} exists for every $I \in \mathbb{S}$.

    Similarly, $f$ is said to admit a right slice derivative at a
    nonreal point $q = u + Iv \in U$ if
    \begin{equation}\label{e1}
        \partial_S f(q) := \lim_{p \to q,\, p \in \mathbb{C}_{I}} (f_{I}(p) - f_{I}(q))(p - q)^{-1}
    \end{equation}
    exists and is finite, and $f$ is said to admit a right slice derivative at
    a real point $q \in U$ if \eqref{e1} exists and is finite for every $I \in \mathbb{S}$.
\end{definition}

\begin{remark}
    Note that in order to
    distinguish between the left and right slice derivative cases, we clarify the following notations.
    \begin{align*}
        &\partial_{S}^{L} f(q)
        := \lim_{p \to q,\, p \in \mathbb{C}_{I}} (p - q)^{-1} (f_{I}(p) - f_{I}(q)),
        \\
        &\partial_{S}^{R} f(q)
        := \lim_{p \to q,\, p \in \mathbb{C}_{I}}  (f_{I}(p) - f_{I}(q))(p - q)^{-1}.
    \end{align*}
   When the function $f: \mathbb{H}^{n} \to X$
    $($where $ n\in \mathbb{Z}^{+} $$)$
    is a multivariate function,
    then
    \[
        \frac{\partial^{L}_{S}}{\partial q_{\ell}} f(q_{1}, q_{2}, \cdots, q_{n}) \quad
        \left(\text{resp. } \frac{\partial^{R}_{S}}{\partial q_{\ell}}
        f(q_{1}, q_{2}, \cdots, q_{n})\right)
    \]
    denotes taking the left $($resp. right$)$ slice derivative of $f$
    with respect to its $\ell$-th variable $q_{\ell}$, where
     $\ell = 1, 2, \cdots, n$.
\end{remark}

In the following, we will present the definition of quaternionic Banach-valued slice regular functions,
and for this purpose, we set forth some necessary concepts.

\begin{definition}[\cite{kexi-co1}]\label{spo1}
    Let $U$ be a domain in $\mathbb{H}$. Then $U$ is called a
    slice domain $($$s$-domain for short$)$ if $U \cap \mathbb{R}$ is nonempty
    and if $U \cap \mathbb{C}_J$ is a domain in $\mathbb{C}_J$ for all $J \in \mathbb{S}$.
    $U$ is called axially symmetric if, for all $q = a + Ib \in U$ $(a,b\in\mathbb{R},I\in\mathbb{S})$,
    the 2-sphere $a + Jb$, $J \in \mathbb{S}$, is contained in $U$.
\end{definition}

\begin{definition}[\cite{kexi-co1}]
    Given $p\in\mathbb{H}$, $p=a+I_{p}b$ $(a,b\in\mathbb{R},I_{p}\in\mathbb{S})$, 
    denote by $[p]$ the set of all
    elements of the form $a+Jb$ when $J$ varies in $\mathbb{S}$.
\end{definition}

\begin{definition}[\cite{kexi-srs0}]\label{kexi-def3.3}
    For any $ I \in \mathbb{S} $, define the isomorphism
    $\phi_{I} : \mathbb{C} \longrightarrow \mathbb{C}_{I}$ as follows:
    \[
    \phi_{I}(a + b i ) := a + bI, \quad a, b \in \mathbb{R}.
    \]
    Let $\mathscr{D}$ be a subset of $\mathbb{C}$,
    symmetric with respect to complex conjugation. Then,
    a subset $\Omega_{\mathscr{D}}$ of $\mathbb{H}$ is said to be
    slice axially symmetric domain associated to $\mathscr{D}$
    if and only if
    \[
        \Omega_{\mathscr{D}} := \bigcup_{I\in \mathbb{S}} \phi_{I}(\mathscr{D})
        = \{a + bI \in \mathbb{H} : a, b \in \mathbb{R}, a + b i  \in \mathscr{D},
        I \in \mathbb{S}\}.
    \]
\end{definition}

Now,
we introduce the notion of
slice regular functions valued  in Banach bilateral $\mathbb{H}$-module.

\begin{definition}[\cite{kexi-srs2}]\label{ke-def5}
    Let $X$ be a Banach bilateral $\mathbb{H}$-module.
    Let $\Omega_{\mathscr{D}} \subseteq \mathbb{H}$ denote a
    slice axially symmetric domain associated to $\mathscr{D}$.
    A function $f : \Omega_{\mathscr{D}} \rightarrow  X$ is said to be a
    left slice function if it takes the form
    \[
        f\big(\phi_{I}(z)\big)
        = f_0(z) + I f_1(z) \quad \text{for } z = u +  v i  \in \mathscr{D}
    \]
    with two functions $f_0, f_1 : \mathscr{D} \rightarrow X$ that satisfy
    the compatibility condition
    \begin{equation}\label{xrtj0}
        f_0(z) = f_0(\overline{z}), \quad f_1(z) = -f_1(\overline{z}).
    \end{equation}
    In addition, if $f_0$ and $f_1$ satisfy the Cauchy-Riemann equations
    \begin{equation}\label{xrtj1}
        \frac{\partial}{\partial u}f_0  - \frac{\partial}{\partial v}f_1  = 0,\quad
        \frac{\partial}{\partial v}f_0  + \frac{\partial}{\partial u}f_1  = 0,
    \end{equation}
    then $f$ is called left slice regular. Similarly,
    a function $f : \Omega_{\mathscr{D}} \to X$ is
    said to be a right slice function if it takes the form
    \[
        f\big(\phi_{I}(z)\big) = f_0(z) + f_1(z)I \quad \text{for }
        z = u +  i v \in \mathscr{D}
    \]
    with two functions $f_0, f_1 : \mathscr{D} \rightarrow X$ that satisfy \eqref{xrtj0}.
    If in addition $f_0$ and $f_1$ satisfy the Cauchy-Riemann equations, then $f$ is called
    right slice regular.

The sets of left and right slice functions on
$\Omega_{\mathscr{D}}$ valued in
$X$
are denoted by $\mathcal{SF}_L(\Omega_{\mathscr{D}},X)$
and $\mathcal{SF}_{R}(\Omega_{\mathscr{D}},X)$, respectively, while
$\mathcal{SH}_L(\Omega_{\mathscr{D}},X)$ and
    $\mathcal{SH}_{R}(\Omega_{\mathscr{D}},X)$ denote the corresponding sets of left and right slice regular functions.
\end{definition}

The proofs of the following two results follow similarly to Lemma 4.1.7 in \cite{kexi-x8} and Theorem 2.19 in \cite{kexi-co3};
we omit them for brevity. For further properties of generalized slice regular functions, one may see \cite{kexi-srs0,kexi-srs2}.

\begin{proposition}[\cite{kexi-x8}]\label{kexi-fenli0}
    Let $\Omega_{\mathscr{D}}$ be a slice axially symmetric domain associated to $\mathscr{D}$
    and let $X$ be a Banach bilateral $\mathbb{H}$-module.
    Let
    $J, I \in \mathbb{S}$ such that $I \perp J$. Suppose that
    $f \in \mathcal{SH}_{L}(\Omega_{\mathscr{D}},X)$,
    then the restriction $f_{J}= f\vert_{\Omega_{\mathscr{D}}\cap \mathbb{C}_{J}}$ satisfies
    \[
        \frac{1}{2} \left( \frac{\partial}{\partial u} f_{J}(z)
        + J \frac{\partial}{\partial v} f_{J}(z) \right) = 0
    \]
    for all $z = u + Jv \in \Omega_{\mathscr{D}} \cap\mathbb{C}_{J}$.
    Furthermore, $F_{J}(z)$ admits the following representation:
    \[
        f_{J}(z) = F_1(z) + F_2(z)I
    \]
    with holomorphic functions
    $F_1,F_2:\Omega_{\mathscr{D}}\cap\mathbb{C}_{J}\rightarrow X_{J}$.

    Suppose that $f\in \mathcal{SH}_{R}(\Omega_{\mathscr{D}},X)$,
    then the restriction $f_{J} = f\vert_{\Omega_{\mathscr{D}} \cap \mathbb{C}_{J}}$ satisfies
    \[
        \frac{1}{2} \left( \frac{\partial}{\partial u} f_{J}(z) +
         \frac{\partial}{\partial v} f_{J}(z)J\right) = 0
    \]
    for all $z = u + Jv \in \Omega_{\mathscr{D}} \cap \mathbb{C}_{J}$. Moreover,
    \[
        f_{J}(z) = F_1(z) + IF_2(z)
    \]
    with holomorphic functions
    $F_1,F_2:\Omega_{\mathscr{D}} \cap \mathbb{C}_{J}\rightarrow X_{J}$.
\end{proposition}

\begin{proposition}[\cite{kexi-co3}]\label{mi0}
    Let $\Omega_{\mathscr{D}}$ be a slice axially symmetric domain
    associated to $\mathscr{D}$
    and let $X$ be a Banach bilateral $\mathbb{H}$-module.
    A function $f : \Omega_{\mathscr{D}} \to X$
    is a left slice function on $\Omega_{\mathscr{D}}$ if and only if for any
    $q = u + Jv \in \Omega_{\mathscr{D}}$, the following representation holds
    \[
        f(q) = \frac{1}{2} \big( f(\overline{z}) + f(z) \big)
        + \frac{1}{2} JI \big( f(\overline{z}) - f(z) \big),
    \]
    where $z = u + Iv$ and $I,J \in \mathbb{S}$.
    Similarly, a function $f : \Omega_{\mathscr{D}} \to X$
    is a right slice function on
    $\Omega_{\mathscr{D}}$ if and only if for any $q = u + Jv \in \Omega_{\mathscr{D}}$
    the following representation holds
    \[
        f(q) = \frac{1}{2} \big( f(\overline{z}) + f(z) \big) + \frac{1}{2}
        \big( f(\overline{z}) - f(z) \big) IJ,
    \]
    where $z = u + Iv$.
\end{proposition}

\subsection{Slice Structure of Banach Bilateral $\mathbb{H}$-Algebra and Slice Regular Semigroup}

In this subsection, using the concept of slice quaternionic operators, we will establish the slice regularity and product properties of quaternionic operator-valued functions in  the structure of Banach bilateral $\mathbb{H}$-modules.

\begin{definition}[\cite{kexi-wq10}]\label{sll1}
    Let $X$ be a Banach bilateral $\mathbb{H}$-module and
    $\mathbb{S}_{c}^{I}(X)=\big\{T\in\mathcal{K}(X): TI=IT\big\}$,
    $I\in \mathbb{S}$,
    a quaternionic operator $T$ is
    called slice quaternionic operator on $X$
    if $T\in\mathbb{C}_{I}^{\mathcal{S}}(X)$, where
    \[
    \mathbb{C}_{I}^{\mathcal{S}}(X) :
    =\big\{T\in\mathcal{K}(X):T=T_{0}+IT_{1},
    T_{l}\in\mathbb{S}_{c}^{I}(X),l=0,1\big\}, \quad I\in\mathbb{S}
    \]
    and $\mathbb{C}_{I}^{\mathcal{S}}(X)$ is called a slice operator space with respect
    to $I\in\mathbb{S}$.
    For convenience, $\mathbb{C}_{I}^{\mathcal{S}}(X)$ is denoted by $\mathbb{C}_{I}^{\mathcal{S}}$.
\end{definition}

\begin{remark}
    Note that for given $I\in\mathbb{S}$,
    $\mathbb{C}_{I}^{\mathcal{S}}(X)$ is isomorphic to some complex operator space, and
    for any $T \in \mathbb{C}_{I}^{\mathcal{S}}(X)$, $T$ is both
    left $\mathbb{C}_{I}$-linear and right $\mathbb{C}_{I}$-linear.
    Moreover,
    for $u=M+NJ\in X$ $($where $M,$ $N \in X_{I})$, we have
    \begin{align*}
        T(u)=&~(T_{0}+IT_{1})(M+NJ)\\
        =&~T_{0}(M+NJ)+(IT_{1})(M+NJ)\\
        =&~\big(T_{0}(M)+IT_{1}(M)\big)+\big(T_{0}(N)+IT_{1}(N)\big)J\\
        =&~T\vert_{X_{I}}(M)+T\vert_{X_{I}}(N)J.
    \end{align*}
    Since $T\in\mathbb{C}_{I}^{\mathcal{S}}(X)$, then
    \begin{align*}
    I(T\vert_{X_{I}}(M)+T\vert_{X_{I}}(N)J)
    =&~IT(u)=T(u)I\\
    =&~(T\vert_{X_{I}}(M)+T\vert_{X_{I}}(N)J)I,
    \end{align*}
    it follows that $T(M+NJ)=T\vert_{X_{I}}(M)$
    for any $M,N\in X_{I}$. Therefore, $T$ is uniquely
    determined by its restriction to $X_{I}$.
\end{remark}

\begin{proposition}
    Let $X$ be a Banach bilateral $\mathbb{H}$-module. Then
    $\mathbb{C}_{I}^{\mathcal{S}}(X)\cap\mathcal{B}^{R}(X)$ is a Banach bilateral
    $\mathbb{C}_{I}$-subalgebra of $\mathcal{B}^{R}(X)$.
\end{proposition}
\begin{proof}
    For any $T,S\in \big(\mathbb{C}_{I}^{\mathcal{S}}(X)\cap\mathcal{B}^{R}(X)\big)$,
    $q\in \mathbb{C}_{I}$ and $v\in X$, we have
    \begin{align*}
        T(v)+S(v)=&~(T_{0}+T_{1}I)(v)+(S_{0}+S_{1}I)(v)\\
        =&~ \big((T_{0}+S_{0})+(T_{1}+S_{1})I\big)(v)\in
        \big(\mathbb{C}_{I}^{\mathcal{S}}(X)\cap\mathcal{B}^{R}(X)\big)
    \end{align*}
    and
    \begin{align*}
        (qT)(v)=\big(qT_{0}+q (T_{1}I)\big)(v)=:(\widetilde{T}_{0}+\widetilde{T}_{1}I)(v)
        \in \big(\mathbb{C}_{L}^{\mathcal{S}}(X)\cap\mathcal{B}^{R}(X)\big),
    \end{align*}
    where $\widetilde{T}_{0}=qT_{0}$, $\widetilde{T}_{1}=qT_{1}$.
    Moreover,
    \begin{align*}
        (T\circ S)(v)=&~ (T_{0}+T_{1}I)\big(S_{0}v+(S_{1}I)v\big) \\
        =&~ T_{0}(S_{0}v) + T_{0}\big(S_{1}I(v)\big)
        + (T_{1}I)(S_{0}v)+ (T_{1}I)\big((S_{1}I)v\big)\\
        =&~(T_{0}S_{0}-T_{1}S_{1})(v)+ \big((T_{0}S_{1}+T_{1}S_{0})I\big)(v)
        \in \big(\mathbb{C}_{I}^{\mathcal{S}}(X)\cap\mathcal{B}^{R}(X)\big).
    \end{align*}
   Hence, $\mathbb{C}_{I}^{\mathcal{S}}(X)\cap\mathcal{B}^{R}(X)$ is a Banach bilateral
    $\mathbb{C}_{I}$-subalgebra of $\mathcal{B}^{R}(X)$.  This completes the proof.
\end{proof}

By means of Banach bilateral
    $\mathbb{C}_{I}$-subalgebra, we can construct some important types of operator spaces in quaternionic settings.

\begin{definition}\label{def0}
    Let $X$ be a Banach bilateral $\mathbb{H}$-module. Define
    \[
        \mathbb{C}^{\mathcal{S}}(X):=\bigcup_{I\in\mathbb{S}}\mathbb{C}_{I}^{\mathcal{S}}(X) \quad
    \text{and}
    \quad
        \mathscr{N}(X):=\big\{T\in\mathcal{K}(X):T\vert_{X_{I}}\in
        \mathbb{C}_{I}^{\mathcal{S}}(X)\big\}.
    \]
    Moreover, define the left and right slice operator spaces as follows:
    \begin{align*}
        &\mathbb{C}^{\mathcal{S}}_{L}(X):=\big\{T(v)=T_{0}(v)+(T_{1}J)(v):
        T_{0},T_{1}\in \mathscr{N} (X), J\in\mathbb{S} , v\in X \big\},  \\
        &\mathbb{C}^{\mathcal{S}}_{R}(X):=\big\{T(v)=T_{0}(v)+(JT_{1})(v):
        T_{0},T_{1}\in \mathscr{N} (X), J\in\mathbb{S}, v\in X \big\}.
    \end{align*}
    The elements in $\mathbb{C}^{\mathcal{S}}_{L}(X)$ and
    $\mathbb{C}^{\mathcal{S}}_{R}(X)$ are referred to as left slice
    operators and right slice operators, respectively.
\end{definition}

\begin{remark}
    Note that for any element
    $T=T_{0}+T_{1}J \in \mathbb{C}^{\mathcal{S}}_{L}(X)$,
    if its components $T_{0}$ and $T_{1}$ are left linear
    $($resp. right linear$)$ over $\mathbb{H}$, 
    then $T$ is also left linear $($resp. right linear$)$ over $\mathbb{H}$.
    For example, for any $S\in\mathbb{C}_{L}^{\mathcal{S}}(X)$ such that
    \[
        S=S_{1}+S_{2}J,
    \]
    where $S_{1},S_{2}\in\mathscr{N}(X)$ are two left linear quaternionic operators, then 
    \begin{align*} 
        S(qx)=&~(S_{0}+S_{1}J)(qx)=qS_{0}(x)+T_{1}(qx)J\\
        =&~qS_{0}(x)+qS_{1}(x)J=qS_{0}(x)+q(S_{1}J)(x)\\
        =&~qS(x) \quad (q\in\mathbb{H},x\in X),
    \end{align*} 
    which implies that $S$ is left linear quaternionic operator.
    Analogously, the same conclusion holds true for any element
    in $\mathbb{C}^{\mathcal{S}}_{R}(X)$.
\end{remark}

    By Proposition \ref{kexi-fenli0}, we can give the following definition
    of slice regularity for operator-valued functions,
    which is equivalent to Definition \ref{ke-def5}.
\begin{definition}
    Let $U$ be a slice axially symmetric domain associated to $\mathscr{D}$ and $X$ be a
    Banach bilateral $\mathbb{H}$-module. Suppose that $f:U\rightarrow \mathcal{B}^{R}(X)$ is an
    operator-valued function, then $f$ is called left $($resp. right$)$ slice regular if
    for any $z\in\mathscr{D}$,
    \[
        f\big(\phi_{I}(z)\big)=F\big(\phi_{I}(z)\big)+G\big(\phi_{I}(z)\big)J\quad
        \big(\text{resp. }f(\phi_{I}(z))=F(\phi_{I}(z))+JG(\phi_{I}(z))\big),
    \]
    where $I\perp J$ and $F,G:U\cap\mathbb{C}_{I}\rightarrow \mathbb{C}_{I}^{\mathcal{S}}(X)\cap\mathcal{B}^{R}(X)$, $I\in\mathbb{S}$.
\end{definition}

In Definition 2.1.16 from \cite{kexi-co3},
Colombo et al. introduced the hyperholomorphic product
structure for quaternion-valued
slice hyperholomorphic functions, while in Definition 5.3 from \cite{kexi-srs0},
Ghiloni et al. gave the $\odot$-product for right slice functions valued in
a Banach bilateral $\mathbb{A}$-algebra $\mathscr{L}^{R}(X)$,
which can be applied to quaternionic setting.

\begin{definition}[\cite{kexi-srs0}]
    Let $\Omega_{\mathscr{D}}$ be a slice axially symmetric domain
    associated to $\mathscr{D}$, let
    $X$ be a Banach bilateral $\mathbb{H}$-module.
    Consider two left slice functions
    $f : \Omega_{\mathscr{D}} \longrightarrow \mathcal{B}^{R}(X)$ and
    $g : \Omega_{\mathscr{D}} \longrightarrow \mathcal{B}^{R}(X)$ such that
    \[
        f\big(\phi_{I}(z)\big)=F_{1}(z)+IF_{2}(z) ,\quad g\big(\phi(z)\big)=G_{1}(z)+IG_{2}(z),
    \]
    where $F_{1}$ and $F_{2}$, as well as $G_{1}$ and $G_{2}$,
    satisfy \eqref{xrtj0}, respectively.
    Then define the left slice product as follows:
    \[
        (f\odot^{L} g)\big(\phi_{I}(z)\big)
        =\big(F_{1}(z)\circ G_{1}(z)-F_{2}(z)\circ G_{2}(z)\big)
        +I\big(F_{2}(z)\circ G_{1}(z)+F_{1}(z)\circ G_{2}(z)\big),
    \]
    where $z\in\mathscr{D}$, $I\in\mathbb{S}$ and
    $ A \circ B $ denotes the composition of operators $ A $ and $ B $,
    i.e.,  $(A \circ B)x = (AB)x = A(Bx)$.

    Similarly, for
    two right slice functions
    $f : \Omega_{\mathscr{D}} \longrightarrow \mathcal{B}^{R}(X)$ and
    $g : \Omega_{\mathscr{D}} \longrightarrow \mathcal{B}^{R}(X)$ such that
    \[
        f\big(\phi_{I}(z)\big)=F_{1}(z)+F_{2}(z)I ,\quad g\big(\phi(z)\big)=G_{1}(z)+G_{2}(z)I,
    \]
    where $F_{1}$ and $F_{2}$, as well as $G_{1}$ and $G_{2}$,
    satisfy \eqref{xrtj0}, respectively.
    Then define the right slice product as follows:
    \[
        (f\odot^{R} g)\big(\phi_{I}(z)\big)
        =\big(F_{1}(z)\circ G_{1}(z)-F_{2}(z)\circ G_{2}(z)\big)
        +\big(F_{2}(z)\circ G_{1}(z)+F_{1}(z)\circ G_{2}(z)\big)I,
    \]
    where $z\in\mathscr{D}$, $I\in\mathbb{S}$.
\end{definition}

\begin{remark}\label{kexi-rem3}
    Note that if
    $f: \Omega_{\mathscr{D}} \rightarrow \mathcal{B}^{R}(X)$
    is a constant operator-valued function,
    then $f$ is both left slice and
    right slice.
    Assume that $f(\phi_I(z)) = A$ and
    $g: \Omega_{\mathscr{D}} \rightarrow \mathcal{B}^{R}(X)$
    is a right slice function, then
    \begin{align*}
        (A \odot^{R} g)\big(\phi_{I}(z)\big)=&~
       ( f\odot^{R}g)\big(\phi_{I}(z)\big)\\
       =&~\big(F_{1}(z)\circ G_{1}(z)-F_{2}(z)\circ G_{2}(z)\big)
        +\big(F_{2}(z)\circ G_{1}(z)+F_{1}(z)\circ G_{2}(z)\big)I \\
        =&~A \circ G_{1}(z) + \big(A\circ G_{2}(z)\big) I\\
        =&~ A\circ g\big(\phi_{I}(z)\big).
    \end{align*}
\end{remark}

\begin{remark}\label{kexi-jishu-03}
Note that if
$A:\Omega_{\mathscr D}\to \mathcal{B}^{R}(X)$
is a right slice operator-valued function, and
$
u:\Omega_{\mathscr D}\to X
$
is a right slice function satisfying
\begin{equation}\label{kexi-jishu-11}
u(q)\in \mathcal{D}\big(A(q)\big), \qquad q\in\Omega_{\mathscr D},
\end{equation}
then, in the subsequent discussion, the notation $A(q)u(q)$ is understood in the sense that
\[
 A \big(\phi_{I}(z)\big)u\big(\phi_{I}(z)\big)
        :=\big(A_{1}(z) u_{1}(z)-A_{2}(z) u_{2}(z)\big)
        +\big(A_{2}(z) u_{1}(z)+A_{1}(z) u_{2}(z)\big)I,
\]
where $q=\phi_{I}(z)$ $(z\in\mathscr{D},I\in\mathbb{S})$ and
\[
    A\big(\phi_{I}(z)\big)=A_{1}(z)+A_{2}(z)I,\quad u\big(\phi_{I}(z)\big)=u_{1}(z)+u_{2}(z)I .
\]
It can be verified that if $A(q)$ and $u(q)$ are right slice regular functions, then
$f(q):=A(q)u(q)$ remains a right slice regular function.
\end{remark}

Below we will introduce the concept of the quaternionic
$\sigma$-ball (see Gentili and Stoppato \cite{kexi-jishu}), which will play a foundational role in geometric theories such as the slice regularity of functions and the extension of evolution dynamical systems. Later, we will further extend it to the framework of associative real $*$-algebras.

\begin{definition}[\cite{kexi-jishu}]
    For all $p,q\in\mathbb{H}$, define
    \[
        \omega(q,p)=\sqrt{\big(\mathrm{Re}(q)-\mathrm{Re}(p)\big)^{2}+\big(|\mathrm{Im}(q)|+|\mathrm{Im}(p)|\big)^{2}}
    \]
    and
    \begin{equation*}
        \sigma(q,p)=\left\{
            \begin{aligned}
                &|q-p|\quad &&\text{if }p,q \text{ lie on the same complex plane } \mathbb{C}_{I},\\
                & \omega(q,p)&& \text{otherwise}.
            \end{aligned}
        \right.
    \end{equation*}
    In addition, define the $\sigma$-ball $\Sigma(p,R)$ as follows:
    \[
        \Sigma(p,R):=\{q\in\mathbb{H}:\sigma(q,p)<R,R\in\mathbb{R}^{+}\}.
    \]
\end{definition}

Next, to continue our discussion exactly,
we introduce the following $\ast_{R}$-product
(see Gentili and Stoppato \cite{kexi-jishu} for the left version):
\[
    (f\ast_{R} g)(q)=\sum_{\ell=0}^{n}a_{\ell}b_{n-\ell}\sum_{n\in\mathbb{N}} q^{n},
\]
where
\[
    f(q)=\sum_{n\in\mathbb{N}} a_{n}q^{n},\quad g(q)=\sum_{n\in\mathbb{N}}b_{n}q^{n}\quad
    (a_{n},b_{n}\in X).
\]

\begin{remark}\label{kexi-re0-1}
Let $X$ be a Banach bilateral $\mathbb{H}$-module.
In \cite{kexi-xty}, the authors proved the local power series expansion of
left slice regular functions that valued in a
Banach bilateral $\mathbb{H}$-module by introducing the multiplicative
structure given in \cite{kexi-jishu}:
\[
    (f\ast_{L} g)(q)=\sum_{n\in\mathbb{N}} q^{n}\sum_{\ell=0}^{n}a_{\ell}b_{n-\ell},
\]
where
\[
    f(q)=\sum_{n\in\mathbb{N}} q^{n}a_{n},\quad g(q)=\sum_{n\in\mathbb{N}}q^{n}b_{n}\quad
    (a_{n},b_{n}\in X).
\]
Based on this,
for the right slice regular
operator-valued function
$ A:U\subset \mathbb{H}\rightarrow \mathcal{B}^{R}(X)$, one has
\[
    A(q)=\sum_{n\in\mathbb{N}}\big((\partial_{S}^{R})^{(n)}A(p)\big)\frac{1}{n!}
    (q-p)^{\ast_{R}n},
\]
where $(q-p)^{\ast_{R} n}=(q-p){\ast_{R}} \cdots {\ast_{R}} (q-p)$ and $q\in \Sigma(p,\rho)$ for some
$\rho>0$
and
$(\partial_{S}^{R})^{(n)}A(p)$ denotes the $n$-th right slice derivative of $A(q)$.
Therefore, for $x\in\mathcal{D}\big(A(q)\big)$, by Theorem 13 in \cite{kexi-jishu},
we have
$A(q)x$
is right slice regular.

Furthermore, if $A:\Omega_{\mathscr{D}}\rightarrow\mathcal{B}^{R}(X)$ and
$u:\Omega\rightarrow X$ are right slice regular functions
such that \eqref{kexi-jishu-11} holds and for $p\in\Omega_{\mathscr{D}}$,
we have the following
\begin{align*}
    A(q)=\sum_{\kappa=0}^{\infty}A_{\kappa}(q-p)^{\ast_{R}\kappa},\quad
    u(q)=\sum_{\ell=0}^{\infty}a_{\ell}(q-p)^{\ast_{R}\ell}
    \quad (q\in\Omega_{\mathscr{D}}).
\end{align*}
Then we define
\begin{align*}
    A(q)u(q):=
    \sum_{n=0}^{\infty}\bigg(\sum_{\ell=0}^{\infty}A_{n-\ell}a_{\ell}\bigg)(q-p)^{\ast_{R}n}.
\end{align*}
Thus, by the right version of Theorem 13 in \cite{kexi-jishu} again, $A(q)u(q)$ remains
a right slice regular function.

\end{remark}

\begin{remark}\label{kexi-rem1}
   Note that if $U\subset \mathbb{H}$ is a slice
    axially symmetric domain, then for a slice regular operator-valued function
    $A:U\rightarrow\mathcal{B}^{R}(X)$
    and a given vector $u_0 \in \mathcal{D}\big(A(q)\big)$, for all $q=\eta+\xi I$,
    we denote $A_{u_0}(q) = A(q)u_0$. Then it follows from Remark 2.1.13 in
    \cite{kexi-co3} and Remark \ref{kexi-re0-1} that
\begin{align*}
\partial^{R}_{S} \big(A(q)u_{0}\big)
=&~\partial_{S}^{R}\big(A_{u_{0}}(q)\big)
= \partial_{S}^{R}\big(A_{u_{0}}(\eta+\xi I)\big)
=\frac{\partial}{\partial \eta}\big(A_{u_{0}}(\eta+\xi I)\big)\\
=&~\lim_{\Delta \eta\rightarrow 0}
\frac{A_{u_{0}} (\eta+\Delta\eta+\xi I)-A_{u_{0}}(\eta+\xi I)}{\Delta \eta}\\
=&~\lim_{\Delta\eta\rightarrow 0}
\frac{A(\eta+\Delta\eta+\xi I)u_{0}-A(\eta+\xi I)u_{0}}{\Delta \eta}\\
=&~\lim_{\Delta\eta\rightarrow 0}
\frac{A(\eta+\Delta\eta+\xi I) -A(\eta+\xi I) }{\Delta \eta}u_{0}
=\big(\partial_{S}^{R}A(q)\big)u_{0}.
\end{align*}
\end{remark}

Let $\varphi\in(0,\pi]$, the so called spherical sector $\Sigma_{\varphi}$ is defined as follows:
\[
    \Sigma_{\varphi}:=\big\{s\in\mathbb{H}\backslash\{0\}:\mathrm{arg}(s)<\varphi\big\}.
\]
Moreover, we will set
\[
    \Sigma_{\varphi}':=\Sigma_{\varphi}\cup\{0\} .
\]

\begin{definition}[\cite{kexi-srs0}]\label{kexi-lian}
    Let $X$ be a Banach bilateral $\mathbb{H}$-module and let $\theta \in (0, \pi/2]$. 
    Then a map
    $\mathcal{S} : \Sigma_{\theta}'\rightarrow \mathcal{B}^{R}(X)$
    will be called a quaternionic
    right slice regular semigroup $($of angle $\theta$$)$ if the restriction
    $\mathcal{S}\vert_{\Sigma_{\theta}} : \Sigma_{\theta}\rightarrow \mathcal{B}^{R}(X)$
    is right slice regular and
    \begin{align*}
        &\mathcal{S}(p + q) = \mathcal{S}(p) \odot_p^{R} \mathcal{S}(q) \quad \forall p, q
        \in \Sigma_{\theta} \text{ with } p + q \in \Sigma_{\theta} \text{ and } pq = qp, \\
    &\mathcal{S}(0) = \mathcal{I}, \\
    &\lim_{q \to 0} \mathcal{S}|_{\Sigma_{\theta'}}(q)x
    = x \quad \forall \theta' \in (0, \theta), \quad \forall x \in X,
\end{align*}
where $\odot^{R}_{p}$ denotes the right slice product
emphasizing the first variable $p$ $($see \cite{kexi-srs0} p. 5009$)$.
    If $\mathcal{S}$ is a quaternionic right slice regular semigroup of
    angle $\theta$ for some
    $\theta \in (0, \pi/2]$,
    for convenience, $T$ will be referred to as
    quaternionic right slice regular semigroup. On the other hand,
    $T$ is called a bounded quaternionic right slice regular semigroup if
    \[
        \forall \delta' \in (0, \delta) \quad \exists M_\delta \in
        [1, \infty) : \sup_{q \in \Sigma_{\delta'}} \|\mathcal{S}(q)\| \leq M_{\delta'}.
    \]
    Furthermore, for the right slice regular semigroup 
    $\mathcal{S}:\Sigma_{\theta}'\rightarrow\mathcal{B}^{R}(X)$, if 
    an operator $A : D(A) \to X $ is the generator of the strongly continuous
    quaternionic semigroup $ \mathcal{S}|_{[0,\infty)}$, then $T$
    is called the generator of $ \mathcal{S}$, or that $A$ generates $\mathcal{S}$.
\end{definition}

For the sake of clarity,
    we will use $\odot^{R}$ to denote the right slice product unless otherwise specified
    in the subsequent discussion.

\begin{remark}\label{kexi-rem4}
    Note that assume $A$ is a
    dense defined spherical $\theta$-sectorial operator with
    vertex $\omega$, if there exists a constant $M>0$ such that
    \[
        \|S_{L}^{-1}(q,A)\|_{\mathcal{B}^{R}(X)} \leq \frac{M}{|q-\omega|}\quad \forall q\in\omega+\Sigma_{\pi/2+\theta},
    \]
    from Theorem 7.1 in \cite{kexi-srs0}, it follows that
    $A$ generates a quaternionic right slice regular semigroup
    $\mathcal{S}:\Sigma_{\theta}\rightarrow \mathcal{B}^{R}(X)$.
    Furthermore, Theorem 7.2 in \cite{kexi-srs0} ensures the
    uniqueness of this semigroup, which is given as follows:
    \begin{equation}\label{scy0}
        \mathcal{S}(q):=\frac{1}{2\pi} \int_{\omega+\Gamma(J;r;\eta)}S_{L}^{-1}(\alpha,A)
        J^{-1}~\mathrm{d}\alpha \exp^{\alpha}(q)\quad (q\in\Sigma_{\theta}).
    \end{equation}
   Following the construction in \cite{kexi-srs0} p. 5021
    $($see also \cite{kexi-srs1} p. 2669$)$, the integration path
    $\omega+\Gamma(J;r,\eta)\subset \Sigma_{\theta}\cap\mathbb{C}_{J}$ is chosen to be a $C^{1}$-path.
By a $C^{1}$-path in $\Sigma_{\theta}\cap\mathbb{C}_{J}$ we mean a $C^{1}$-mapping
$
\gamma : [a,b] \to \Sigma_{\theta}\cap\mathbb{C}_{J}
$
which is either injective, or else fails to be injective only at the endpoints, in the sense that
\[
\gamma(t)=\gamma(s) \implies (t,s)=(a,b)\ \text{or}\ (t,s)=(b,a).
\]
Finally, the exponential function $\exp^{\alpha}(q)$
    was defined by
Definition 5.1 in \cite{kexi-srs0}.
    Taking the right derivative on both sides of \eqref{scy0}, we obtain that
    \begin{align}\label{kexi-fangsuo-0}
        \partial_{S}^{R}\mathcal{S}(q)
        =\frac{1}{2\pi} \int_{\omega+\Gamma(J;r;\eta)}S_{L}^{-1}(\alpha,A)\alpha
        J^{-1}~\mathrm{d}\alpha~  \exp^{\alpha}(q),
    \end{align}
    by Theorem 2.18 in \cite{kexi-co1} and Remark \ref{kexi-rem3},  it follows that
    \begin{align}\label{kexi-tianjia-1}
        \partial_{S}^{R}\mathcal{S}(q)
        =\frac{1}{2\pi} \int_{\omega+\Gamma(J;r;\eta)} A  S_{L}^{-1}(\alpha,A)
        J^{-1}~\mathrm{d}\alpha~  \exp^{\alpha}(q)
        =A \circ\mathcal{S}(q)=A \odot^{R}\mathcal{S}(q).
    \end{align}
\end{remark}

\section{Evolution Dynamics and Cauchy Problems with Quaternionic
Variables}\label{sec4}

In this section, we develop the theory of evolution dynamics within the framework of Banach bilateral
$\mathbb{H}$-algebras. This allows us to establish a complete theory of Cauchy problems for evolution equations involving quaternionic variables, including the extension of solutions and the existence and uniqueness theory. Within the framework of Banach bilateral
$\mathbb{H}$-algebras, we introduce the rigorous notion of quaternionic evolution systems. For the first time, we obtain the geometric formulation of the extension and uniqueness of solutions to Cauchy problems with quaternionic variables on the closure of quaternionic
$\sigma$-ball, and we establish the corresponding main theorems.

\subsection{
Cauchy Problem and Evolution Systems in Banach Bilateral
$\mathbb{H}$-Algebras}

\subsubsection{Existence and Uniqueness for
Right $\sigma$-Analytic Solutions in Banach Bilateral
$\mathbb{H}$-Module}

In this subsection,
we will study the following first-order Cauchy problem with quaternionic variables
as follows:
\begin{equation}\label{kexi-e4.1}
    \left\{\begin{aligned}
     &\partial_{S}^{R} u (q)
    =A(q)u(q)+f(q), \\
     &u(s_{0})=u_{0},
    \end{aligned}\right.\qquad s_{0},q\in \Sigma_{\theta}.
\end{equation}
    where $A(q): \mathcal{D}(A(q)) \subset X \to X$ is a right linear closed
    quaternionic operator on $X$, $\theta\in(0,\pi/2]$,
    $f:\Sigma_{\theta}\rightarrow X$ is a given slice regular function and $X$ is a Banach bilateral $\mathbb{H}$-module.

    \begin{definition}
    A function $u: \Sigma_{\theta} \rightarrow X$ is
    called a classical solution
    of \eqref{kexi-e4.1} if  $u(q)$ satisfies \eqref{kexi-e4.1} and
    $u(q) \in \mathcal{D}\big(A(q)\big)$  is right slice regular for each $q\in \Sigma_{\theta}$.
\end{definition}

We begin by addressing the homogeneous case:
\begin{equation}\label{kexi-e2}
    \left\{\begin{aligned}
     &\partial_{S}^{R}u (q)=A(q)u(q),  \\
        &u(s_{0})=u_{0},
    \end{aligned}\right.
\end{equation}
where $s_{0},q\in\Sigma_{\theta}$.

To establish the existence and uniqueness of solutions to
\eqref{kexi-e2}, we will introduce the notion of $\sigma$-analyticity on the numerical range (see \cite{kexi-jishu}).

\begin{definition}[\cite{kexi-jishu}]
  Let $U$ be a $\sigma$-ball and $X$ be a
   Banach bilateral $\mathbb{H}$-module. A function $f:U\rightarrow X$ is called left $($resp. right$)$
    $\sigma$-analytic at $p\in U$ if there exists power series
    $\sum_{n\in\mathbb{N}}(q-p)^{{\ast_{L}} n}a_{n}$
    $($resp. $\sum_{n\in\mathbb{N}}a_{n}(q-p)^{{\ast_{R}} n}$$)$
    converging in a $\sigma$-neighborhood
    $U'$ of $p$ in $U$ such that $f(q)=\sum_{n\in\mathbb{N}}(q-p)^{{\ast_{L}} n}a_{n}$
    $($resp. $f(q)=\sum_{n\in\mathbb{N}}a_{n}(q-p)^{{\ast_{R}} n})$ for all $q\in U'$. Furthermore,
    $f$ is called left $($resp. right$)$ $\sigma$-analytic if it is left $($resp. right$)$
    $\sigma$ at all $p\in U$.
\end{definition}

In the following, within the quaternionic
$\sigma$-ball, we establish, for the first time, a local existence and uniqueness for
$\sigma$-analytic solutions of the evolution equation \eqref{kexi-e2} with quaternionic variables. On the Banach bilateral
$\mathbb{H}$-module, and employing the slice geometry theory of quaternionic space, we extend the
$\sigma$-analytic solution beyond the boundary of the closure of the quaternionic
$\sigma$-ball, thereby proving a global existence and uniqueness theorem for
$\sigma$-analytic solutions of \eqref{kexi-e2}.

\begin{theorem}[The local existence and uniqueness theorem]\label{kexi-t4.1}
    Let $X$ be a Banach bilateral $\mathbb{H}$-module and
   $\theta\in(0,\pi/2]$.
    Assume that the operator-valued function
    $A: \Sigma_{\theta}\to \mathcal{B}^{R}(X)$
    is right slice regular and $A(q)$ is uniformly bounded closed right linear quaternionic 
    operator for each $q\in\Sigma_{\theta}$, then
    \eqref{kexi-e2} has a unique
    right $\sigma$-analytic solution $u$
    in
    $\overline{\Sigma(s_{0},R)}\subset \Sigma_{\theta}$.
\end{theorem}

\begin{proof}
    Without loss of generality, we will assume that
    $s_{0}\in\mathbb{C}_{I}$ $(I\in\mathbb{S})$, and
    for convenience, we will use the notation
    $\hat{\Sigma}:=\overline{\Sigma(s_{0},R)}\cap\mathbb{C}_{I}$.
    Let $\mathcal{H} ( \hat{\Sigma},X)$
    be the set of all right
    hyperholomorphic function from $\hat{\Sigma}$ to $X$. Then
    $\mathcal{H}  (\hat{\Sigma},X )$ equipped with the usual norm
    forms a Banach space.
    Since $A(q)$ is uniformly bounded, now
    let $\alpha_{\mathscr{D}} = \sup_{q\in\Sigma_{\theta}} \|A(q)\|_{\mathcal{B}^{R}(X)}$
    and $T$ a mapping
    from $\mathcal{H}  (\hat{\Sigma},X)$  to itself given by
    \begin{equation*}
        (Tu)(q) = u_{0} + \int_{\gamma_{[s_{0},q,I]}} A(\xi)u(\xi)~\mathrm{d}\xi,
    \end{equation*}
    where
    $q \in  \hat{\Sigma}  \subset \Sigma_{\theta}$ and
    $\gamma_{[s_{0},q,I]} :[a,b]\rightarrow  \Sigma_{\theta} \cap\mathbb{C}_{I}$ is
    $C^{1}$-path from
    $s_{0}$ to $q$,
    i.e., $\gamma_{[s_{0},q,I]} (a)=s_{0}$ and $\gamma_{[s_{0},q,I]} (b)=q$.
    Denote $\|u\|_{\infty} = \sup_{q\in \Sigma_{\theta}} \|u(q)\|_{X}$, then
    \begin{equation*}
        \|Tu(q) - Tv(q)\|_{X} \leq \alpha_{\mathscr{D}}
        \vert \gamma_{[s_{0},q,I]} \vert \|u - v\|_{\infty},
    \end{equation*}
    where $\vert \gamma_{[s_{0},q,I]} \vert$ denotes the arc length of $\gamma_{[s_{0},q,I]}$ on the complex plane $\mathbb{C}_{I}$.
    It follows by induction that
    \[
    \|T^{n}u(q) - T^{n}v(q)\|_{X} \leq \frac{\alpha_{\mathscr{D}}^{n}
    \vert {\gamma_{[s_{0},q,I]}} \vert^{n}}{n!}
    \|u - v\|_{\infty},
    \]
     which implies
    \[
        \|T^{n}u- T^{n}v\|_{\infty} \leq \frac{\alpha_{\mathscr{D}}^{n}
        \vert {\gamma_{[s_{0},q,I]}}  \vert^{n}}{n!}
        \|u - v\|_{\infty}.
    \]
    Since $\alpha_{\mathscr{D}}^{n}\vert {\gamma_{[s_{0},q,I]}} \vert^{n}/n! < 1$ as
    $n\rightarrow\infty$, then
    the Banach contraction principle
    yields that $T$ possesses a unique fixed point $u$
    within the right hyperholomorphic function space
    $\mathcal{H}  (\hat{\Sigma},X )$, which implies \eqref{kexi-e2} has a unique solution.
    Considering the following integral equation
    \begin{equation}\label{kexi-solu1}
        u(q)=u_{0}+\int_{\gamma_{[s_{0},q,I]}}A(\xi)u(\xi)~\mathrm{d}\xi.
    \end{equation}
    Then \eqref{kexi-solu1} is the corresponding integral equation of \eqref{kexi-e2} restricted to $\hat{\Sigma}$.
    By the right version of the Banach bilateral
    $\mathbb{H}$-module form of Corollary 2 in \cite{kexi-jishu},
    $u(q)$ is right $\sigma$-analytic function restricted to $\hat{\Sigma}$. Thus,
    $u(q)$ $(q\in \hat{\Sigma})$ have the following Laurent series expansion
    \begin{align}\label{kexi-jishu-01}
        u(q)=\sum_{n=0}^{\infty}a_{n}(q-s_{0})^{\ast_{R}n}\quad (a_{n}\in X).
    \end{align}
    Since $A(q)$ is right slice regular, then we can assume that
    \[
        A(q)=\sum_{n=0}^{\infty}A_{n}(q-s_{0})^{\ast_{R}n},
    \]
    and hence \eqref{kexi-jishu-11} can be rewritten as:
    \begin{equation}\label{kexi-jishu-02}
        u(q)=u_{0}+\sum_{n=0}^{\infty}\frac{1}{n+1}
        \bigg(\sum_{\ell=0}^{n}A_{\ell}a_{\ell}\bigg)(q-s_{0})^{\ast_{R}(n+1)}\quad
        (q\in\hat{\Sigma}).
    \end{equation}
    According to \eqref{kexi-jishu-01} and \eqref{kexi-jishu-02}, we have
    \[
        u(q)=u_{0}+\sum_{n=1}^{\infty}\bigg(\frac{1}{n}\sum_{\ell=0}^{n-1}A_{n-1-\ell}a_{\ell}\bigg)
        (q-s_{0})^{\ast_{R}n}.
    \]

    Next, by virtue of Proposition \ref{mi0},  we can transform the geometric slices of the locally unique solution on plane $\mathbb{C}_{I}$ to plane $\mathbb{C}_{J}$, constructing the quaternionic variable solution that satisfies the evolution equation \eqref{kexi-e2} as follows
    \begin{equation}\label{kexi-jishu-04}
        u\big(\phi_{J}(z)\big)=\left\{
        \begin{aligned}
            &u_{0}+\sum_{n=1}^{\infty}\bigg(\frac{1}{n}\sum_{\ell=0}^{n-1}A_{n-1-\ell}a_{\ell}\bigg)
            (q-s_{0})^{\ast_{R}n},
            \quad &&J=I,\\
            &\frac{1}{2}u\big(\phi_{I}(\overline{z})\big)(1+IJ)+\frac{1}{2}u\big(\phi_{I}(z)\big)(1-IJ),
            &&J\neq I,
        \end{aligned}\right.
    \end{equation}
    where $\phi_{J}(z)=q\in \overline{\Sigma(s_{0},R)}$.
    Obviously, for the form constructed at $J\neq I$ in \eqref{kexi-jishu-04},
    we can also transform it into a series form in the same way as \eqref{kexi-jishu-02}.
    Therefore, the function $u(q)$ given by \eqref{kexi-jishu-04}
    is right $\sigma$-analytic in $\overline{\Sigma(s_{0},R)}$
    and that it is unique. Therefore, \eqref{kexi-e2} has a
    unique right $\sigma$-analytic solution in $\overline{\Sigma(s_{0},R)}$. This completes the proof.
\end{proof}

After establishing the local existence and uniqueness of the solution to the evolution equation \eqref{kexi-e2} on the quaternionic $\sigma$-ball $\overline{\Sigma(s_{0},R)}$, we can extend it beyond the closure of $\overline{\Sigma(s_{0},R)}$ and obtain a global quaternionic solution on the entire spherical sector region $\Sigma_{\theta}$, thereby establishing the continuation theorem for the solution.

\begin{theorem}\label{kexi-t4.2}
Let $X$ be a
Banach bilateral $\mathbb{H}$-module
and let $\Sigma_{\theta}\subset\mathbb{H}$ be an axially symmetric slice domain.
Assume that $ A(q) $
satisfies the hypotheses of Theorem \ref{kexi-t4.1},
then the right $\sigma$-analytic solution of \eqref{kexi-e2} on $\overline{\Sigma(s_{0},R)}$
can be extended to a global solution on $\Sigma_{\theta}$.
\end{theorem}

\begin{proof}
    By Theorem \ref{kexi-t4.1}, \eqref{kexi-e2} has a
    unique right $\sigma$-analytic solution $u_{1}$
    on $\overline{\Sigma(s_{0},R)}\subset \Sigma_{\theta}$.

    Now, let $s_{1}$ be an arbitrary boundary point of
    $\overline{\Sigma(s_{0},R)}\cap\mathbb{C}_{I}$,
    and consider $\overline{\Sigma(s_{1},R_{1})} \subset \Sigma_{\theta}$.
    Then, it follows from Theorem \ref{kexi-t4.1} that
    \eqref{kexi-e2} has a unique right $\sigma$-analytic solution $u_{2}$ on
    $\overline{\Sigma(s_{1},R_{1})}$. Note that at the point $s_{1}$, we have
$
u_{1}(s_{1}) = u_{2}(s_{1}).
$
Therefore, by the uniqueness of the solution, we obtain $u_{1}(q) = u_{2}(q)$ for
$q\in \overline{\Sigma(s_{0},R)}\cap\overline{\Sigma(s_{1},R_{1})}$. Set
\begin{equation*}
    u(q)=\left\{
        \begin{aligned}
            &u_{1}(q),\quad && q\in \overline{\Sigma(s_{0},R)}\\
            &u_{2}(q), && \overline{\Sigma(s_{1},R_{1})}\backslash \overline{\Sigma(s_{0},R)}.
        \end{aligned}\right.
\end{equation*}
Thus, we have constructed a solution of \eqref{kexi-e2} on a larger domain
$\overline{\Sigma(s_{0},R)}\cup \overline{\Sigma(s_{1},R_{1})}$.
Next, choose an arbitrary point $s_{2}$ on the boundary of the
$\big(\overline{\Sigma(s_{1},R_{1})}\cup\overline{\Sigma(s_{0},R)}\big)\cap\mathbb{C}_{I}$, and by repeating the above procedure,
we obtain a solution of \eqref{kexi-e2} on another, larger domain.
By iterating this process, the domain of the solution can be
consistently extended to the entire $\Sigma_{\theta}$.
According to a result analogous to Corollary 2 in \cite{kexi-jishu}, the conclusion follows.
  This completes the proof.
\end{proof}

\begin{figure}[htbp]
    \centering
        \centering
        \includegraphics[width=0.4\textwidth]{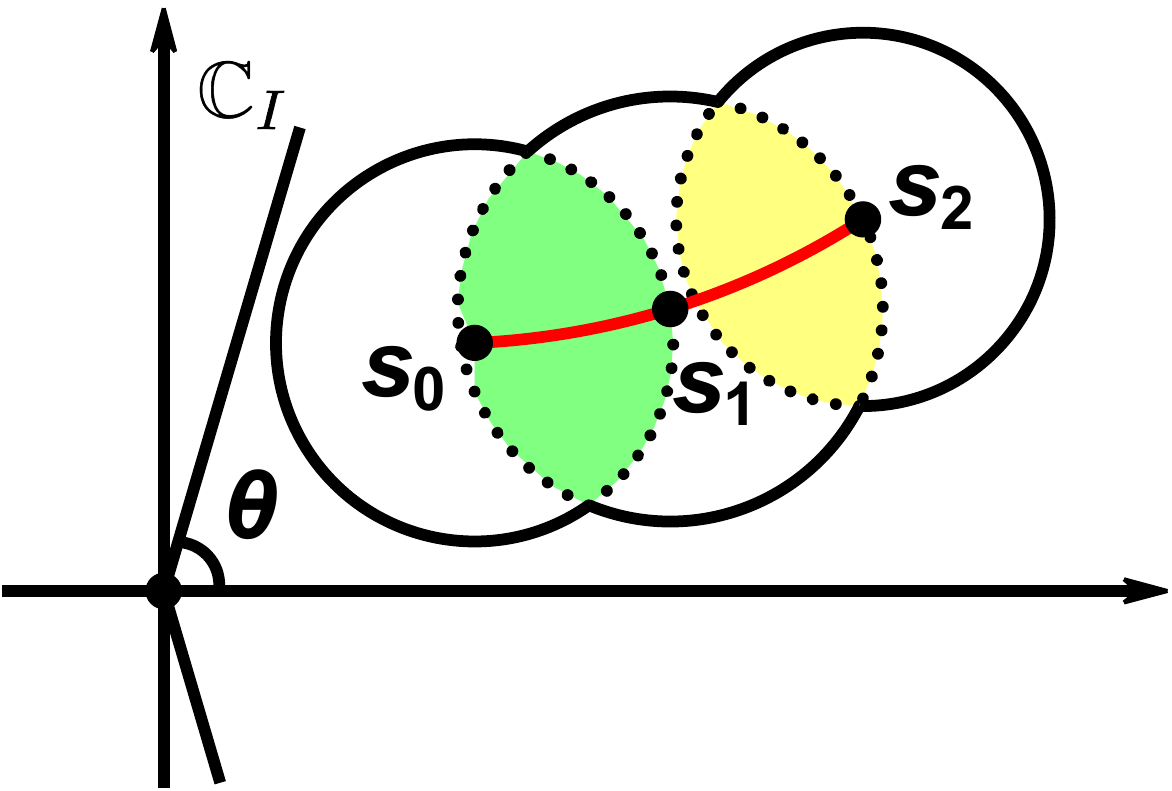}
        \caption{The case for $|\mathrm{Im}(s_{\ell})|\geq R$, $\ell=1,2,\cdots$}
        \label{kexi-yantuo333}
\end{figure}

\begin{figure}[htbp]
    \centering
    \begin{minipage}[b]{0.45\textwidth}
        \centering
        \includegraphics[width=\textwidth]{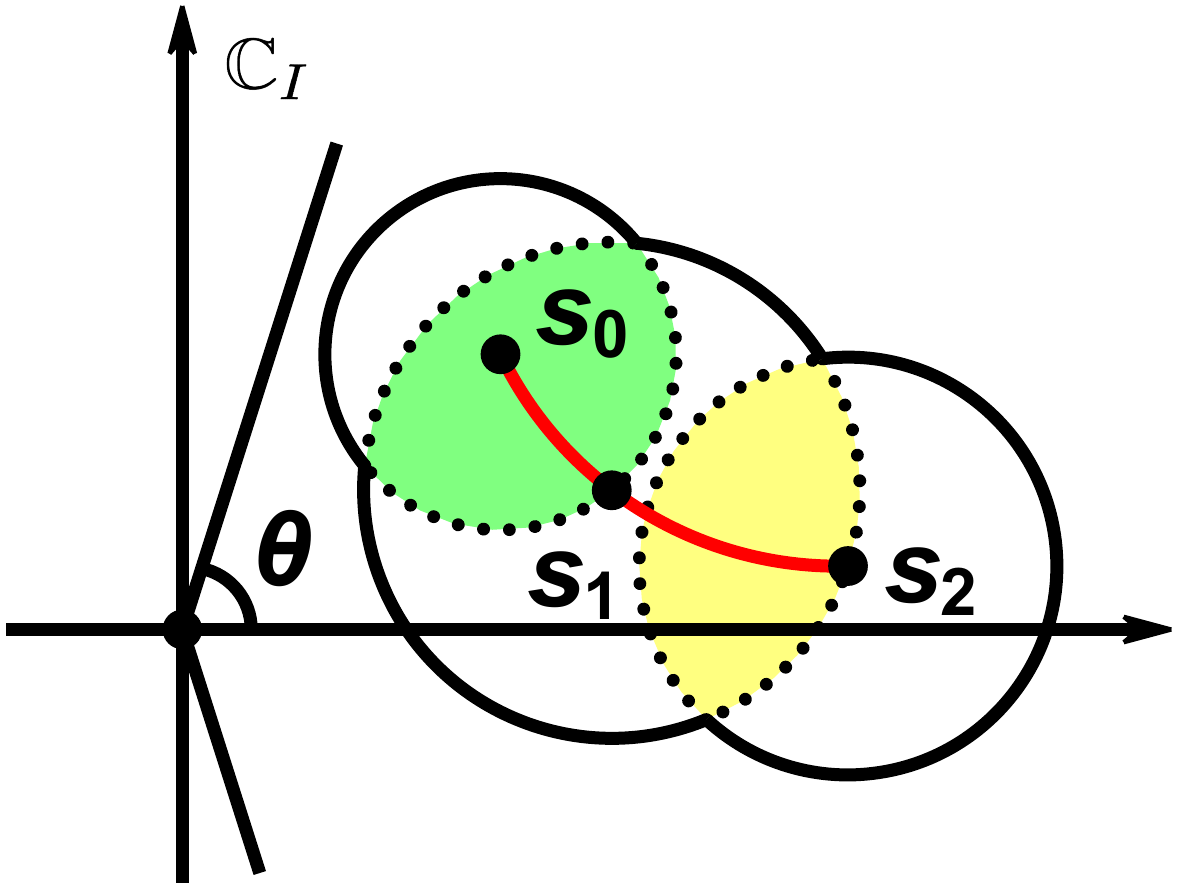}
    \end{minipage}
 \hfill
    \begin{minipage}[b]{0.45\textwidth}
        \centering
        \includegraphics[width=\textwidth]{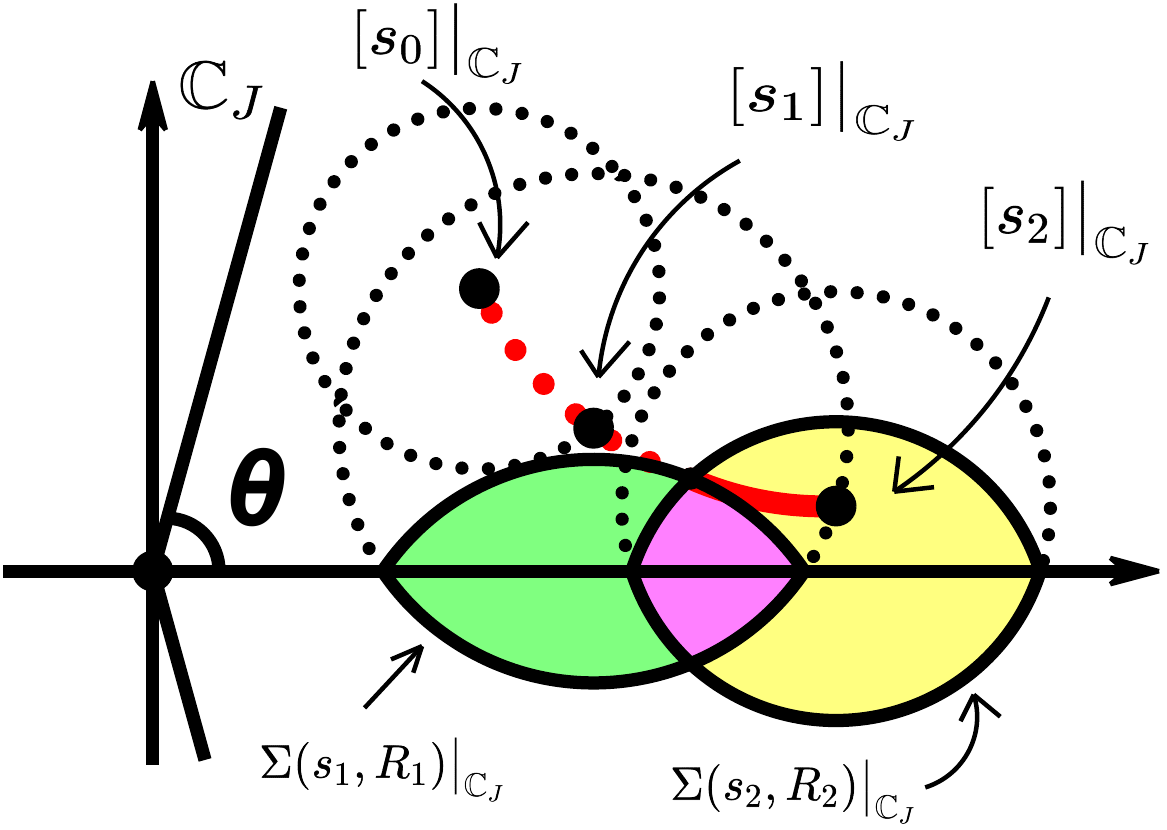}
    \end{minipage}
    \caption{The case for $0<|\mathrm{Im}(s_{\ell})|<R$, $\ell=1,2,\cdots$} \label{kexi-yantuo335}
\end{figure}

\begin{figure}[htbp]
    \centering
    \begin{minipage}[b]{0.45\textwidth}
        \centering
        \includegraphics[width=\textwidth]{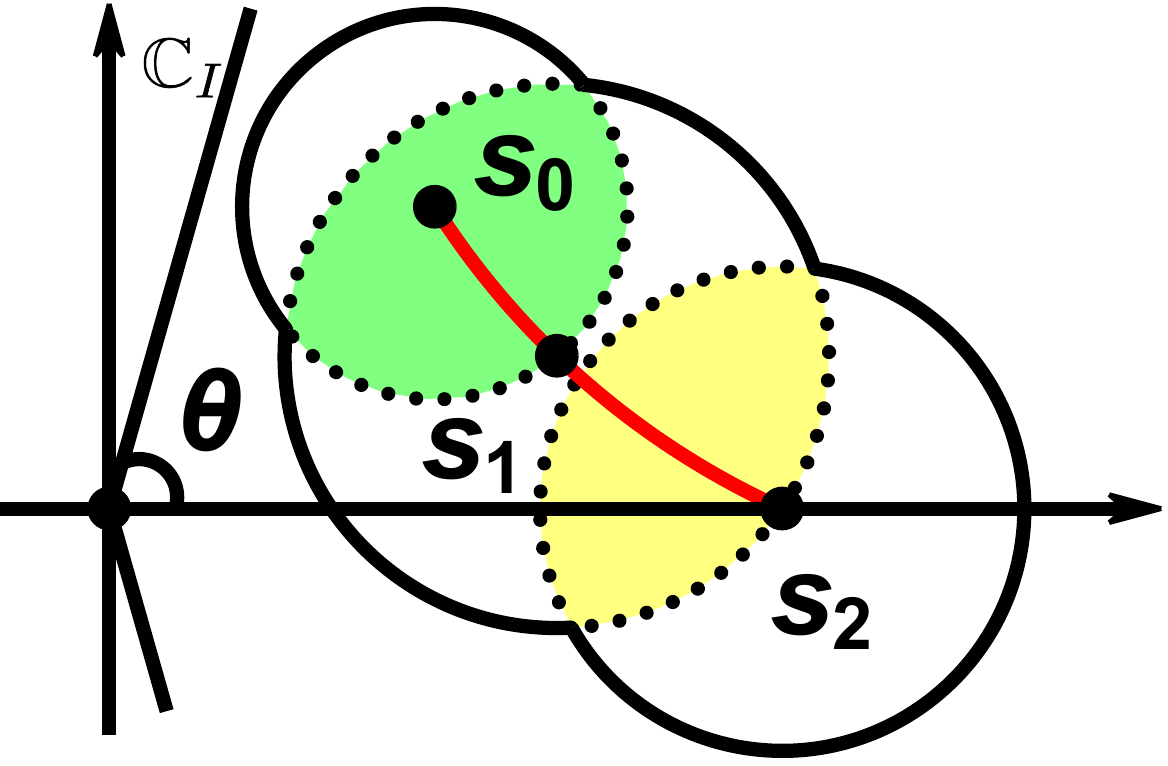}
    \end{minipage}
 \hfill
    \begin{minipage}[b]{0.5\textwidth}
        \centering
        \includegraphics[width=\textwidth]{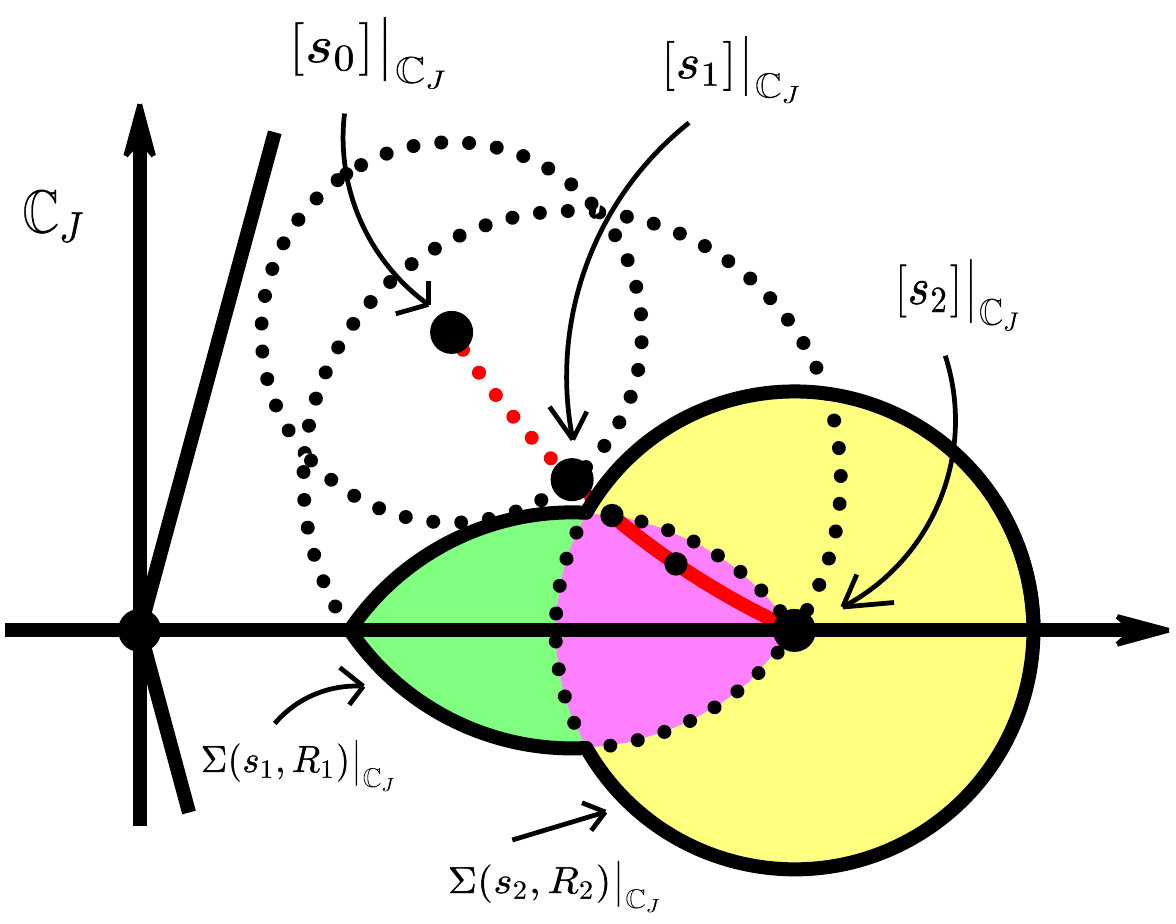}
    \end{minipage}
    \caption{The case for $0<|\mathrm{Im}(s_{1})|<R$ and
$ |\mathrm{Im}(s_{\ell})|=0$ $(\ell=2,3\cdots)$, $\ell=2,3,\cdots$} \label{kexi-yantuo336}
\end{figure}

Below, we describe the continuation of the solution of
\eqref{kexi-e2} through the geometric properties of the quaternionic $\sigma$-ball and the $\sigma$-analyticity of functions defined on it (for
more details, one may refer to Gentili and Stoppato \cite{kexi-jishu}).
Without loss of generality, we begin with the case
$
|\mathrm{Im}(s_{0})| \ge 0,
$
while the other cases can be treated analogously.
The continuation procedure can be divided into the following three cases:
\begin{itemize}
    \item[\rm(i)] The extended $\sigma$-balls $\overline{\Sigma(s_{\ell},R_{\ell})}$
    $(\ell=1,2,\cdots)$ satisfy
$|\mathrm{Im}(s_{\ell})|\geq R$ (see Figure \ref{kexi-yantuo333}).
    \item[\rm(ii)] The extended $\sigma$-balls $\overline{\Sigma(s_{\ell},R_{\ell})}$
    $(\ell=1,2,\cdots)$ satisfy the condition
$0<|\mathrm{Im}(s_{\ell})|<R$
(see Figure \ref{kexi-yantuo335}).
    \item[\rm(iii)]  The extended $\sigma$-balls $\overline{\Sigma(s_{\ell},R_{\ell})}$
    satisfy $0<|\mathrm{Im}(s_{1})|<R$ and
$ |\mathrm{Im}(s_{\ell})|=0$ $(\ell=2,3\cdots)$
(see Figure \ref{kexi-yantuo336}).
\end{itemize}

In the sequel,  we shall provide a unified representation for the solution of
\eqref{kexi-e2}, for which we need the following result.

\begin{proposition}\label{kexi-po-1}
Let $X$ be a Banach bilateral $\mathbb{H}$-module and let $\Omega_{\mathscr{D}}$
be a slice axially symmetric domain associated to $\mathscr{D}$.
Assume that $f$ is a right slice regular function on $\Omega_{\mathscr{D}}$.
For any $K\in\mathbb{S}$ and $q=x+Ky\in\Omega_{\mathscr{D}}$, set
\[
F(q)=\int_{\gamma_{[a+bK,q,K]}}f(\xi)~\mathrm{d}\xi\quad (a,b\in\mathbb{R}),
\]
where
$\gamma_{[s,q,K]} :[a,b]\rightarrow \Omega_{\mathscr{D}}\cap\mathbb{C}_{K}$
is a $C^{1}$-path   from
    $s$ to $q$ on $\mathbb{C}_{K}$.
Then for any $x+Jy\in \Omega_{\mathscr{D}}$ $(J\in\mathbb{S})$,
\begin{align*}
F(x+Jy)=&~\int_{\gamma_{[s,q,J]}}f(\xi)~\mathrm{d}\xi\\
=&~\frac{1}{2}F(x+Iy) (1-IJ)+\frac{1}{2}F(x-Iy) (1+IJ), \quad I\in\mathbb{S}.
\end{align*}
\end{proposition}

\begin{proof}
    Since $f$ is right slice regular, then
    \[
        f_{K}(x+Ky)=\alpha(x,y)+\beta(x,y)K,\quad K\in\mathbb{S},
    \]
    where $\alpha,\beta: \mathscr{D}\rightarrow X$. Suppose that $s=a+Jb$,
    consider the integral over $\mathbb{C}_{J}$:
    \begin{align*}
        \int_{\gamma_{[a+Jb,x+Jy,J]}}f(\xi)~\mathrm{d}\xi
        =&~ \int_{\gamma_{[a+Jb,x+Jb,J]}}f(\xi)~\mathrm{d}\xi
            +\int_{\gamma_{[x+Jb,x+Jy,J]}}f(\xi)~\mathrm{d}\xi\\
        =&~\int_{a}^{x}\alpha(\eta,b)~\mathrm{d}\eta+\int_{a}^{x}\beta(\eta,b)~\mathrm{d}\eta
        \cdot J
        +\int_{b}^{y}\alpha(x,\zeta)~\mathrm{d}\zeta\cdot J-\int_{b}^{y}\beta(x,\zeta)~\mathrm{d}\zeta.
    \end{align*}
        Let
    \begin{align*}
        P(x,y)=\int_{a}^{x}\alpha(\eta,b)~\mathrm{d}\eta
        -\int_{b}^{y}\beta(x,\zeta)~\mathrm{d}\zeta,\quad
        Q(x,y)=\int_{a}^{x}\beta(\eta,b)~\mathrm{d}\eta
        +\int_{b}^{y}\alpha(x,\zeta)~\mathrm{d}\zeta.
    \end{align*}
    Similarly, we obtain that
    \begin{align*}
         \int_{\gamma_{[a+Ib,x+Iy,I]}}f(\xi)~\mathrm{d}\xi=P(x,y)+Q(x,y)I, \quad
        \int_{\gamma_{[a+Ib,x-Iy,I]}}f(\xi)~\mathrm{d}\xi=P(x,y)-Q(x,y)I.
    \end{align*}
    Then $F(x+Iy)=P+QI$, $F(x-Iy)=P-QI$, and hence
    \begin{align*}
        &~\frac{1}{2}F(x+Iy) (1-IJ)+\frac{1}{2}F(x-Iy) (1+IJ)\\
        =&~\frac{1}{2}\big((P+QI)(1-IJ)\big) +
        \frac{1}{2}\big((P-QI)(1+IJ)\big)\\
        =&~P+QJ=F(x+yJ).
    \end{align*}
    This completes the proof.
\end{proof}

\begin{remark}
Note that the spherical sector $\Sigma_{\theta}$, where $\theta \in (0,\pi]$,
can be naturally interpreted as a slice axially symmetric domain $\Omega_{\mathscr{D}}$
in the sense of the Definition \ref{kexi-def3.3}.
Indeed, for any element $q = x + Iy \in \Sigma_{\theta}$ with $I \in \mathbb{S}$, there exists a complex number $z = x + iy$ belonging to the complex sector $\mathscr{D} = \{z \in \mathbb{C} : |\arg(z)| < \theta\}$ such that $\phi_{I}(z) = q$. Consequently, throughout the subsequent discussion, unless explicitly stated otherwise, the notation $\Omega_{\mathscr{D}}$ will be specifically employed to denote the spherical sector $\Sigma_{\theta}$, where the angle $\theta$ is taken from the interval $(0, \pi/2]$.
\end{remark}
    According to the right-side version of Theorem 2.3 in \cite{kexi-xty},
    along with Theorem \ref{kexi-t4.2}, Proposition \ref{mi0}, and Proposition \ref{kexi-po-1}, we can immediately obtain the following conclusion.
\begin{proposition}
    Let $X$ be a Banach bilateral $\mathbb{H}$-module.
    Assume that $A(q)$
    satisfies the hypotheses of Theorem \ref{kexi-t4.1} and $s_{0}=a+bJ=\phi_{J}(z_{0})$
    $(J\in\mathbb{S})$, then the classical solution
    of \eqref{kexi-e2} admits the following unified form:
    \begin{equation}\label{kexi-jie0}
        u\big(\phi_{I}(z)\big)=u_{I,0}+\int_{\gamma_{[\phi_{I}(z_{0}),\phi_{I}(z),I]}}
        A(\xi)u(\xi)~\mathrm{d}\xi,
    \end{equation}
    where $\phi_{I}(z_{0})=a+bI$, $z\in\mathscr{D}$ and
    \[
        u_{I,0}=u\big(\phi_{I}(z_{0})\big)=
        \left\{\begin{aligned}
        &u_{0}, &&I=J,\\
        &\frac{1}{2}u(\overline{s_{0}})(1+JI)+\frac{1}{2}u(s_{0})(1-JI),\quad &&I\neq J.
        \end{aligned}\right.
    \]
\end{proposition}

The following definitions and results will be used in our later discussion.
\begin{definition}
    For $\theta\in(0,\pi/2]$ and $p\in\mathbb{C}_{I}$ with $I\in\mathbb{S}$, define the 
    set
    \[
        \mathscr{X}_{\theta}(p):=\big\{q\in\mathbb{C}_{I}:\mathrm{Re}(q)\geq \mathrm{Re}(p)\text{ and }
        |\mathrm{Im}(q)-\mathrm{Im}(p)|\leq \big(\mathrm{Re}(q)-\mathrm{Re}(p)\big)\tan\theta\big\}.
    \]
\end{definition}

\begin{figure}[htbp]
        \centering{
        \includegraphics[scale=0.25]{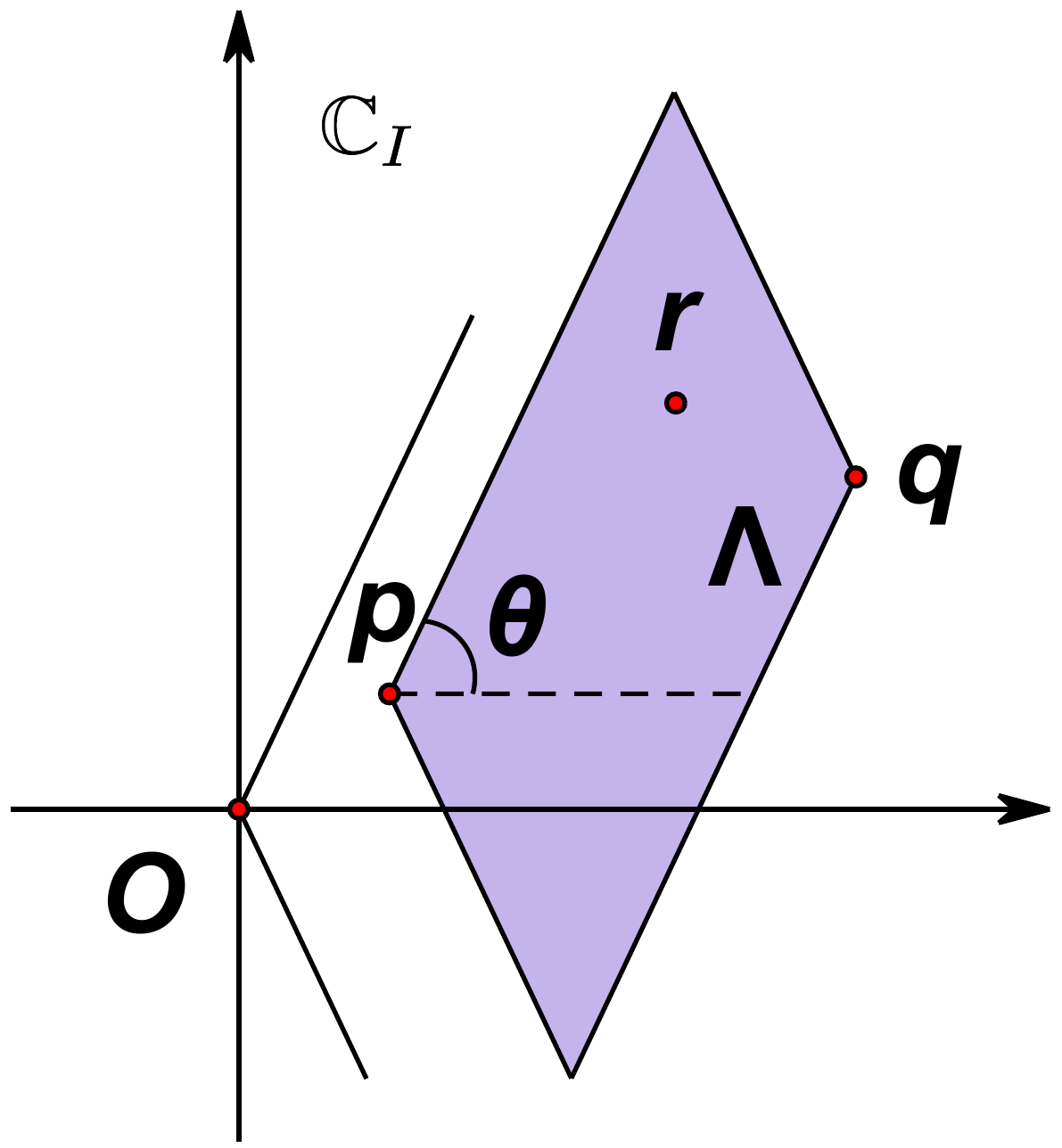}}
        \caption{Schematic diagram of the region $\Lambda$}
        \label{xi3}
    \end{figure}

\begin{remark}
    Note that if $q\in\mathscr{X}_{\theta}(r)$ and $r\in\mathscr{X}_{\theta}(p)$,
    then $q\in \mathscr{X}_{\theta}(p)$.
    Moreover, if $ q \in \mathscr{X}_{\theta}(p) $,
    then there exists $ \Lambda \subset \Sigma_{\theta}$
    such that for all $ r \in \Lambda $,
    we have $ r \in \mathscr{X}_{\theta}(p) $ and
    $ q \in \mathscr{X}_{\theta}(r) $ $(\text{see Figure \ref{xi3}})$.
\end{remark}

\subsubsection{Quaternionic Evolution Systems in Banach Bilateral
$\mathbb{H}$-Algebras}

In this part, we introduce the notion of quaternionic evolution system
of quaternionic Cauchy problem \eqref{kexi-e2}.

\begin{definition}\label{kexi-def4}
    Let $u(q)$ be the solution of \eqref{kexi-e2} in $\Omega_{\mathscr{D}}$
    and $s_{0}=\phi_{J}(z_{0})$ $(z_{0}\in\mathscr{D},J\in\mathbb{S})$.
    A two-parameter operator-valued function
    $\mathcal{U}:\Omega_{\mathscr{D}}\times \Omega_{\mathscr{D}}\rightarrow \mathcal{B}^{R}(X)$
    is called the quaternionic evolution system of \eqref{kexi-e2} if
    $(q,w)\mapsto\mathcal{U}(q,w)$ is right hyperholomorphic with respect to both parameters
    $q$ and $w$ as $q,w\in\mathbb{C}_{K}$ for every $K\in\mathbb{S}$
    and the following relations hold
    \[
        \mathcal{U} \big(\phi_{I}(z),\phi_{I}(z_{0})\big) u\big(\phi_{I}(z_{0})\big)
        =u\big(\phi_{I}(z)\big) \quad
        \big(I\in\mathbb{S},z\in  \mathscr{D},\phi_{I}(z)\in\mathscr{X}_{\theta}(\phi_{I}(z_{0}))\big)
    \]
    and
    \[
        \mathcal{U}\big(\phi_{I}(z_{1}),\phi_{I}(z_{1})\big)=\mathcal{I}, \quad
       \mathcal{U}\big(\phi_{I}(z_{2}),\phi_{I}(z_{1})\big)\odot ^{R}
       \mathcal{U}\big(\phi_{I}(z_{1}),\phi_{I}(z_{0})\big)
       =\mathcal{U}\big(\phi_{I}(z_{2}),\phi_{I}(z_{0})\big),
    \]
    where $z_{1},z_{2}\in\mathscr{D}$,
    $\phi_{I}(z_{1})\in\mathscr{X}_{\theta}\big(\phi_{I}(z_{0})\big)$
    and $\phi_{I}(z_{2})\in\mathscr{X}_{\theta}\big(\phi_{I}(z_{1})\big)$ .
\end{definition}

In Definition \ref{kexi-def4}, we have defined a quaternionic evolution system within the quaternionic Banach bilateral $\mathbb{H}$-algebra. Under the right multiplication operation
$\odot ^{R}$, this quaternionic evolution system satisfies the following important properties.

\begin{theorem}\label{kkkp}
    Let $s_{0}=\phi_{J}(z_{0})$ $(J\in\mathbb{S})$. Then
    $\mathcal{U}(q,w)$  given by Definition \ref{kexi-def4} is bounded
    right linear operator and satisfies the following properties:
    \begin{itemize}
        \item[\rm(i)] There exists a constant $C \geq 1$ such that
        \[
            \big\Vert\mathcal{U}\big(q,\phi_{I}(z_{0})\big)\big\Vert_{\mathcal{B}^{R}(X)}
            \leq \exp\bigg(\int_{\gamma_{[\phi_{I}(z_{0}),q,I]} } \Vert A(\xi)
            \Vert_{\mathcal{B}^{R}(X)} ~\vert\mathrm{d}\xi\vert\bigg),
        \]
        where $q\in\mathscr{X}_{\theta}\big(\phi_{I}(z_{0})\big)$ and $I\in\mathbb{S}$.
        \item[\rm(ii)] $\frac{\partial_{S}^{R}}{\partial q}\mathcal{U}\big(q,\phi_{I}(z_{0})\big)u_{I,0}=\big(A(q)\odot^{R}\mathcal{U}\big(q,\phi_{I}(z_{0})\big)\big)u_{I,0}$, where $q\in\mathscr{X}_{\theta}\big(\phi_{I}(z_{0})\big)$
        and $I\in\mathbb{S}$.
        \item[\rm(iii)] $\frac{\partial_{S}^{R}}{\partial w}\mathcal{U}(q,w)u_{I,0}=\big(-\mathcal{U}(q,w)\odot^{R} A(w)\big)u_{I,0}$, where
        $w=\phi_{I}(z_{0})$, $q\in\mathscr{X}_{\theta}\big(\phi_{I}(z_{0})\big)$, $u_{I,0}\in X$ and $I\in\mathbb{S}$.
    \end{itemize}
\end{theorem}

\begin{proof}
    (i) By \eqref{kexi-jie0}, we obtain that
\begin{equation*}
    \big\|u\big(\phi_{I}(z)\big)\big\|_{X} \leq
        \|u_{I,0}\|_{X} +
        \int_{\gamma_{[\phi_{I}(z_{0}),\phi_{I}(z),I]}} \|A(\xi)\|_{\mathcal{B}^{R}(X)}\|u(\xi)\|_{X} \,\vert \mathrm{d}\xi\vert  .
\end{equation*}
Now,
for $q=\phi_{I}(z)\in \mathscr{X}_{\theta}\big(\phi_{I}(z_{0})\big)$
$(z_{0},z\in\mathscr{D},I\in\mathbb{S})$, let
\[
M(t)=\big\Vert u(\gamma_{[\phi_{I}(z_{0}),q,I]}(t))\big\Vert_{X},\quad
 N(t)=\big\Vert A(\gamma_{[\phi_{I}(z_{0}),q,I]}(t))\big\Vert_{\mathcal{B}^{R}(X)}\cdot \big\vert \gamma_{[\phi_{I}(z_{0}),q,I]}'\big\vert,
 \]
where  $t\in [a,b]$.
Then
\[
    M(b)\leq \|u_{I,0}\|_{X} +\int_{a}^{b} N(\tau)M(\tau)~ \mathrm{d}\tau.
\]
By Lemma 7.1.2 in \cite{kexi-8qq}, we obtain that
\begin{align*}
    &~\big\|\mathcal{U}\big(q,\phi_{I}(z_{0})\big) u_{I,0} \big\|_{X}
    = \big\|u(q)\big\|_{X} \\
    \leq&~
    \|u_{I,0}\|_{X}
    \exp\left( \int_{a}^{b} \big\|A\big(\gamma_{[\phi_{I}(z_{0}),q,I]}(\tau)\big)\big\|_{\mathcal{B}^{R}(X)}\cdot\big\vert\gamma_{[\phi_{I}(z_{0}),q,I]}'\big\vert ~\mathrm{d}\tau  \right) \\
    =&~\|u_{I,0}\|_{X}
    \exp\left( \int_{\gamma_{[\phi_{I}(z_{0}),q,I]}} \|A(\xi)\|_{\mathcal{B}^{R}(X)}~  |\mathrm{d}\xi| \right).
\end{align*}
Thus, $\mathcal{U}\big(q,\phi_{I}(z_{0})\big)$ is bounded and satisfies (i).

(ii) For $q=\phi_{I}(z)\in\mathscr{X}_{\theta}\big(\phi_{I}(z_{0})\big)$
$(z_{0},z\in\mathscr{D},I\in\mathbb{S})$
and $u_{I,0}\in X$,
by \eqref{kexi-jie0} and the right hyperholomorphy of
$\mathcal{U}\big(q,\phi_{I}(z_{0})\big)$,
we obtain that $\mathcal{U}\big(q,\phi_{I}(z_{0})\big)u_{I,0}$
is the unique solution of the integral equation
\begin{equation}\label{kexi-e1.10}
   \mathcal{U}\big(q,\phi_{I}(z_{0})\big)u_{I,0}
   = u_{I,0} + \int_{\gamma_{[\phi_{I}(z_{0}),q,I]}} \big(A(\xi)\odot^{R}
    \mathcal{U}(\xi, \phi_{I}(z_{0}))\big){u_{I,0}}~  \mathrm{d}\xi.
\end{equation}
Taking the right slice derivative of \eqref{kexi-e1.10} with respect to $ q $,
the conclusion (ii) holds.

(iii)
Let $w=\phi_{I}(z_{0})$,
computing the right slice derivative of \eqref{kexi-e1.10} with respect to $w$ we obtain that
\begin{equation*}
    \frac{\partial_{S}^{R}}{\partial w} \mathcal{U}(q, w) {u_{I,0}}
    = -A(w)u_{I,0} + \int_{\gamma_{[w,q,I]}}
    \Big(A(\xi)\odot^{R}\frac{\partial_{S}^{R}}{\partial w} \mathcal{U}(\xi, w) \Big)
   u_{I,0}~
    \mathrm{d}\xi.
\end{equation*}
Let $Y(q)=\frac{\partial_{S}^{R}}{\partial w}\mathcal{U}(q,w){u_{I,0}}$. Then
\begin{equation*}
    Y(q)=-A(w) {u_{I,0}}+\int_{\gamma_{[w,q,I]}} A(\xi) Y (\xi)~\mathrm{d}\xi,
\end{equation*}
which implies that
\[
    \partial_{S}^{R}Y(q)= A(q) Y (q),\quad  Y(w)=-A(w){u_{I,0}}.
\]
By
Theorem \ref{kexi-t4.1}, we obtain that
\[
    Y(q)=\mathcal{U}(q,w)  Y(w)=\mathcal{U}(q,w)  \big(-A(w){u_{I,0}}\big)
    =\big(-\mathcal{U}(q,w)\odot^{R}A(w)\big){u_{I,0}},
\]
i.e.,
\[
\frac{\partial_{S}^{R}}{\partial w} \mathcal{U}(q,w){u_{I,0}}
= \big(-\mathcal{U}(q,w) \odot^{R}A(w)\big){u_{I,0}}.
\]
The completes the proof.
\end{proof}

Through Definition \ref{kexi-def4} of the evolution system in the quaternionic Banach bilateral $\mathbb{H}$-algebra and the related Main Theorem \ref{kkkp},
the integral representation of the solution to the nonhomogeneous equation \eqref{kexi-e4.1} with quaternionic variables in the quaternionic Banach bilateral $\mathbb{H}$-module can be rigorously obtained immediately.
\begin{proposition}\label{kexi-po-01}
    Let $X$ be a Banach bilateral $\mathbb{H}$-module and
    $f\in \mathcal{SH}_{R}(\Omega_{\mathscr{D}},X)$.
    Assume that the uniformly bounded operator-valued function
    $A: \Omega_{\mathscr{D}} \to \mathcal{B}^{R}(X)$
    is right slice regular and $s_{0}=\phi_{J}(z_{0})$ $(z_{0}\in\mathscr{D},J\in\mathbb{S})$.
    If there is a quaternionic evolution
    system $\mathcal{U}(q,w)$ of the
    quaternionic Cauchy problem \eqref{kexi-e2} such that the inclusion
    $u_{I,0}\in\mathcal{D}\big(A(\phi_{I}(z_{0}))\big)$ implies that
    $\mathcal{U}\big(q,\phi_{I}(z_{0})\big) u_{I,0}\in\mathcal{D}\big(A(q)\big)$,
    $q\in \mathscr{X}_{\theta}\big(\phi_{I}(z_{0})\big)$, and
    $\mathcal{U}\big(q,\phi_{I}(z_{0})\big) u_{I,0}$
    admits right slice derivative both in $q$
    and $w$ $\big(w=\phi_{I}(z_{0})\big)$ satisfying
    \begin{align*}
         \frac{\partial_{S}^{R}}{\partial q}\mathcal{U}\big(q,\phi_{I}(z_{0})\big) u_{I,0}
        = &~\big( A(q) \odot^{R} \mathcal{U}(q,\phi_{I}(z_{0})) \big) u_{I,0},\\
        \frac{\partial_{S}^{R}}{\partial w}\mathcal{U}(q,w) u_{I,0}
        = &~\big(-\mathcal{U}(q,w)\odot^{R} A(w) \big) u_{I,0},
    \end{align*}
    then every classical solution $u$ of \eqref{kexi-e4.1},
    with $u_{I,0}\in\mathcal{D}\big(A(\phi_{I}(z_{0}))\big)$,
    restricted to $\big\{p\in\Omega_{\mathscr{D}}:p\in\mathscr{X}_{\theta}(\phi_{I}(z_{0})),I\in\mathbb{S}\big\}$
    is given by
    \begin{equation*}
        u(q)=
        \mathcal{U}\big(q,\phi_{I}(z_{0})\big) u_{I,0}+\int_{\gamma_{[\phi_{I}(z_{0}),q,I]}}
        \mathcal{U}(q,\xi)  f(\xi)~ \mathrm{d}\xi,\quad
         q \in \mathscr{X}_{\theta}\big(\phi_{I}(z_{0})\big).
    \end{equation*}
\end{proposition}

\subsection{Evolution System in Banach Bilateral
$\mathbb{H}$-Algebras under the Hyperbolic-Like Conditions}\label{subsec5}

In this subsection, we will establish the construction of the quaternionic evolution
system linked to the quaternionic Cauchy problem on
a Banach bilateral $\mathbb{H}$-module $X$.
For convenience, we will suppose that
$\mathfrak{D}$ is densely and continuously imbedded in $X$, i.e., $\mathfrak{D}$
forms a dense subspace of $X$ and there exists a constant $C$ such that
\begin{equation*}
    \Vert x\Vert_{X} \leq  C\Vert x\Vert_{\mathfrak{D}}\quad \text{for }x\in \mathfrak{D},
\end{equation*}
and denote it by $\mathfrak{D} \hookrightarrow X$.

\subsubsection{The Homogeneous Evolution Equation under the Hyperbolic-Like Conditions}

For the purposes of the subsequent analysis, the notion of the stable family is hereby presented below.
\begin{definition}\label{kexi-fd0}
    Let $\theta\in(0,\pi/2]$ and $\{A(q)\}_{q\in\Sigma_{\theta}}$ be a family of dense defined
    spherical $\theta$-sectorial operators with vertex $\omega\in\mathbb{R}$.
    Then $\{A(q)\}_{q\in\Sigma_{\theta}} $ is called stable if there are constants
    $M\geqslant 1$ $($called the stability constant$)$ such that
    \begin{equation}\label{kexi-en-0}
        \bigg\Vert \prod_{l=1}^{n}S_{L}^{-1}
        \big(\lambda,A(q_{l})\big)\bigg\Vert_{\mathcal{B}^{R}(X)}
        \leq  \frac{M}{\vert \lambda -\omega\vert^{n}},
        \quad   \lambda\in(\Sigma_{\theta+\pi/2}+\omega).
    \end{equation}
\end{definition}

The following result immediately follows from Theorem 7.1 in \cite{kexi-srs0}.

\begin{theorem}\label{kexi-t5}
    Let $\theta\in(0,\pi/2]$ and $\{A(q)\}_{q\in\Sigma_{\theta}}$ be a family of dense defined
spherical $\theta$-sectorial operators with vertex $\omega$.
    If \eqref{kexi-en-0} holds, then $A(q)$ generates a quaternionic right
    slice regular semigroup
    $\{\mathcal{S}_{q}(w)\}_{w\in\Sigma_{\theta}'}$ such that
    \begin{equation*}
        \left\Vert  \prod_{l=1}^{n}\mathcal{S}_{q_{l}}(p_{l})
        \right\Vert _{\mathcal{B}^{R}(X)}\leq
        M \exp\left\{\omega\sum_{l=1}^{n}\mathrm{Re}(p_{l})\right\} ,
        \quad p_{l}\in \Sigma_{\theta}'
    \end{equation*}
    and
    \begin{equation*}
        \left\Vert  \prod_{l=1}^{n}S_{L}^{-1}\big(\lambda_{l},A(q_{l})\big)
        \right\Vert _{\mathcal{B}^{R}(X)}\leq   M\prod_{l=1}^{n}\vert\lambda_{l}-\omega\vert^{-1},
        \quad   \lambda_{l}\in(\Sigma_{\theta+\pi/2}+\omega).
    \end{equation*}
\end{theorem}

Next, we shall establish a perturbation result, which will be used to prove that the family of perturbed spherical $\theta$-sectorial operators is stable.

\begin{theorem}\label{kexi-qxm0}
    Let $A$ be a spherical $\theta$-sectorial operator with vertex $\omega$ for some
    $\theta\in (0, \pi/2]$. If $P$ is quaternionic linear operator
    such that the following conditions hold:
    \begin{itemize}
\item[\rm(a)] $P\in\mathcal{B}^{R}(X)$ and $\mathcal{D}(A) \subset \mathcal{D}(P)$,
\item[\rm(b)] $R_{\lambda}(P,A) Q_\lambda(A)^{-1} : X \rightarrow X$ for all $\lambda \in \rho_S(A)$,
\item[\rm(c)] $\|R_{\lambda}(P,A) Q_\lambda(A)^{-1}\|_{\mathcal{B}^{R}(X)} < 1$,
for all $\lambda \in \rho_S(A)$,
\end{itemize}
where $Q_{\lambda}(A) =A^{2}+2\mathrm{Re}(\lambda)A+|\lambda|^{2}\mathcal{I}$ and
$R_{\lambda}(P,A):=2\mathrm{Re}(\lambda)P-P^{2}-AP-PA$.
Then $A+P$ is also a spherical $\theta$-sectorial operator with vertex $\omega$.
\end{theorem}

\begin{proof}
    Let $\lambda\in(\omega+\Sigma_{\theta+\pi/2})\subseteq \rho_{S}(A)$. Then
    \begin{align*}
    Q_\lambda(A + P) &= (A + P)^2 - 2\lambda_0(A + P) + |\lambda|^2\mathcal{I} \\
    &=  A^2 - 2\mathrm{Re}(\lambda)A + |\lambda|^2\mathcal{I} -
    2\mathrm{Re}(\lambda)P + P^2 + AP + PT\\
    &= Q_{\lambda}(A)-R_{\lambda}(P,A)
    =\big(\mathcal{I}-R_{\lambda}(P,A)Q_{\lambda}(A)^{-1}\big)Q_{\lambda}(A).
\end{align*}
Since $\|R_{\lambda}(P,A) Q_\lambda(A)^{-1}\|_{\mathcal{B}^{R}(X)} < 1$,
applying Lemma 6.6 (Neumann series) in \cite{kexi-jishu-06}, then we obtain that
$\big(\mathcal{I}-R_{\lambda}(P,A)Q_{\lambda}(A)^{-1}\big)$ is invertible.
Thus, it follows that $Q_{\lambda}(A+P)$ is invertible,
which implies that $\lambda\in \rho_{S}(A+P)$.
Therefore, $A+P$ is also a spherical $\theta$-sectorial operator with vertex
$\omega$. This completes the proof.
\end{proof}

\begin{corollary}\label{kexi-wd0}
    Let $\theta\in(0,\pi/2]$ and $\{A(q)\}_{q\in\Sigma_{\theta}}$ be a family of
    spherical $\theta$-sectorial operators with vertex $\omega$ on $X$.
    If $\{A(q)\}_{q\in \Sigma_{\theta}}$ is stable and
    $ P(q)$ satisfies the conditions {\rm (a), (b), (c)}
    in Theorem \ref{kexi-qxm0} for each $q\in\Sigma_{\theta}$. In addition,
    if the following
    condition holds:
    \begin{itemize}
    \item[\rm(d)]
$\|(\overline{\lambda}\mathcal{I}-A(q)-P(q))\|_{\mathcal{B}^{R}(X)}\|(\overline{\lambda}\mathcal{I}-A(q))^{-1}\|_{\mathcal{B}^{R}(X)}\leq C$ for each $q\in \Sigma_{\theta}$,
where $C>0$ is a constant,
    \end{itemize}
    then $\big\{A(q)+P(q)\big\}_{q\in \Sigma_{\theta}}$ is a stable family.
\end{corollary}
\begin{proof}
    By Theorem \ref{kexi-qxm0}, $A(q)+P(q)$ is a
    spherical $\theta$-sectorial operator with vertex $\omega$
    for each $q\in \Sigma_{\theta}$. On the other hand,
    for $\lambda\in(\omega+\Sigma_{\theta+\pi/2})\subseteq \rho_{S}(A)$,
    it follows from Proposition 2.10 in \cite{kexi-xxx9} that
    \[
        S_{L}^{-1}\big(\lambda,A(q)+P(q)\big)
        =\bigg[\sum_{m\geq 0}Q_{\lambda}\big(A(q)\big)^{-1}\big(R_{\lambda}(P(q),A(q))Q_{\lambda}(A(q))^{-1}\big)^{m}
        \bigg] \big(\overline{\lambda}\mathcal{I}-A(q)-P(q)\big).
    \]
    Hence
    \begin{align*}
        &~\big\|S_{L}^{-1}\big(\lambda,A(q)+P(q)\big)\big\|_{\mathcal{B}^{R}(X)}\\
        \leq &~
        \bigg\|Q_{\lambda}\big(A(q)\big)^{-1}\big(\overline{\lambda}\mathcal{I}-A(q)\big)
        \big(\overline{\lambda}\mathcal{I}-A(q)\big)^{-1}\\
        &~\cdot
        \bigg[\sum_{m\geq 0}\big(R_{\lambda}(P(q),A(q))Q_{\lambda}(A(q))^{-1}\big)^{m}\bigg]
        \big(\overline{\lambda}\mathcal{I}-A(q)-P(q)\big)\bigg\|_{\mathcal{B}^{R}(X)}\\
        \leq &~
        \big\|S_{L}^{-1}\big(\lambda,A(q)\big)\big\|_{\mathcal{B}^{R}(X)}\cdot
        \bigg\|\sum_{m\geq 0}\big(R_{\lambda}(P(q),A(q))Q_{\lambda}(A(q))^{-1}\big)^{m}\bigg\|_{\mathcal{B}^{R}(X)}\\
        &~\cdot \big\|\overline{\lambda}\mathcal{I}-A(q)-P(q)\big\|_{\mathcal{B}^{R}(X)}\cdot
        \big\|(\overline{\lambda}\mathcal{I}-A(q))^{-1}\big\|_{\mathcal{B}^{R}(X)}.
    \end{align*}
  Since the condition (c) holds, then there exists a constant $C_{1}>0$ such that
    \[
        \bigg\|\sum_{m\geq 0}\Big(R_{\lambda}\big(P(q),A(q)\big)Q_{\lambda}\big(A(q)\big)^{-1}\Big)^{m}\bigg\|_{\mathcal{B}^{R}(X)}\leq C_{1}.
    \]
    Therefore,
    \[
        \big\|S_{L}^{-1}\big(\lambda,A(q)+P(q)\big)\big\|_{\mathcal{B}^{R}(X)}\leq \frac{CC_{1}M}{|\lambda-\omega|},
    \]
    where $M$ is the stability constant,
    which implies that
    \[
        \Big\|\prod_{l=1}^{n} S_{L}^{-1}(\lambda,A(q_{l})+P(q_{l}))\Big\|_{\mathcal{B}^{R}(X)}
        \leq \frac{(CC_{1}M)^{n}}{|\lambda-\omega|^{n}}:=\frac{\widetilde{M}}{|\lambda-\omega|^{n}}.
    \]
    Thus, $\big\{A(q)+P(q)\big\}_{q\in\Sigma_{\theta}}$ is stable family. This completes the proof.
\end{proof}

\begin{definition}\label{kexi-def2}
    Let $X$ be a Banach bilateral $\mathbb{H}$-module and
    let $\{A(q)\}_{q\in\Sigma_{\theta}}$ is a stable family.
    The notation $\mathscr{P}\big(A(q)\big)$ is used to
    denote the set of quaternionic operators that satisfy
    the conditions {\rm (a), (b), (c)} in Theorem \ref{kexi-qxm0} and
    {\rm (d)} in Corollary \ref{kexi-wd0}.
\end{definition}

To further establish the properties of stable quaternionic generators,
for a Banach bilateral $\mathbb{H}$-module,
we introduce the notion of $A$-admissible
subspace.

\begin{definition}
    Let $\{\mathcal{S} (q)\}_{q\in\Sigma_{\theta}'}$ be a
    quaternionic right slice regular semigroup
    with generator $A$.
    A subspace $Y$ of $X$ is called $A$-admissible if it is an invariant
    subspace of $\mathcal{S}(q)$,
    and the restriction of $\mathcal{S}(q)$ to $Y$ is a quaternionic
    right slice regular semigroup in $Y$.
\end{definition}

In the following part, we shall focus on the quaternionic
Cauchy problem as follows:
\begin{equation}\label{kexi-fd5}
    \begin{cases}\displaystyle
        \partial_{S}^{R}u (q)= A (q)u(q),\\
      u(s_{0})=u_{0},
    \end{cases}
\end{equation}
where $s_{0},q\in \Sigma_{\theta}$.
In the following discussion, we will always assume that $\theta\in (0,\pi/2]$,
$s_{0}= \phi_J(z_0) $ for some $ J \in \mathbb{S} $ and $ z_0 \in \mathscr{D} $.

To find the quaternionic evolution system and the solutions for \eqref{kexi-fd5},
we need to assume some additional conditions.
For each $q\in\Sigma_{\theta}$, let $A(q)$ be a spherical $\theta$-sectorial operator on $X$.
We shall adopt the following assumptions.
    \begin{itemize}
        \item[$(\mathfrak{H}_{1})$] $\{A(q)\}_{q\in\Sigma_{\theta}} $ is a
        stable family with stability constant $M>0$.
        \item[$(\mathfrak{H}_{2})$] $\mathfrak{D}$ is $A(q)$-admissible for
        $q\in \Sigma_{\theta}$ and $\widetilde{A}(q)$ is
        the restriction of $A(q)$ in $\mathfrak{D}$ such that
        $\big\{\widetilde{A}(q)\big\}_{q\in\Sigma_{\theta}}$
        is a stable family in $\mathfrak{D}$ with stability
        constant $\widetilde{M}$.
        \item[$(\mathfrak{H}_{3})$]
        For $q\in \Sigma_{\theta}$,
        $\mathfrak{D}\subset\mathcal{D}\big(A(q)\big)$,
        $q\mapsto A(q)\in\mathcal{B}^{R}(\mathfrak{D},X)$
        is right slice regular and $A(q)$ is uniformly bounded.
    \end{itemize}

The following result demonstrates that if
$\{A(q)\}_{q\in\Sigma_{\theta}}$ satisfies assumptions
$(\mathfrak{H}_{1})$-$(\mathfrak{H}_{3})$, then one can associate a unique quaternionic evolution system
$\mathcal{U}(q,w)$ of the quaternionic Cauchy problem \eqref{kexi-fd5}.

\begin{theorem}\label{ht1}
    Let the family $\{A(q)\}_{q\in \Sigma_{\theta}}$
    satisfy the stipulated hypotheses
    $(\mathfrak{H}_{1})$-$(\mathfrak{H}_{3})$. Then
    there exists a unique quaternionic evolution
    system $ \mathcal{U}( q, w)$,
    $(q,w)\in \Sigma_{\theta} \times \Sigma_{\theta}$,
    satisfying the following conditions$:$
    \begin{itemize}
        \item[$(\mathfrak{E}_{1})$]
        $\Vert\mathcal{U}(q,w)\Vert_{\mathcal{B}^{R}(X)} \leq  M e^{\omega\cdot\mathrm{Re}(q-w)}$, where $w=\phi_{I}(z_{0})$,
        $q\in\mathscr{X}_{\theta}\big(\phi_{I}(z_{0})\big)$ and $I\in\mathbb{S}$.

        \item[$(\mathfrak{E}_{2})$]
        $\frac{\partial_{S}^{R}}{\partial q}\mathcal{U}\big(q,\phi_{I}(z_{0})\big) u_{I,0} =\big(A(q)\odot^{R}\mathcal{U}(q,\phi_{I}(z_{0}))\big)u_{I,0}$, where
        $u_{I,0}\in \mathfrak{D}$, $q \in \mathscr{X}_{\theta}\big(\phi_{I}(z_{0})\big)$ and $I\in\mathbb{S}$.

        \item[$(\mathfrak{E}_{3})$]
        $\frac{\partial_{S}^{R}}{\partial w}\mathcal{U}(q,w) u_{I,0}=\big(-\mathcal{U}(q,w)\odot^{R} A(w)\big) u_{I,0}$, where $u_{I,0}\in \mathfrak{D}$, $w=\phi_{I}(z_{0})$,
        $q\in \mathscr{X}_{\theta}\big(\phi_{I}(z_{0})\big)$ and $I\in\mathbb{S}$.
    \end{itemize}
\end{theorem}

\begin{proof}
    For given $I \in \mathbb{S}$ and
    $p\in \mathscr{X}_{\theta}\big(\phi_{I}(z_{0})\big)$,
    let $\gamma_{[\phi_{I}(z_{0}),p,I]}:[a,b]\rightarrow \mathbb{C}_{I}$
    be a $C^{1}$-path from $\phi_{I}(z_{0})$ to $p$,
    satisfying
    $\gamma_{[\phi_{I}(z_{0}),p,I]}(a)=\phi_{I}(z_{0})$ and $\gamma_{[\phi_{I}(z_{0}),p,I]}(b)=p$.
    For convenience, we will assume $0\leq a\leq b$ and denote
    $\gamma_{[\phi_{I}(z_{0}),p,I]}$ by $\gamma$.

    We initiate the process by approximating the family of quaternionic operator
    $\{A_{n}(q)\}$ $(q\in\gamma ([a,b]))$
    using piecewise constant families
    $\{A(q)\}$ $(q\in\gamma ([a,b])) $, for $n=1,2,\ldots$, which are defined as follows.

    Since $(\mathfrak{H}_{1})$ is satisfied,
    applying Theorem 7.1 in \cite{kexi-srs0}, we obtain that $A(q)$ generates
    a quaternionic right slice regular semigroup
    $\big\{\mathcal{S}_{q}(w)\big\}_{w\in\Sigma_{\theta}'}$.
    Let $t_{\kappa}^{n}=\frac{\kappa(b-a)}{n}$,
    where $\kappa=0,1,\ldots,n$. The operator family $\{A_{n}(q)\}_{q\in\gamma ([a,b])}$ is then
    constructed by setting:
    \begin{equation*}
        \left\{\begin{aligned}
            &A_{n}\big(\gamma (t)\big)=A\big(\gamma (t_{\kappa}^{n})\big),\quad\text{for }
            t_{\kappa}^{n}\leq t< t_{\kappa+1}^{n},\quad \kappa=0,1,\ldots,n-1,\\
            &A_{n}\big(\gamma (b)\big)=A\big(\gamma (b)\big).
        \end{aligned}\right.
    \end{equation*}
    Since $q\mapsto A(q)$ is right slice regular,
    then $q\mapsto A(q)$ is continuous, and hence
    \begin{equation}\label{kexi-wo7}
        \Vert A(q)-A_{n}(q)\Vert_{\mathcal{B}^{R}(\mathfrak{D},X)}
        \rightarrow 0\quad\text{as}\quad n\rightarrow\infty
    \end{equation}
    uniformly in $t\in [a,b]$, where $q=\gamma (t)\in \gamma ([a,b])$.
    According to the definition of $A_{n}(q)$ and the assumptions of
    the theorem, it follows that for $n\geqslant 1$, $\big\{A_{n}(q)\big\} $ $(q\in\gamma ([a,b])) $
    is a stable family restricted to $\gamma$
    in $X$ with stability constant $M$ and $\big\{\widetilde{A}_{n}(q)\big\} $
    $\big(q\in\gamma ([a,b])\big)$
    is a stable family restricted to $\gamma$ in $\mathfrak{D}$ with stability
    constant $\widetilde{M}$.

    Now, let us define a two-parameter family of quaternionic operators
    $\mathcal{U}_{n}(q,w)$ with respect to $\gamma$ as follows:
    \begin{equation}\label{kexi-iq0}
        \mathcal{U}_{n}\big(\gamma (t),\gamma(\tau)\big)
        =\left\{\begin{aligned}
             \mathcal{S}_{\gamma(t_{\varsigma }^{n})}\big(\gamma(t)-\gamma(\tau)\big)
            \qquad~~~ &\text{for} \quad t^{n}_{\varsigma}
            \leq \tau \leq t \leq t_{\varsigma +1}^{n},\\
             \mathcal{S}_{\gamma(t_{\kappa }^{n})}\big(\gamma(t)-\gamma(t_{\kappa }^{n})\big)\odot^{R}
            \bigg[&{}^{\odot^{R}}\prod_{\varsigma =\ell+1}^{\kappa -1}
            \mathcal{S}_{\gamma(t_{\varsigma }^{n})} \left(\gamma\left(\frac{b-a}{n}\right) \right)\bigg]
            \odot^{R}\mathcal{S}_{\gamma(t_{\ell}^{n})}\big(\gamma(t_{\ell+1}^{n})-\gamma(\tau)\big)\\
            &\text{for} \quad \kappa >\ell, t_{\kappa }^{n}
            \leq t \leq t_{\kappa +1}^{n},t_{\ell}^{n}\leq \tau
            \leq t_{\ell+1}^{n} ,
        \end{aligned}\right.
    \end{equation}
    where $\big\{\mathcal{S}_{q}(w)\big\}_{w\in\Sigma_{\theta}'}$
    is the quaternionic right slice regular group generated by
    $A(q)$, $ a\leq\tau\leq t\leq b $, and
    \begin{align*}
       &~ {}^{\odot^{R}}\prod_{\varsigma =\ell+1}^{\kappa -1}
            \mathcal{S}_{\gamma(t_{\varsigma }^{n})}
            \left(\gamma\left(\frac{b-a}{n}\right) \right)\\
            =&~\mathcal{S}_{\gamma(t_{\ell+1}^{n})}\left(\gamma\left(\frac{b-a}{n}\right) \right)
            \odot^{R}
            \mathcal{S}_{\gamma(t_{\ell+2}^{n})}\left(\gamma\left(\frac{b-a}{n}\right) \right)
            \odot^{R}\cdots\odot^{R}
            \mathcal{S}_{\gamma(t_{\kappa-1}^{n})}\left(\gamma\left(\frac{b-a}{n}\right) \right).
    \end{align*}
    It is evident that $\mathcal{U}(q,w)$ is a quaternionic evolution system
    restricted to $\gamma \big([a,b]\big)$ such that
    \begin{equation}\label{kexi-wo10}
        \mathcal{U}_{n}(q,q)=\mathcal{I},
        ~~\mathcal{U}_{n}(q,w)=\mathcal{U}_{n}(q,r)\odot^{R}\mathcal{U}_{n}(r,w) ,
    \end{equation}
    where $q,r,w\in\gamma\big([a,b]\big)$ with $q\in\mathscr{X}_{\theta}(r)$,
    $r\in\mathscr{X}_{\theta}(w)$.
    In the same time, we note that
        \begin{equation}\label{kexi-wo11}
        (q,w)\mapsto \mathcal{U}_{n}(q,w)\text{ is right hyperholomorphic on }
        \gamma ([a,b])\text{ respect to $q$ and $w$}.
    \end{equation}
    Therefore, for any $I\in\mathbb{S}$ and
    $ p \in \mathscr{X}_{\theta}\big(\phi_{I}(z_{0})\big)$,
    repeating the above steps yields a quaternionic evolution system over
    the entire $\Sigma_{\theta}\times\Sigma_{\theta}$.
    By Theorem \ref{kexi-t5}, we obtain that
    \begin{equation}\label{kexi-wo9}
        \Vert \mathcal{U}_{n}(q,w)\Vert_{\mathcal{B}^{R}(X)} \leq
        Me^{\omega\cdot\mathrm{Re}(q-w)} \quad\big(w=\phi_{I}(z_{0})\big),
    \end{equation}
    where $q \in \mathscr{X}_{\theta} \big(\phi_{I} (z_{0})\big)$.
    On the other hand, from the assumption $(\mathfrak{H}_{2})$, it follows that
    \begin{equation*}
        \mathcal{U}_{n}(q,w)
        \mathfrak{D}\subset \mathfrak{D}, \quad q \in\mathscr{X}_{\theta} \big(\phi_{I} (z_{0})\big).
    \end{equation*}
     Now,
    since $\mathfrak{D}\subset\mathcal{D}\big(A(q)\big)$ for
    $q\in \mathscr{X}_{\theta} \big(\phi_{I} (z_{0})\big)$, then by
    the definition of $\mathcal{U}_{n}(q,w)$
    and Remarks \ref{kexi-re0-1} and
    \ref{kexi-rem1}, for $u_{I,0}\in \mathfrak{D}$, we obtain that
    \begin{align}\label{kexi-m0}
        \frac{\partial_{S}^{R}}{\partial q}\mathcal{U}_{n}(q,w)u_{I,0}=
         \big(A_{n}(q)\odot^{R}\mathcal{U}_{n}(q,w)\big) u_{I,0} \quad \big(q\in \mathscr{X}_{\theta}(w)\big)
    \end{align}
    for $q=\gamma(t)$, $t\neq t_{\varsigma}^{n},~\varsigma=0,1,\ldots,n$, and
\begin{align}\label{kexi-m1}
        \frac{\partial_{S}^{R}}{\partial w}\mathcal{U}_{n}(q,w)u_{I,0}=
        \big(-\mathcal{U}_{n}(q,w)\odot^{R}A_{n}(w)\big)u_{I,0}.
    \end{align}
    From Theorem \ref{kexi-t5} and hypotheses $(\mathfrak{H}_{2})$, it follows that
    \begin{equation}\label{kexi-qi4}
        \Vert\mathcal{U}_{n}(q,w)\Vert_{\mathcal{B}^{R}(\mathfrak{D},X)}
         \leq  \widetilde{M}e^{ \omega \cdot \mathrm{Re}(q-w) }\quad
        \text{for }  q \in \mathscr{X}_{\theta}(w).
    \end{equation}
    For $u_{I,0}\in \mathfrak{D}$,
    consider the map $r\mapsto \big(\mathcal{U}_{n}(q,r)\odot^{R}\mathcal{U}_{m}(r,w)\big)u_{I,0}$,
    where $r\in \gamma\big([a,b]\big) $ with
    $q\in  \mathscr{X}_{\theta}(r)$ and $r\in \mathscr{X}_{\theta}(w)$. By
    \eqref{kexi-m0} and \eqref{kexi-m1}, we obtain that except for a finite number of values of
    $r=\gamma(\tau')$ $(a\leq \tau'\leq t)$, this map admits right slice
    derivative at $r$ and by Remark 2.1.13 in \cite{kexi-co3},
    we obtain that
    \begin{align*}
        &~\mathcal{U}_{n}(q,w)u_{I,0}-\mathcal{U}_{m}(q,w)u_{I,0}\\
        =&~-\int_{\gamma_{[w,q,I]}}
         \frac{\partial^{R}_{S}}{\partial r}\big(\mathcal{U}_{n}(q,r)
         \odot^{R}\mathcal{U}_{m}(r,w)\big)u_{I,0}~\mathrm{d}r  \\
        =&~\int_{\gamma_{[w,q,I]}}
        \big(\mathcal{U}_{n}(q,r)\odot^{R}\big(A_{n}(r)-A_{m}(r)\big)\odot^{R}
        \mathcal{U}_{m}(r,w)\big)
        u_{I,0}~\mathrm{d}r ,
    \end{align*}
    which implies that
    \begin{align}\label{kexi-wo0}
        \begin{split}
        &~\big\Vert\mathcal{U}_{n}(q,w) u_{I,0}-\mathcal{U}_{m}(q,w) u_{I,0}\big\Vert_{X} \\
        \leq&~  M\widetilde{M} e^{\omega\cdot\mathrm{Re}(q-w)}
        \Vert u_{I,0}\Vert_{\mathfrak{D}}\int_{\gamma_{[w,q,I]}}
        \Vert A_{n}(r)-A_{m}(r)\Vert_{\mathcal{B}^{R}(\mathfrak{D},X)}~|\mathrm{d}r| .
        \end{split}
    \end{align}
    According to \eqref{kexi-wo7} and \eqref{kexi-wo0}, $\mathcal{U}_{n}(q,w) u_{I,0}$ converges in $X$,
    uniformly on
    $(q,w)\in \Sigma_{\theta} \times \Sigma_{\theta}$,
    $q\in\mathscr{X}_{\theta}(w)$, as $n\rightarrow\infty$.
    Since $\mathfrak{D}$ is dense in $X$, by
    the convergence of $\mathcal{U}_{n}(q,w) u_{I,0}$
    and \eqref{kexi-wo9}, we obtain that $\mathcal{U}_{n}(q,w)$ converges
    strongly in $X$, uniformly on
    $(q,w)\in \Sigma_{\theta} \times \Sigma_{\theta}$, $q\in\mathscr{X}_{\theta}(w)$. Let
    \begin{equation}\label{kexi-qi5}
        \mathcal{U}(q,w) v=\lim\limits_{n\rightarrow\infty}
        \mathcal{U}_{n}(q,w) v\quad \text{for }
        v\in X ,q\in\mathscr{X}_{\theta}(w) .
    \end{equation}
    Then, by \eqref{kexi-wo10}, \eqref{kexi-wo11} and
    \eqref{kexi-wo9}, $(\mathfrak{E}_{1})$ is satisfied.

    Next, we will show that $(\mathfrak{E}_{2})$ and $(\mathfrak{E}_{3})$.
    Let
    $r\mapsto\big(\mathcal{U}_{n}(q,r)\odot^{R} \mathcal{S}_{\xi}(r-w)\big) u_{I,0}$
    be a function for
    $u_{I,0}\in \mathfrak{D}$, where $r\in\mathscr{X}_{\theta}(w)$ and
    $q\in\mathscr{X}_{\theta}(r)$.
    Note that this function admits right slice derivative except for a finite number of values of
    $r$ and hence
    \begin{align*}
        \begin{split}
        \mathcal{U}_{n}(q,w) u_{I,0}-\mathcal{S}_{\xi}(q-w) u_{I,0}
        =&~-\int_{\gamma_{[w,q,I]}} \frac{\partial^{R}_{S}}{\partial r}
        \big(\mathcal{U}_{n}(q,r)\odot^{R}
        \mathcal{S}_{\xi}(r-w)\big) u_{I,0}~\mathrm{d}r\\
            =&~\int_{\gamma_{[w,q,I]}}
            \big(\mathcal{U}_{n}(q,r)\odot^{R}
            \big(A_{n}(r)-A(\xi)\big)\odot^{R}
            \mathcal{S}_{\xi}(r-w)\big)
             u_{I,0}
            ~\mathrm{d}r ,
        \end{split}
    \end{align*}
    which implies that
    \begin{equation*}
        \big\Vert \mathcal{U}_{n}(q,w) u_{I,0}-\mathcal{S}_{\xi}(q-w) u_{I,0}\big\Vert_{X} \leq
        M\widetilde{M}e^{\omega\cdot\mathrm{Re}(q-w)}
        \Vert u_{I,0}\Vert_{\mathfrak{D}}\int_{\gamma_{[w,q,I]}}
        \Vert A_{n}(r)-A(\xi)\Vert_{\mathcal{B}^{R}(\mathfrak{D},X)}~|\mathrm{d}r|.
    \end{equation*}
    Let $n\rightarrow\infty$, we have
    \begin{equation}\label{kexi-wo1}
        \Vert\mathcal{U}(q,w) u_{I,0}-\mathcal{S}_{\xi}(q-w) u_{I,0}\Vert_{X} \leq  M\widetilde{M}e^{\omega\cdot\mathrm{Re}(q-w)}
        \Vert u_{I,0}\Vert_{\mathfrak{D}}\int_{\gamma_{[w,q,I]}}
        \Vert A(r)-A(\xi)\Vert_{\mathcal{B}^{R}(\mathfrak{D},X)}~|\mathrm{d}r| .
    \end{equation}
    Choosing $\xi=w$ in \eqref{kexi-wo1}, multiply it by $(q-w)^{-1}$ and
    letting $q\rightarrow w$, by the right slice regularity of $q\mapsto A(q)$ in the
    $\mathcal{B}^{R}(\mathfrak{D},X)$-norm, we obtain that
    \begin{equation}\label{kexi-wo2}
        \underset{q\in\mathbb{C}_{ I },q\rightarrow w}{\lim\sup}\Vert
        \mathcal{U}(q,w) u_{I,0}-\mathcal{S}_{w}(q-w) u_{I,0}\Vert_{X}(q-w)^{-1} =0.
    \end{equation}
    Note that
    $\mathcal{S}_{w}(q-w) u_{I,0}$ admits right slice derivative at $q$,
    by \eqref{kexi-wo2},
    $\mathcal{U}(q,w) u_{I,0}$ admits right slice
    derivative and they have the same derivatives at $q$,
    which indicates that $(\mathfrak{E}_{2})$ is satisfied.

    According to $(\mathfrak{E}_{2})$ and the right slice regularity of $\mathcal{U}(q,w)$
    in $X$, we obtain that
    \begin{align*}
             \frac{\partial_{S}^{R}}{\partial w}\mathcal{U}(q,w) u_{I,0}
            =&~\Big(\frac{\partial_{S}^{R}}{\partial w} \mathcal{U}(q,w)\Big)u_{I,0}\\
            =&~\Big(\lim\limits_{h\in\mathbb{C}_{I},h\rightarrow 0 }
            \big(\mathcal{U}(q,w+h) -\mathcal{U}(q,w) \big)h^{-1}\Big)u_{I,0}\\
            =&~\Big(\lim\limits_{h\in\mathbb{C}_{I},h\rightarrow 0 }\mathcal{U}(q,w+h)\odot^{R}
            \big(\mathcal{I}-\mathcal{U}(w+h,w)  \big) h^{-1}\Big)u_{I,0} \\
            =&\big(-\mathcal{U}(q,w)\odot^{R}A(w)\big)u_{I,0},
    \end{align*}
    which yields $\mathcal{U}(q,w)$ satisfies $(\mathfrak{E}_{3})$.

    Finally, we need to show that $\mathcal{U}(q,w)$, $q \in \mathscr{X}_{\theta}(w)$, is
    unique.
    Suppose $W(q,w)$ is another quaternionic
    evolution system satisfying $(\mathfrak{E}_{1})$-$(\mathfrak{E}_{3})$.
    Let $u_{I,0}\in \mathfrak{D}$, and let
    $r\mapsto \big(W(q,r)\odot^{R}\mathcal{U}_{n}(r,w)\big) u_{I,0}$ be a function,
    where $r\in \mathscr{X}_{\theta}(w)$ and $q\in \mathscr{X}_{\theta}(r)$.
    Since $W(q,w)$ adheres to the condition $(\mathfrak{E}_{3})$, and in light of the construction of
    $\mathcal{U}_{n}(q,w)$, we deduce that this function
    admits right slice derivative
    everywhere except at a finite number of points $r$.
    Integrating the derivative subsequently yields the following expression
    \begin{equation*}
        \begin{aligned}
        &~W(q,w) u_{I,0}-\mathcal{U}_{n}(q,w) u_{I,0}\\
        =&~\int_{\gamma_{[w,q,I]}}\big(W(q,r) \odot^{R}
        \big(A(r)-A_{n}(r)\big) \odot^{R}\mathcal{U}_{n}(r,w)\big) u_{I,0}~\mathrm{d}r,
        \end{aligned}
    \end{equation*}
    thus,
    \begin{equation}\label{wo6}
        \begin{aligned}
        &~\Vert W(q,w) u_{I,0}-\mathcal{U}_{n}(q,w) u_{I,0}\Vert_{X}\\
        \leq&~
        M\widetilde{M}e^{\omega\cdot\mathrm{Re}(q-w)}\Vert u_{I,0}\Vert_{\mathfrak{D}}
        \int_{\gamma_{[w,q,I]}}
        \Vert A(r)-A_{n}(r)\Vert_{\mathcal{B}^{R}(\mathfrak{D},X)}~|\mathrm{d}r|.
        \end{aligned}
    \end{equation}
    Taking the limit as $n\rightarrow\infty$ in \eqref{wo6},
    and utilizing the result established in
    \eqref{kexi-wo7}, it is concluded that the identity
    $W(q,w) u_{I,0}=\mathcal{U}(q,w) u_{I,0}$
    holds for all initial values
    $u_{I,0}\in \mathfrak{D}$ and $q\in\mathscr{X}_{\theta}(w)$.
    Given that $\mathfrak{D}$ constitutes a dense subspace of $X$,
    and since both $\mathcal{U}(q,w)$ and $W(q,w)$ comply with
    $(\mathfrak{E}_{1})$,
    we definitively establish the equality: $\mathcal{U}(q,w)=W(q,w)$. This completes the proof.
\end{proof}

To obtain the solutions of \eqref{kexi-fd5} via quaternionic evolution system, we will introduce a notion of solutions
of \eqref{kexi-fd5} called $\mathfrak{D}$-valued solution.

\begin{definition}
    A continuous function $u\in C( \Sigma_{\theta};\mathfrak{D})$ is called
    a $\mathfrak{D}$-valued solution of the quaternionic
    Cauchy problem \eqref{kexi-fd5} if
    $u: \Sigma_{\theta}  \to X$
    is right slice regular and satisfies \eqref{kexi-fd5} in $X$.
\end{definition}

 By the Definition \ref{kexi-def4} and Theorem \ref{ht1},
the following result can be obtained immediately.

\begin{theorem}\label{tht0}
    Let $\{A(q)\}_{q\in \Omega_{\mathscr{D}}}$
    satisfy the stipulated hypotheses
    $(\mathfrak{H}_{1})$-$(\mathfrak{H}_{3})$.
    Suppose that  the quaternionic Cauchy problem \eqref{kexi-fd5} has a
    $\mathfrak{D}$-valued solution $u$, then $u$ is unique and
    \begin{equation*}
        u(q)=\mathcal{U}\big(q,\phi_{I}(z_{0})\big) u_{I,0} \quad
        \big(q\in \mathscr{X}_{\theta}(\phi_{I}(z_{0})),I\in\mathbb{S}\big),
    \end{equation*}
    where
    $\mathcal{U}(q,w)$ is the quaternionic evolution system given by Theorem \ref{ht1}.
\end{theorem}

Now, we shall establish the conditions to guarantee $\mathcal{U}(q,w) u_{I,0}$
is the $\mathfrak{D}$-valued solution
of \eqref{kexi-fd5}.

\begin{theorem}\label{kexi-ht2}
    Let $\{A(q)\}_{q\in \Omega_{\mathscr{D}}}$
    satisfy the stipulated hypotheses
    $(\mathfrak{H}_{1})$-$(\mathfrak{H}_{3})$.
    Let $\mathcal{U}(q,w)$ be the quaternionic evolution system given in
    Theorem \ref{ht1}.
         Suppose that the following conditions hold:
    \begin{itemize}
        \item[$(\mathfrak{E}_{4})$] $\mathcal{U}\big(q,\phi_{I}(z_{0})\big) \mathfrak{D} \subset\mathfrak{D}$, where $I\in\mathbb{S}$ and
        $q\in\mathscr{X}_{\theta}\big(\phi_{I}(z_{0})\big)$.
        \item[$(\mathfrak{E}_{5})$] For any $u_{0} \in\mathfrak{D}$,
        $\mathcal{U}\big(q,\phi_{I}(z_{0})\big) u_{I,0}$
        is right hyperholomorphic in $\mathfrak{D}$, where  $q\in\mathscr{X}_{\theta}\big(\phi_{I}(z_{0})\big)$.
    \end{itemize}
    Then for each $u_{0}\in \mathfrak{D}$,
    $\mathcal{U}\big(q,\phi_{I}(z_{0})\big) u_{I,0}$ is the
    unique $\mathfrak{D}$-valued solution of the Cauchy problem \eqref{kexi-fd5}, which is restricted to
    $\big\{p\in\Sigma_{\theta}:p\in\mathscr{X}_{\theta}(\phi_{I}(z_{0})),I\in\mathbb{S}\big\}$.
\end{theorem}

\begin{proof}
    By Theorem \ref{tht0},
    the uniqueness of $\mathfrak{D}$-valued solutions of \eqref{kexi-fd5}
    is immediate. It remains to show that
    if $q\in \mathscr{X}_{\theta}\big(\phi_{I}(z_{0})\big)$ and $u_{0}\in \mathfrak{D}$, then
    $u(q)= \mathcal{U}\big(q,\phi_{I}(z_{0})\big) u_{I,0}$ is a
    $\mathfrak{D}$-valued solution of \eqref{kexi-fd5}.
    In view of the assumptions of theorem, we obtain that
    $u(q)\in \mathfrak{D}$.
    Moreover, $u(q)$ is right hyperholomorphic
    in the $\mathfrak{D}$-norm for $z\in \mathscr{D} $ and each $I\in \mathbb{S}$.

    Next, we need to show that $u$ satisfies \eqref{kexi-fd5}.
    Since $u(q)=\mathcal{U}\big(q,\phi_{I}(z_{0})\big) u_{I,0}\in\mathfrak{D}$ for
    $z \in\mathscr{D}$ and $I\in\mathbb{S}$,
    it follows from Theorem \ref{ht1} that
    \begin{align}\label{kexi-iq3}
        \begin{split}
             \frac{\partial^{R}_{S}}{\partial q}
            \mathcal{U}\big(q,\phi_{I}(z_{0})\big) u_{I,0}
            =&~\Big(\frac{\partial_{S}^{R}}{\partial q}\mathcal{U}\big(q,\phi_{I}(z_{0})\big)\Big)u_{I,0}\\
            =&~\Big(\lim\limits_{h\in\mathbb{C}_{I},h\rightarrow 0 }
             \big(\mathcal{U}(q+h,\phi_{I}(z_{0})) -\mathcal{U}(q,\phi_{I}(z_{0}))  \big) h^{-1}\Big)u_{I,0}\\
            =&~\Big(\lim\limits_{h\in\mathbb{C}_{I},h\rightarrow 0 }\big(
            (\mathcal{U}(q+h,q) -\mathcal{I} )\odot^{R}\mathcal{U}(q,\phi_{I}(z_{0})) \big)h^{-1}\Big)u_{I,0}\\
            =&~\big(A(q)\odot^{R} \mathcal{U} (q,\phi_{I}(z_{0}))\big)u_{I,0}.
        \end{split}
    \end{align}
    Since $q\mapsto \mathcal{U}\big(q,\phi_{I}(z_{0})\big) u_{I,0}$
    is right hyperholomorphic in the
    $\mathfrak{D}$-norm and $q\mapsto A(q)$ is right slice regular, then
    the right side of \eqref{kexi-iq3} is right hyperholomorphic in $X$ for each $I\in\mathbb{S}$.
    Thus, $\mathcal{U}\big(q,\phi_{I}(z_{0})\big) u_{I,0}$ is right slice regular in $X$
    for $q\in\{p\in\Sigma_{\theta}:p\in\mathscr{X}_{\theta}(\phi_{I}(z_{0})),I\in\mathbb{S}\}$.
    By \eqref{kexi-iq3}, we obtain that
    \[
        \frac{\partial_{S}^{R}}{\partial q} \mathcal{U}\big(q,\phi_{I}(z_{0})\big) u_{I,0}
        =\big(A(q)\odot^{R}\mathcal{U}\big(q,\phi_{I}(z_{0})\big)\big) u_{I,0}\quad \text{for } q\in \mathscr{X}_{\theta}\big(\phi_{I}(z_{0})\big).
    \]
    Therefore, from Theorem \ref{kexi-t4.2}, it follows that
    $\mathcal{U}\big(q,\phi_{I}(z_{0})\big)$ can be extended to $\mathfrak{D}$-valued
    solution on $\Sigma_{\theta}$.
    This completes the proof.
\end{proof}

\begin{remark}
    It follows from Theorem \ref{kexi-ht2} that if the quaternionic evolution system
    $\mathcal{U}(q,w)$ given
    by Theorem \ref{ht1} satisfies the conditions $(\mathfrak{E}_{4})$ and $(\mathfrak{E}_{5})$, then
    the quaternionic Cauchy problem \eqref{kexi-fd5} has a unique $\mathfrak{D}$-valued solution
    restricted to $\big\{p\in\Sigma_{\theta}:p\in\mathscr{X}_{\theta}(\phi_{I}(z_{0})),I\in\mathbb{S}\big\}$
    given by
    $\mathcal{U}\big(\phi_{I}(z),\phi_{I}(z_{0})\big) u_{I,0}$ for any $u_{0}\in \mathfrak{D}$,
    $z \in\mathscr{D}$
    and $I\in\mathbb{S}$.
    To obtain a quaternionic evolution system $\mathcal{U}(q,w)$ that satisfies
    $(\mathfrak{E}_{1})$-$(\mathfrak{E}_{5})$, we will consider the following condition:
    \begin{description}
        \item[$(\mathfrak{H}_{2})^{+}$]  There is a family
        $\{Q(q)\}_{q\in\Sigma_{\theta}} $ of isomorphisms of
        $\mathfrak{D}$ onto $X$ such that
        $q\mapsto Q(q)\in\mathcal{B}^{R}(\mathfrak{D},X)$ and $q\mapsto Q(q)^{-1}\in\mathcal{B}^{R}(X,\mathfrak{D})$ are right slice regular.
        Moreover,
        \begin{equation*}
            (Q\odot ^{R}A\odot^{R}Q^{-1})(q)=A(q)+B(q),
        \end{equation*}
        where $B(q)\in\mathscr{P}(A(q))$ and
        $Q(q)$ satisfies $\partial_{S}^{R}Q(q) \odot^{R}Q ^{-1} (q)\in\mathscr{P}\big(A(q)\big)$
        $($see Definition \ref{kexi-def2}$)$
        and $q\mapsto B(q)\in\mathcal{B}^{R}(X)$ is right slice regular.
    \end{description}
\end{remark}

According to Theorem \ref{kexi-qxm0}, we can immediately obtain the following conclusion.

\begin{lemma}\label{iq4}
 The conditions $(\mathfrak{H}_{1})$ and $(\mathfrak{H}_{2})^{+}$ suffice to ensure
    $(\mathfrak{H}_{2})$.
\end{lemma}

By applying Gr{\"o}nwall's inequality,
we immediately obtain the following result.

\begin{lemma}\label{iq5}
    Let $\mathcal{U}(q,w)$ be a quaternionic evolution system
    such that $\big\Vert\mathcal{U}\big(\phi_{I}(z),\phi_{I}(z_{0})\big)\big\Vert_{\mathcal{B}^{R}(X)} \leq  M$
    for $z \in \mathscr{D} $ and $I\in\mathbb{S}$.
    If $q\mapsto W(q) $ is right slice regular and $W(q)$ is bounded
    for every $q\in \Omega_{\mathscr{D}}$, then for $u_{0}\in X$,
    there exists a unique
    family of bounded
    right hyperholomorphic operators-valued function
    $\big(\phi_{I}(z),\phi_{I}(z_{0})\big)\mapsto V\big(\phi_{I}(z),\phi_{I}(z_{0})\big)$
    $(z \in \mathscr{D})$ on $X$ such that
    \begin{align*}
         &~V\big(\phi_{I}(z),\phi_{I}(z_{0})\big) u_{I,0}\\
        = &~ \mathcal{U}\big(\phi_{I}(z),\phi_{I}(z_{0})\big) u_{I,0}
         +\int_{\gamma_{[\phi_{I}(z_{0}),\phi_{I}(z),I]}}\big(V(\phi_{I}(z),r)
        \odot^{R}W(r)\odot^{R}\mathcal{U}(r,\phi_{I}(z_{0}))\big) u_{I,0}~\mathrm{d}r .
    \end{align*}
\end{lemma}

Next, we are ready to present the main result of this subsection.
\begin{theorem}\label{kexi-fd7}
    Let $\{A(q)\}_{q\in\Sigma_{\theta}}$ satisfy the conditions $(\mathfrak{H}_{1})$,
    $(\mathfrak{H}_{2})^{+}$
    and $(\mathfrak{H}_{3})$.
    Then there exists a unique
    quaternionic evolution system
    $\mathcal{U}(q,w)$ in $X$
    satisfying $(\mathfrak{E}_{1})$-$(\mathfrak{E}_{5})$.
\end{theorem}

\begin{proof}
    From Lemma \ref{iq4}, it follows that $(\mathfrak{H}_{1})$-$(\mathfrak{H}_{3})$ are satisfied.
    According to Theorem \ref{ht1}, there exists a unique quaternionic
    evolution system $\mathcal{U}(q,w)$ satisfying $(\mathfrak{E}_{1})$-$(\mathfrak{E}_{3})$.
    Let $u_{0}\in \mathfrak{D}$ and denote the right slice derivative of
    $Q(q)u_{I,0}$ by $\dot{Q}(q)u_{I,0}$ $\big(q=\phi_{I}(z),z\in\mathscr{D}\big)$, and set
    \begin{equation*}
        \mathcal{C}(q)=\dot{Q}(q)\odot^{R}Q(q)^{-1}.
    \end{equation*}
    Then $\mathcal{C}(q) $ is right slice regular.
    From Lemma \ref{iq5}, it follows that for $u_{0}\in X$,
    there exists a unique bounded solution $V\big(\phi_{I}(z),\phi_{I}(z_{0})\big)$ of
    the integral equation
    \begin{align}\label{qi6}
        \begin{split}
         &~V\big(q,\phi_{I}(z_{0})\big) u_{I,0}\\
        = &~\mathcal{U}\big(q,\phi_{I}(z_{0})\big)u_{I,0}
        +\int_{\gamma_{[\phi_{I}(z_{0}),q,I]}}\big(V(q,r)\odot^{R} \big(B(r)
        +\mathcal{C}(r)\big)\odot ^{R}\mathcal{U}(r,\phi_{I}(z_{0}))\big) u_{I,0}~\mathrm{d}r,
        \end{split}
    \end{align}
    where $q=\phi_{I}(z)\in\mathscr{X}_{\theta}\big(\phi_{I}(z_{0})\big)$.
   Since $\mathcal{U}\big(q,\phi_{I}(z_{0})\big) \mathfrak{D}\subset \mathfrak{D}$,
    by the construction of \eqref{kexi-iq6},
    then $\mathcal{U}(q,w)$ satisfies $(\mathfrak{E}_{4})$.
    Next, according to the right hyperholomorphy of $V\big(q,\phi_{I}(z_{0})\big) u_{I,0}$
    and the properties of $Q\big(\phi_{I}(z_{0})\big)$ and $Q(q)^{-1}$,
    we obtain that $\mathcal{U}\big(q,\phi_{I}(z_{0})\big) u_{I,0}$
    is right hyperholomorphic
    in $\mathfrak{D}$ for $q\in  \mathscr{X}_{\theta}\big(\phi_{I}(z_{0})\big) $.
    Thus, $\mathcal{U}(q,w)$ satisfies  $(\mathfrak{E}_{5})$.

    In the following, we will show that
    \begin{equation}\label{kexi-iq6}
        \mathcal{U}\big(q,\phi_{I}(z_{0})\big)=Q(q)^{-1}\odot^{R}V\big(q,\phi_{I}(z_{0})\big)\odot^{R}Q(\phi_{I}
        \big(z_{0})\big).
    \end{equation}
    It follows from $(\mathfrak{H}_{3}) $
    that for each $u_{0}\in X$, $Q(q)^{-1} u_{I,0}$ $\big(q=\phi_{I}(z),z\in\mathscr{D}\big)$ admits right slice
    derivative in $\mathfrak{D}$ and
    \begin{equation}\label{iq7}
       \partial_{S}^{R}\big(Q(q)^{-1} u_{I,0}\big)
        =\big(-Q(q)^{-1}\odot^{R}\dot{Q}(q)\odot^{R}Q(q)^{-1}\big)u_{I,0}.
    \end{equation}
    Let
    \begin{equation*}
        Q(q,r)=\mathcal{U}(q,r)\odot ^{R}Q(r)^{-1} \quad \big(q\in\mathscr{X}_{\theta}(r),
        r\in\mathscr{X}_{\theta} (\phi_{I}(z_{0}) )\big).
    \end{equation*}
    According to \eqref{iq7} and $(\mathfrak{E}_{3})$, for each $u_{I,0}\in X$,
    it follows that $r\mapsto Q(q,r) u_{I,0}$ admits right slice derivative in $X$ and
    \begin{align*}
         &\frac{\partial^{R}_{S}}{\partial r}Q(q,r) u_{I,0}\\
         =&~-\big(\mathcal{U}(q,r)\odot ^{R}A(r)\odot^{R}Q(r)^{-1}\big) u_{I,0}\\
         &~- \big(\mathcal{U}(q,r)\odot ^{R}Q(r)^{-1}\odot^{R}\dot{Q}(r)\odot^{R}Q(r)^{-1}\big) u_{I,0}\\
        =&~\big(-\mathcal{U}(q,r)\odot ^{R}A(r)\odot ^{R}Q(r)^{-1}\big) u_{I,0}
        -\big(Q(q,r)\odot ^{R}\mathcal{C}(r)\big) u_{I,0}.
    \end{align*}
    On the other hand, for every $u_{I,0}\in\mathfrak{D}$,
    by the condition $(\mathfrak{H}_{2})^{+} $, we have
    \[
        \big(A(r)\odot^{R}Q(r)^{-1}\big) u_{I,0}
        =\big(Q(r)^{-1}\odot^{R}\big(A(r)+B(r)\big)\big)
         u_{I,0}.
    \]
    Thus, for $u_{I,0}\in \mathfrak{D}$, we obtain that
    \begin{align}\label{kexi-iq9}
        \begin{split}
          \frac{\partial_{S}^{R}}{\partial r}Q(q,r) u_{I,0}
        =&~\big(-\mathcal{U}(q,r)\odot ^{R}Q(r)^{-1}\odot ^{R}\big(A(r)+B(r)
           +\mathcal{C}(r)\big)\big) u_{I,0}\\
        =&~\big(-Q(q,r)\odot^{R}\big(A(r)+B(r)+\mathcal{C}(r)\big)\big) u_{I,0}.
        \end{split}
    \end{align}
    Let $\mathcal{U}_{n}(q,r)$
    be the quaternionic operator constructed in \eqref{kexi-iq0},
    it follows from \eqref{kexi-m0} that
    \begin{equation}\label{kexi-iq8}
        \frac{\partial_{S}^{R}}{\partial r}\mathcal{U}_{n}\big(r,\phi_{I}(z_{0})\big) u_{I,0}
        =\big(A_{n}(r)\odot^{R}\mathcal{U}_{n}(r,\phi_{I}(z_{0}))\big) u_{I,0}\quad
        \big( u_{I,0}\in \mathfrak{D}, I\in\mathbb{S} \big),
    \end{equation}
    for all $r\in\mathscr{X}_{\theta}\big(\phi_{I}(z_{0})\big)$
    except for a finite number of values of $r$.
    According to \eqref{kexi-iq9} and \eqref{kexi-iq8}, we obtain that
    \begin{align}\label{kexi-qi0}
        \begin{split}
        &~\frac{\partial_{S}^{R}}{\partial r}\big(Q(q,r)\odot^{R}\mathcal{U}_{n}\big(r,\phi_{I}(z_{0})\big)\big) u_{I,0}\\
        =&~\big(-Q(q,r)\odot^{R}\big(A(r)+B(r)
         +\mathcal{C}(r)-A_{n}(r)\big)
        \odot^{R}\mathcal{U}_{n}\big(r,\phi_{I}(z_{0})\big)\big) u_{I,0}.
    \end{split}
    \end{align}
    Integrating \eqref{kexi-qi0} along $\gamma_{[\phi_{I}(z_{0}),q,I]}$ yields
    \begin{align}\label{kexi-qi1}
        \begin{split}
        &\big(Q(t)^{-1}\odot^{R}\mathcal{U}_{n}\big(q,\phi_{I}(z_{0})\big)\big) u_{I,0}-
        Q\big(q,\phi_{I}(z_{0})\big) u_{I,0}\\
        =&-\int_{\gamma_{[\phi_{I}(z_{0}),q,I]}}
        \big(Q(q,r) \odot^{R}\big(A(r)+B(r)+\mathcal{C}(r)-A_{n}(r)\big)
        \odot^{R}\mathcal{U}_{n}\big(r,\phi_{I}(z_{0})\big)\big) u_{I,0}~\mathrm{d}r.
        \end{split}
    \end{align}
    It follows from $(\mathfrak{E}_{1})$ and \eqref{kexi-qi4} that
    \begin{align}\label{kexi-qi2}
        &~\left\Vert\int_{\gamma_{[\phi_{I}(z_{0}),q,I]}}
        \big(Q(q,r)\odot^{R}\big(A(r)-A_{n}(r)\big)\odot^{R}\mathcal{U}_{n}\big(r,\phi_{I}(z_{0})\big) \big)u_{I,0}~\mathrm{d}r\right\Vert_{X} \\
        \leq &~ M\widetilde{M}e^{\omega\cdot\mathrm{Re}(q-\phi_{I}(z_{0}))}M_{1}M_{2},
    \end{align}
    where
    \begin{align*}
        &M_{1}=\sup\limits_{r\in\gamma_{[\phi_{I}(z_{0}),q,I]}}\Vert Q^{-1}(r)\Vert_{\mathcal{B}^{R}(X,\mathfrak{D})}
        \Vert u_{I,0}\Vert_{\mathfrak{D}},\\
        &M_{2}=\int_{\gamma_{[\phi_{I}(z_{0}),q,I]}}\Vert A(r)-A_{n}(r)\Vert_{\mathcal{B}^{R}(\mathfrak{D},X)}~|\mathrm{d}r|.
    \end{align*}
    Let $n\rightarrow\infty$ in \eqref{kexi-qi1}, according to
    \eqref{kexi-qi2} and \eqref{kexi-qi5},
    for $u_{I,0}\in \mathfrak{D}$, it follows that
    \begin{align}\label{kexi-qi3}
        \begin{split}
        &~\big(Q(q)^{-1}\odot^{R}\mathcal{U}(q,\phi_{I}(z_{0}))\big) u_{I,0}-Q\big(q,\phi_{I}(z_{0})\big) u_{I,0}\\
        =&~-\int_{\gamma_{[\phi_{I}(z_{0}),q,I]}}
        \big(Q(q,r)\odot^{R}\big(B(r)+\mathcal{C}(r)\big)
        \odot^{R}\mathcal{U}(r,\phi_{I}(z_{0}))\big)u_{I,0}~\mathrm{d}r.
    \end{split}
    \end{align}
    Note that all the quaternionic operators in \eqref{kexi-qi3} are bounded in $X$ and
    $\mathfrak{D}$ is dense in $X$, then
    \eqref{kexi-qi3} holds for every $u_{I,0}\in X$, which implies that
    \begin{equation}\label{kexi-qi7}
        \begin{aligned}
        Q\big(q,\phi_{I}(z_{0})\big) u_{I,0}= &~ \big(Q(q)^{-1}\odot^{R}\mathcal{U}\big(q,\phi_{I}(z_{0})\big)\big) u_{I,0} \\
         &~+\int_{\gamma_{[\phi_{I}(z_{0}),q,I]}}
         \big(Q(q,r)
        \odot^{R}(B(r)+\mathcal{C}(r))\odot^{R}\mathcal{U}
        \big(r,\phi_{I}(z_{0})\big)\big) u_{I,0}~\mathrm{d}r.
    \end{aligned}
    \end{equation}
    On the other hand, left-multiplying \eqref{qi6} by $Q(q)^{-1}$ gives
    \begin{align}\label{kexi-qi8}
        \begin{split}
        &~\big(Q(q)^{-1}\odot^{R}V(q,\phi_{I}(z_{0}))\big) u_{I,0}
        = \big(Q(q)^{-1}\odot^{R}\mathcal{U}(q,\phi_{I}(z_{0}))\big) u_{I,0}\\
        &~+\int_{\gamma_{[\phi_{I}(z_{0}),q,I]}}
        \big(Q(q)^{-1}\odot^{R}V(q,r)
        \odot^{R} (B(r)+\mathcal{C}(r))\odot^{R}\mathcal{U}(r,\phi_{I}(z_{0}))\big) u_{I,0}~\mathrm{d}r.
        \end{split}
    \end{align}
    In view of \eqref{kexi-qi7} and the uniqueness of the solution of \eqref{kexi-qi8}, we conclude that
    \[
        \mathcal{U}\big(q,\phi_{I}(z_{0})\big)\odot ^{R}Q(\phi_{I})^{-1}
        =Q\big(q,\phi_{I}(z_{0})\big)=Q(q)^{-1}\odot^{R}V\big(q,\phi_{I}(z_{0})\big),
    \]
    which implies that \eqref{kexi-iq6} holds. This completes the proof.
\end{proof}

From Theorems \ref{kexi-t4.2}, \ref{kexi-ht2} and \ref{kexi-fd7}, the following corollary is immediate.

\begin{corollary}
    Let $\{A(q)\}_{q\in\Sigma_{\theta}}$ satisfy the conditions $(\mathfrak{H}_{1})$, $(\mathfrak{H}_{2})^{+}$
    and $(\mathfrak{H}_{3})$.
    Then for every $u_{0}\in \mathfrak{D}$, the quaternionic Cauchy problem
    \begin{equation*}
        \left\{\begin{aligned}
            &\partial_{S}^{R}u(q)=A(q)u(q), \\
            &u(s_{0})=u_{0},
        \end{aligned}\right.\qquad s_{0},q\in\Sigma_{\theta}
    \end{equation*}
    has a unique $\mathfrak{D}$-valued solution on $\Sigma_{\theta}$.
\end{corollary}

A special case of Theorem \ref{kexi-fd7} in which $\mathcal{D}\big(A(q)\big)=\mathfrak{D}$ is independent
of $q\in\Sigma_{\theta}$ will be established in the following theorem.

\begin{theorem}\label{ia4}
    Let $\{A(q)\}_{q\in\Sigma_{\theta}}$ be a family of stable spherical $\theta$-sectorial operators
    with vertex $\omega$ on $X$.
    If $\mathcal{D}\big(A(q)\big)=\mathfrak{D}$
    is independent of $q$ and for $u_{0}\in \mathfrak{D}$, $q\mapsto A(q)$
    is right slice regular, and there exists $\lambda >\omega$ such that
    \[
        \partial_{S}^{R} A(q) \odot^{R}\big(\lambda \mathcal{I} -A(q)\big)^{-1}\in \mathscr{P}
        \big(A(q)\big),
    \]
    then there exists a unique quaternionic evolution system $\mathcal{U}(q,w)$
    satisfying conditions $(\mathfrak{E}_{1})$-$(\mathfrak{E}_{5})$,
    where $\mathfrak{D}$ is equipped
    with the norm $\Vert\cdot\Vert_{\mathfrak{D}}$ given by
    \begin{equation}\label{ia2}
        \Vert u_{0}\Vert_{\mathfrak{D}}=\Vert u_{0} \Vert_{X}
        +\Vert A(0) u_{0}\Vert_{X}\quad \big(u_{0}\in \mathfrak{D}=\mathcal{D}(A)\big).
    \end{equation}
\end{theorem}

\begin{proof}
    By Theorem \ref{kexi-fd7},
    it is sufficient to show that $\{A(q)\}_{q\in\Sigma_{\theta}}$ satisfies the conditions
    $(\mathfrak{H}_{1})$, $(\mathfrak{H}_{2})^{+}$ and $(\mathfrak{H}_{3})$.
    The conditions $(\mathfrak{H}_{1})$ and $(\mathfrak{H}_{3})$ are given directly in our assumptions.

    On the other hand, by the assumption, for $ \lambda $, the quaternionic operator
    $Q(q)= \big(\lambda\mathcal{I}-A(q)\big)$
    is an isomorphism from $\mathfrak{D}$ onto $X$, and since
    $q\mapsto A(q)$ is right slice regular,
    by Remark \ref{kexi-re0-1}, $A(q)u_{0}$ is right slice regular. Therefore,
    $Q(q) u_{0}$ is right slice regular in $X$ for each
    $u_{0}\in \mathfrak{D}$. Furthermore, since
    \[
        (Q \odot^{R}A \odot^{R}Q ^{-1})(q)=A(q),
    \]
    i.e., we take $B(q)\equiv 0$ in condition $(\mathfrak{H}_{2})^{+}$,
    which implies that $\{A(q)\}_{q\in\Sigma_{\theta}}$ satisfies the condition $(\mathfrak{H}_{2})^{+}$. This completes the proof.
\end{proof}

\subsubsection{The Inhomogeneous Evolution Equation under the Hyperbolic-Like Conditions}

In this part, let us consider the following inhomogeneous quaternionic Cauchy problem
\begin{equation}\label{ia0}
    \left\{\begin{aligned}
    &\partial_{S}^{R}u(q)=A(q)u(q)+f(q),\\
    &u(s_{0})=u_{0}
    \end{aligned}\right.\qquad s_{0},q\in\Sigma_{\theta}
\end{equation}
such that the conditions $(\mathfrak{H}_{1})$-$(\mathfrak{H}_{3})$ hold.

\begin{theorem}\label{iq1}
    Let $\{A(q)\}_{q\in\Omega_{\mathscr{D}}}$ satisfy
    the assumptions $(\mathfrak{H}_{1})$-$(\mathfrak{H}_{3})$  and let
    $f\in \mathcal{SH}_{R}\big(\Omega_{\mathscr{D}},X\big)$.
    If the quaternionic Cauchy problem \eqref{ia0}
    has a $\mathfrak{D}$-valued solution
    $u$, then this solution is unique and its restriction to
    $\big\{p\in\Sigma_{\theta}:p\in\mathscr{X}_{\theta}(\phi_{I}(z_{0})),I\in\mathbb{S}\big\}$
    is given by
    \begin{equation}\label{ia1}
        u\big(\phi_{I}(z)\big)= \mathcal{U}\big(\phi_{I}(z),\phi_{I}(z_{0})\big) u_{I,0}+
        \int_{\gamma_{[\phi_{I}(z_{0}),\phi_{I}(z),I]}} \mathcal{U}\big(\phi_{I}(z),r\big) f(r)~\mathrm{d}r,
    \end{equation}
    where $z\in\mathscr{D}$, $\phi_{I}(z)\in \mathscr{X}_{\theta}\big(\phi_{I}(z_{0})\big)$ and
    $\mathcal{U}(q,w)$ is the quaternionic evolution system given by Theorem \ref{ht1}.
\end{theorem}

\begin{proof}
    Suppose that $\mathcal{U}_{n}(q,w)$ is
    the quaternionic evolution
    system given in \eqref{kexi-iq0} and $u$ is a $\mathfrak{D}$-solution of \eqref{ia0}.
    By the properties of $\mathcal{U}_n(q,w)$ and $u$, it follows
    that the function
    $r\mapsto\mathcal{U}_n(q,r) u(r)  $ $\big(q\in \mathscr{X}_{\theta}(r)\big)$
    is right hyperholomorphic in $X$ except for a finite number of values of $r$ and
    \begin{equation}\label{tht1}
        \begin{aligned}
         \frac{\partial_{S}^{R}}{\partial r} \mathcal{U}_{n}(q,r) u(r)
        =  \big(-\mathcal{U}_{n}(q,r )\odot^{R}A_{n}(r)\big) u(r)
         +\big(\mathcal{U}_{n}(q,r)\odot^{R} A(r)\big) u(r)
         +\mathcal{U}_{n}(q,r) f(r).
        \end{aligned}
    \end{equation}
    Integrating \eqref{tht1} along $\gamma_{[\phi_{I}(z_{0}),q,I]}$, we have
    \begin{equation*}
        \begin{aligned}
        u(q)= &~\mathcal{U}_{n}\big(q,\phi_{I}(z_{0})\big) u_{I,0}+\int_{\gamma_{[\phi_{I}(z_{0}),q,I]}}
        \mathcal{U}_{n}(q,r) f(r)~\mathrm{d}r\\
        &~+\int_{\gamma_{[\phi_{I}(z_{0}),q,I]}}
        \big(\mathcal{U}_{n}(q,r) \odot^{R}
        (A(r)-A_{n}(r) )\big) u(r)~\mathrm{d}r.
        \end{aligned}
    \end{equation*}
    Letting $C=\max_{r\in\gamma_{[\phi_{I}(z_{0}),q,I]}}\Vert u(r)\Vert_{\mathfrak{D}}$,
    it follows from \eqref{kexi-wo9} that
    \begin{equation}\label{tht2}
        \begin{aligned}
            &~\bigg\Vert u(q)-\mathcal{U}_{n}\big(q,\phi_{I}(z_{0})\big)
            u_{I,0}-\int_{\gamma_{[\phi_{I}(z_{0}),q,I]}}\mathcal{U}_{n}
            (q,r) f(r)~\mathrm{d}r\bigg\Vert_{X}\\
            \leq   &~ M e^{\omega \cdot\mathrm{Re} (q-\phi_{I}(z_{0}) )}C\int_{\gamma_{[\phi_{I}(z_{0}),q,I]}}
            \big\Vert A(r)-A_{n}(r)\big\Vert_{\mathcal{B}^{R}(\mathfrak{D},X)}~|\mathrm{d}r| .
        \end{aligned}
    \end{equation}
    Letting $n\rightarrow\infty$ in \eqref{tht2}, by \eqref{kexi-wo7} and \eqref{kexi-qi5}.
    Then \eqref{ia1} holds for all $ I\in\mathbb{S}$. The
    uniqueness of $u$ follows from Theorem \ref{ht1} and \eqref{ia1}. This completes the proof.
\end{proof}

Now we will establish the existence and uniqueness of $\mathfrak{D}$-valued solutions for the quaternionic inhomogeneous Cauchy problem.
\begin{theorem}\label{ia8}
    Let $\{A(q)\}_{q\in\Omega_{\mathscr{D}}}$
    satisfy the assumptions
    $(\mathfrak{H}_{1})$-$(\mathfrak{H}_{3})$ and $(\mathfrak{E}_{4})$-$(\mathfrak{E}_{5})$ given
    by Theorem \ref{kexi-ht2}. If
    $f \in \mathcal{SH}_{R}\big(\Omega_{\mathscr{D}}, \mathfrak{D})$, then for every $u_{0}\in \mathfrak{D}$,
    the quaternionic Cauchy problem
    \eqref{ia0} has a unique $\mathfrak{D}$-valued solution $u$ and
    its restriction to
    $\big\{p\in\Sigma_{\theta}:p\in\mathscr{X}_{\theta}(\phi_{I}(z_{0})),I\in\mathbb{S}\big\}$
    is given by \eqref{ia1}.
\end{theorem}

\begin{proof}
    From Theorem \ref{kexi-ht2}, it follows that
    $\mathcal{U}\big(q,\phi_{I}(z_{0})\big) u_{I,0}$ is a $\mathfrak{D}$-valued solution
    of the homogeneous case \eqref{kexi-fd5} restricted to
    $\big\{p\in\Sigma_{\theta}:p\in\mathscr{X}_{\theta}(\phi_{I}(z_{0})),I\in\mathbb{S}\big\}$.
    To show that
    \eqref{ia1}
    is a $\mathfrak{D}$-valued solution of \eqref{ia0}
    restricted to
    $\big\{p\in\Sigma_{\theta}:p\in\mathscr{X}_{\theta}(\phi_{I}(z_{0})),I\in\mathbb{S}\big\}$,
    we need to prove that
    \[
        w\big(\phi_{I}(z)\big)=\int_{\gamma_{[\phi_{I}(z_{0}),\phi_{I}(z),I]}} \mathcal{U}\big(\phi_{I}(z),r\big) f(r)~\mathrm{d}r\quad
        \big(\phi_{I}(z)\in\mathscr{X}_{\theta}(\phi_{I}(z_{0})),z\in\mathscr{D},I\in\mathbb{S}\big)
    \]
    is a $\mathfrak{D}$-valued solution of \eqref{ia0}
    restricted to
    $\big\{p\in\Sigma_{\theta}:p\in\mathscr{X}_{\theta}(\phi_{I}(z_{0})),I\in\mathbb{S}\big\}$
    with the initial value $w(s_{0})=u_{0}=0$.
    By the assumptions on $f$ and $\mathsf{(\mathfrak{E}_{4})}$, it follows that
    $w\big(\phi_{I}(z)\big)\in \mathfrak{D}$
    for $z\in \mathscr{D} $. On the other hand, since
    $\mathsf{(\mathfrak{E}_{5})}$ is satisfied, then
    $r \mapsto \mathcal{U}(q,r) f(r) $ $\big(q\in \mathscr{X}_{\theta}(r)\big)$
    is right hyperholomorphic in $\mathfrak{D}$, which implies
    that $q\mapsto w(q)$ is right hyperholomorphic
    in $\mathfrak{D}$, where $q=\phi_{I}(z)$, and that
     $r \mapsto \big(A(q) \odot^{R}\mathcal{U}(q,r)\big)f(r) $
    is right hyperholomorphic in $X$ for $z\in\mathscr{D} $, $I\in\mathbb{S}$.
    Therefore, by the right hyperholomorphy of
    $
        r \mapsto \big(A(q)\odot^{R} \mathcal{U}(q,r)\big) f(\phi_{I}(z_{1}))
    $
    we obtain that $w(q)$ is right hyperholomorphic in $X$ and
    \[
       \partial_{S}^{R}w(q)=A(q)
        w(q)+f(q),\quad
        \forall  q\in \mathscr{X}_{\theta}\big(\phi_{I}(z_{0})\big)
    \]
    holds in $X$.
    Finally, by Theorem \ref{iq1}, it follows that the
    $\mathfrak{D}$-valued solutions of \eqref{ia0} is unique.
    This completes the proof.
\end{proof}

In the sequel, we will establish the existence and uniqueness of
the classical solution of the quaternionic
Cauchy problem \eqref{ia0}.
\begin{theorem}
    Let $\{A(q)\}_{q\in\Omega_{\mathscr{D}}}$ be a family of stable spherical $\theta$-sectorial operators
    with vertex $\omega$ on $X$
    such that $ \mathcal{D}\big(A(q)\big)$
    is independent of $q$ and
    $u_{0}\in \mathcal{D}\big(A(q)\big)$. Let $q\mapsto A(q) $
    be right slice regular and there exists $\lambda >\omega$ large enough such that
    \[
        \partial_{S}^{R} A(q)\odot^{R} \big(\lambda \mathcal{I} -A(q)\big)^{-1}\in \mathscr{P}
        \big(A(q)\big).
    \]
    If $f\in \mathcal{SH}_{R}\big(\Omega_{\mathscr{D}},X\big)$, then for
    every $u_{0}\in\mathcal{D}\big(A(q)\big)$, the quaternionic Cauchy problem \eqref{ia0} has a unique classical solution $u$ whose restriction to
    $\big\{p\in\Sigma_{\theta}(\phi_{0}):p\in\mathscr{X}_{\theta}(\phi_{I}(z_{0})),I\in\mathbb{S}\big\}$ is
    given by
    \begin{equation}\label{ia5}
        u\big(\phi_{I}(z)\big)= \mathcal{U}\big(\phi_{I}(z),\phi_{I}(z_{0})\big)  u_{I,0}+
        \int_{\gamma_{[\phi_{I}(z_{0}),\phi_{I}(z),I]}} \mathcal{U}\big(\phi_{I}(z),r\big) f(r)~\mathrm{d}r,
    \end{equation}
    where $z \in\mathscr{D}$, $I\in\mathbb{S}$.
\end{theorem}

\begin{proof}
    Let $\mathcal{D}\big(A(q)\big)$ endow with the graph norm of $A(0)$ given by \eqref{ia2} and denote this Banach bilateral $\mathbb{H}$-module by $\mathfrak{D}$.
    According to the assumptions for $\lambda>\omega$, we have
    $Q(q)= \big(\lambda\mathcal{I}-A(q)\big)$
    is an isomorphism of $\mathfrak{D}$ onto $X$ such that
    $Q(q) u_{I,0}$ $\big(q=\phi_{I}(z),z\in\mathscr{D}\big)$ is right
    hyperholomorphic in $X$ for
    $u_{0}\in \mathfrak{D}$.
    Let
    $\dot{Q}(q) u_{I,0}$ be the right slice derivative of $Q(q) u_{I,0}$,
    it follows that $\dot{Q}(q)\in\mathcal{B}^{R}(\mathfrak{D},X)$ and
    $\big\{\Vert\dot{Q}(q)\Vert_{\mathcal{B}^{R}(\mathfrak{D},X)}\big\}_{q\in\Sigma_{\theta}}$ is uniformly bounded. It follows from Theorem
    \ref{ia4} that $\mathcal{U}\big(q ,\phi_{I}(z_{0})\big) u_{I,0}$ is the $\mathfrak{D}$-valued solution of the
    quaternionic homogeneous Cauchy problem
    \eqref{kexi-fd5} restricted to
    $\big\{p\in\Sigma_{\theta}(\phi_{0}):p\in\mathscr{X}_{\theta}(\phi_{I}(z_{0})),I\in\mathbb{S}\big\}$. To show that \eqref{ia5} is a classical solution of \eqref{ia0} restricted to
    $\big\{p\in\Sigma_{\theta}(\phi_{0}):p\in\mathscr{X}_{\theta}(\phi_{I}(z_{0})),I\in\mathbb{S}\big\}$,
    it is sufficient to show that
    \[
        w(\phi_{I}(z))
        =\int_{\gamma_{[\phi_{I}(z_{0}),\phi_{I}(z),I ]}} \mathcal{U}\big(\phi_{I}(z),r\big)  f(r)~\mathrm{d}r \quad (I \in\mathbb{S})
    \]
    is a classical solution of \eqref{ia0}
    restricted to
    $\big\{p\in\Sigma_{\theta}(\phi_{0}):p\in\mathscr{X}_{\theta}(\phi_{I}(z_{0})),I\in\mathbb{S}\big\}$ with initial value
    $w(s_{0})=0$. Note that
    $Q(r)^{-1} f(r)$ admits right slice derivative
    in $\mathfrak{D}$ and that
    \begin{align}\label{ia6}\begin{split}
         \partial_{S}^{R}
          Q ^{-1}(r)f(r)
        = - \big(Q ^{-1}(r)\odot ^{R}\dot{Q}(r) \odot ^{R}Q ^{-1}(r)\big)f(r)
         + Q ^{-1}(r) \partial_{S}^{R}  f (r)
        =  Q ^{-1}(r) g(r),\end{split}
    \end{align}
    where $g(r) = \partial_{S}^{R}f(r) - \big(\dot{Q}(r)\odot^{R} Q ^{-1}(r)\big) f (r)$.
    Considering the right slice derivative of $\big(\mathcal{U}(q,r)\odot^{R}Q(r)^{-1}\big) f(r)$
    $\big(\phi_{I}(z)=q\in\mathscr{X}_{\theta}(r)\big)$
    with respect to $r$, by
    $(\mathfrak{E}_{3})$ and \eqref{ia6}, we obtain that
    \begin{equation}\label{kexi-jishu-05}
    \begin{aligned}
         &~\frac{\partial_{S}^{R}}{\partial r}
         \big(\mathcal{U}(q,r)\odot^{R}Q(r)^{-1}\big) f(r)\\
        =&~\big(-\mathcal{U}(q,r)\odot ^{R}A(r)\odot^{R}Q(r)^{-1}\big) f(r)
         +\big(\mathcal{U}(q,r)
        \odot^{R} Q(r)^{-1}\big) g(r)\\
        =&~\mathcal{U}(q,r) f(r) +
        \big(\mathcal{U}(q,r)\odot^{R} Q(r)^{-1}\big) \big[g(r)-\lambda f(r)\big].
    \end{aligned}
\end{equation}
    Integrating \eqref{kexi-jishu-05} along $\gamma_{[\phi_{I}(z_{0}),q,I ]}$ yields
    \begin{align}\label{ia7}
        \begin{split}
        w(q)=&~ Q ^{-1}(q) f (q)
         -\bigg[\big(\mathcal{U}(q,\phi_{I}(z_{0}))\odot^{R} Q(\phi_{I}(z_{0}))^{-1}\big) f(\phi_{I}(z_{0})) \\
         &~+\int_{\gamma_{[\phi_{I}(z_{0}),q,I ]}}
         \big(\mathcal{U}(q,r) \odot^{R} Q(r)^{-1}\big)
        \big[g(r)-\lambda f(r)\big]\mathrm{d}r\bigg]\\
        =&~ Q ^{-1}(q) f (q)-\nu(q),
        \end{split}
    \end{align}
    where
    \begin{align*}
        \nu(q)=&~\big(\mathcal{U}(q,\phi_{I}(z_{0}))\odot^{R} Q(\phi_{I}(z_{0}))^{-1}\big) f(\phi_{I}(z_{0}))\\
        &~+\int_{\gamma_{[\phi_{I}(z_{0}),q,I]}}
        \big(\mathcal{U}(q,r) \odot^{R} Q(r)^{-1}\big)
        \big[g(r)-\lambda f(r)\big]\mathrm{d}r.
    \end{align*}
    Since $ Q ^{-1}\big(\phi_{I}(z_{0})\big)f \big(\phi_{I}(z_{0})\big)\in\mathfrak{D}$
    and $r\mapsto  Q ^{-1}(r) \big(g (r)-\lambda_{0}f(r)\big) $ is right hyperholomorphic in
    $\mathfrak{D}$, according to Theorem \ref{ia8}, we have
    \begin{equation*}
        \partial_{S}^{R}\nu(q)
        =A(q)\nu (q)
        + Q ^{-1}(q) \big(g(q) -\lambda f (q)\big) ,
        \quad   q\in  \mathscr{X}_{\theta}\big(\phi_{I}(z_{0})\big) ,I\in\mathbb{S}.
    \end{equation*}
    Thus, it follows from \eqref{ia7} that
    \begin{align*}
         \partial^{R}_{S} w(q)
        =&~  \partial^{R}_{S}
         Q ^{-1}(q)f (q)  - \partial^{R}_{S} \nu(q)\\
        =&~ Q ^{-1}(q)g (q)-Q ^{-1}(q)\big(g (q)-\lambda f (q)\big) - A(q)\nu (q)\\
        =&~ A(q) w(q)+ \lambda Q ^{-1}(q) f(q)
        - \big(A (q)\odot^{R}Q ^{-1}(q)\big) f (q)\\
        =&~A(q)w(q)+f(q).
    \end{align*}
    Since  $ \partial_{S}^{R} \nu(q)$ and
        $ Q ^{-1}(q)g (q) $
    are right hyperholomorphic for every $I\in\mathbb{S}$ in $X$,
    then $\partial_{S}^{R}w(q) $ is right hyperholomorphic
    in $X$ for every $I\in\mathbb{S}$, which implies that
    $w(q)$ $\big(q=\phi_{I}(z),I\in\mathbb{S}\big)$
    is a classical solution of \eqref{ia0}
    restricted to
    $\big\{p\in\Sigma_{\theta}(\phi_{0}):p\in\mathscr{X}_{\theta}(\phi_{I}(z_{0})),I\in\mathbb{S}\big\}$ with $u_{0}=0$.

    Next, we will show the the uniqueness of $u$.
    Suppose that $u_{1}$ is another classical solution of
    \eqref{ia0}. According to the assumptions and Theorem \ref{ia4}, it follows that
    $r\mapsto\mathcal{U}(q,r) u_{1}(r)$ is right slice regular in $X$
    and that
    \begin{equation}\label{ia9}
         \frac{\partial^{R}_{S}}{\partial q}
         \mathcal{U}(q,r) u_{1}(r)
         =\mathcal{U}(q,r) f(r),
    \end{equation}
    integrating \eqref{ia9} along $\gamma_{[\phi_{I}(z_{0}),q,I]}$ yields
    $u_{1}(q)=u(q)$ for $q\in \mathscr{X}_{\theta}\big(\phi_{I}(z_{0})\big)$.
    This completes the proof.
\end{proof}

\subsection{Evolution System in Banach Bilateral
$\mathbb{H}$-Algebras under the Parabolic-Like Conditions}\label{subsec6}

In this subsection, we will study the quaternionic Cauchy problem
\begin{equation*}
    \left\{\begin{aligned}
     &\partial_{S}^{R} u(q)+A(q)u(q)=f(q),
    \\
    &u(s_{0})=u_{0},
    \end{aligned}\right. \quad s_{0},q\in\Sigma_{\theta}
\end{equation*}
under different conditions from those in Subsection \ref{sec5}.
For the sake of convenience, in the discussion of this subsection,
we always assume that $X$ is a Banach bilateral $\mathbb{H}$-module and
$\Omega_{\mathscr{D}}= \Sigma_{\theta} $ for some $\theta\in (0,\pi/2]$.

\subsubsection{The Homogeneous Evolution Equation under the Parabolic-Like Conditions}

Now we first consider the homogeneous case, i.e., the quaternionic Cauchy problem as follows:
\begin{equation}\label{ib0}
    \left\{\begin{aligned}\displaystyle
    &\partial_{S}^{R} u(q)+A(q)u(q)=0 \\
    &u(s_{0})=u_{0} ,
    \end{aligned}\right.\qquad s_{0},q\in\Sigma_{\theta}.
\end{equation}

We begin with a method for constructing quaternionic evolution systems,
which serves merely as a formal convention and assumption.
For each $q\in\Omega_{\mathscr{D}}$, suppose that
$-A(q)$ generates the quaternionic slice regular
semigroup $\{\mathcal{S}_{q}(w)\}_{w\in\Sigma_{\theta}'} $
on
$X$ and let $\mathcal{U}(q,w)$ be a quaternionic evolution system of \eqref{ib0}.
For $z\in\mathscr{D}$, let
\begin{align}\label{ib4}
    \begin{split}
     \mathcal{U}\big(\phi_{I}(z),\phi_{I}(z_{0})\big)
    =&~\mathcal{S}_{s}\big(\phi_{I}(z)-\phi_{I}(z_{0})\big)
    +\mathcal{W}\big(\phi_{I}(z),\phi_{I}(z_{0})\big)\\
    :=&~\mathcal{S}_{s}\big(\phi_{I}(z)-\phi_{I}(z_{0})\big)+
    \int_{\gamma_{[\phi_{I}(z_{0}),\phi_{I}(z),I]}} \mathcal{S}_{\xi}
    \big(\phi_{I}(z)-\xi\big)\odot^{R}\mathscr{R} \big(\xi,\phi_{I}(z_{0})\big)~\mathrm{d}\xi.
    \end{split}
\end{align}
Then
\begin{align*}
    \frac{\partial_{S}^{R}}{\partial q}
    \mathcal{U}\big(q,\phi_{I}(z_{0})\big)=&~-A\big(\phi_{I}(z_{0})\big) \odot^{R}
    \mathcal{S}_{\phi_{I}(z_{0})}\big(q-\phi_{I}(z_{0})\big)
    +\mathscr{R}\big(q,\phi_{I}(z_{0})\big)   \\
    &~-\int_{\gamma_{[\phi_{I}(z_{0}),q,I]}} A(\xi)\odot^{R}\mathcal{S}_{\xi}(q-\xi)
    \odot^{R}\mathscr{R}
    \big(\xi,\phi_{I}(z_{0})\big)~\mathrm{d}\xi,
\end{align*}
where $q=\phi_{I}(z)\in\mathscr{X}_{\theta}\big(\phi_{I}(z_{0})\big)$,  and hence
\begin{align}\label{ib2}
    \begin{split}
     \frac{\partial_{S}^{R}}{\partial q}\mathcal{U}\big(q,\phi_{I}(z_{0})\big)
    +A(q)\odot^{R}\mathcal{U}\big(q,\phi_{I}(z_{0})\big)
    =&~ \mathscr{R}\big(q,\phi_{I}(z_{0})\big) -\mathscr{R} _{1}\big(q,\phi_{I}(z_{0})\big)\\&~
    -\int_{\gamma_{[\phi_{I}(z_{0}),q,I]}}\mathscr{R}_{1}(q,\xi)\odot^{R}\mathscr{R}
    \big(\xi,\phi_{I}(z_{0})\big)~\mathrm{d}\xi,
    \end{split}
\end{align}
where
\begin{equation*}
    \mathscr{R} _{1}\big(q,\phi_{I}(z_{0})\big)  =\big(A(\phi_{I}(z_{0}))-A(q)\big)
    \odot^{R}\mathcal{S}_{\phi_{I}(z_{0})}\big(q-\phi_{I}(z_{0})\big).
\end{equation*}
According to \eqref{ib2}, $\mathscr{R}\big(q,\phi_{I}(z_{0})\big)$ satisfies
the integral equation
\begin{equation}\label{ib3}
    \mathscr{R}\big(\phi_{I}(z ),\phi_{I}(z_{0})\big)
    =\mathscr{R} _{1}\big(\phi_{I}(z ),\phi_{I}(z_{0})\big)
    +\int_{\gamma_{[\phi_{I}(z_{0}),\phi_{I}(z),I]}}\mathscr{R} _{1}\big(\phi_{I}(z),\xi\big)\odot^{R}
    \mathscr{R}\big(\xi,\phi_{I}(z_{0})\big)~\mathrm{d}\xi.
\end{equation}

Suppose that $\big\{A(q)\big\}_{ q\in\Sigma_{\theta}}$ is a family of closed
right linear quaternionic operators on $X$,
we introduce the following conditions $(\mathscr{P}_{1})$-$(\mathscr{P}_{3})$:
\begin{itemize}
    \item[$(\mathscr{P}_{1})$] The domain $\mathcal{D}\big(A(q)\big)=\mathfrak{D}$ of $A(q)$
    is dense in $X$ and independent of $q$, the function $q\mapsto A(q)$ is right
    slice regular and $A(q)$ is uniformly bounded.

    \item[$(\mathscr{P}_{2})$] $-A(q) $ is spherical $\theta$-sectorial
    operator with vertex $\omega=0$ for all $q\in\Omega_{\mathscr{D}}$ and
    $0\in\rho_{S}\big(A(q)\big)$. In addition, for any $\lambda \in -\Sigma_{\theta+\pi/2}$,
    There are constants $M$ and $C$ such that
 \begin{equation}\label{ao0}
        \big\Vert S_{L}^{-1}\big(\lambda,A(q)\big) \big\Vert_{\mathcal{B}^{R}(X)}
         \leq  \frac{M}{\vert\lambda\vert+1 } ,\quad
        \left\Vert \frac{\partial_{S}^{R}}{\partial q} S_{L}^{-1}\big(\lambda,A(q)\big)
        \right\Vert_{\mathcal{B}^{R}(X)} \leq
        \frac{C}{\vert\lambda\vert+1 }.
    \end{equation}
    The $\mathcal{B}^{R}(X)$-valued function
    $A(q) \odot^{R}\mathcal{S}_{\xi}(w)$
    is right slice regular for $w,q,\xi\in \Sigma_{\theta}$,
    where $\{\mathcal{S}_{\xi}(p)\}_{p\in\Sigma_{\theta}'}$ is quaternionic
    right slice regular semigroup generated by $A(\xi)$.

    \item[$(\mathscr{P}_{3})$] There exist constants $0<\alpha \leq  1$ and
    $L>0$ such that
    \begin{equation}\label{ta3}
        \big\Vert \big(A(p)-A(q)\big)\odot^{R}
         A(r)^{-1}\big\Vert_{\mathcal{B}^{R}(X)}
        \leq
         L\vert p-q\vert^{\alpha},
        \quad p,q,r\in \Sigma_{\theta} ,
    \end{equation}
    and
            \begin{equation}\label{ic9}
        \big\Vert A(q)\odot^{R}\big(\mathcal{S}_{\xi_{1}}(w)-\mathcal{S}_{\xi_{2}}(w)\big)
        \big\Vert_{\mathcal{B}^{R}(X)}
         \leq  C|w|^{-1}\vert \xi_{1}-\xi_{2}\vert^{\alpha},
    \end{equation}
    where $w,q,\xi_{1},\xi_{2}\in  \Sigma_{\theta} $.
    Moreover, there exists a constant $\widetilde{M}>0$ such that
    $A(\phi_{I}(z_{4}))\odot^{R} A(\phi_{I}(z_{5}))^{-1}<\widetilde{M}$ for all
    $z_{4},z_{5}\in\mathscr{D}$ and $I\in\mathbb{S}$.

\end{itemize}

\begin{remark}
    Note that the condition $(\mathscr{P}_{2})$ implies that
    $-A(q)$ generates a quaternionic right slice regular semigroup
    $\{\mathcal{S}_{q}(w)\}_{w\in\Sigma_{\theta}'}$ and by Theorem 7.1 in \cite{kexi-srs0}, there exist
    constants $M,M_{1}>0$ such that
    \begin{equation}\label{e6.13}
        \| \mathcal{S}_{q}(w)\|_{\mathcal{B}^{R}(X)}\leq M ,
    \end{equation}
    and
    \begin{equation}\label{kexi-fangsuo-1}
        \big\|S_{L}^{-1}(p,A(q))\big\|_{\mathcal{B}^{R}(X)}\leq \frac{M_{1}}{|p|},\quad \forall p\in
        \Sigma_{\pi/2+\theta}.
    \end{equation}
    Then by \eqref{kexi-fangsuo-0}, \eqref{kexi-tianjia-1} and \eqref{kexi-fangsuo-1}, there exists
    a constant $C>0$ such that
    \begin{equation}\label{e6.14}
        \big\|A(q)\odot^{R}\mathcal{S}_{q}(w)\big\|_{\mathcal{B}^{R}(X)}\leq \frac{C}{|w|}.
    \end{equation}
\end{remark}

From the assumptions $(\mathscr{P}_{1} )$-$(\mathscr{P}_{3} )$, some basic results can be proved which will applied later.
\begin{lemma}\label{ta0}
    Let $\{A(q)\}_{q\in\Sigma_{\theta}}$ be a family of closed right linear
    quaternionic operators which satisfy
    the conditions $(\mathscr{P}_{1})$-$(\mathscr{P}_{3})$. Then
    \begin{equation}\label{kexi-ic7}
        \big\Vert \big(A(q_{1})-A(q_{2})\big)\odot^{R}\mathcal{S}_{\xi}(w)
        \big\Vert_{\mathcal{B}^{R}(X)} \leq
        C\vert w\vert^{-1}\vert q_{1}-q_{2}\vert^{\alpha},
    \end{equation}
    where
        $q_{1},q_{2}\in \Sigma_{\theta}$.
    \begin{equation}\label{kexi-ic8}
        \big\Vert A(q)\odot^{R}\big(\mathcal{S}_{\xi}
        (s_{2})-\mathcal{S}_{\xi}(s_{1})\big)\big\Vert_{\mathcal{B}^{R}(X)}
         \leq  C\frac{\big| s_{1} - s_{2} \big|}{|s_{1}|\cdot|s_{2}|},
    \end{equation}
    where $q,\xi\in \Sigma_{\theta}$ and $s_{1},s_{2}\in\mathbb{C}_{I}$ $(I\in\mathbb{S})$.
    The constant $C$  in \eqref{kexi-ic7} and \eqref{kexi-ic8} is a generic constant.
\end{lemma}

\begin{proof}
    From Theorem 7.1 in \cite{kexi-srs0},
    it follows that $-A(q)$ generates a quaternionic slice regular semigroup
    $\big\{\mathcal{S}_{q}(w)\big\}_{w\in\Sigma_{\theta}'}$. Thus, according to \eqref{ta3} and \eqref{e6.14},
    there exists a constant $M$ such that
    \begin{align*}
        &~\big\Vert \big(A(q_{1})-A(q_{2})\big)\odot^{R}
        \mathcal{S}_{\xi}(w)\big\Vert_{\mathcal{B}^{R}(X)}\\
        \leq   &~\big\Vert \big(A(q_{1})-A(q_{2})\big)\odot^{R}A(\xi)^{-1}
        \big\Vert_{\mathcal{B}^{R}(X)}\cdot
        \Vert A(\xi)\odot^{R}\mathcal{S}_{\xi}(w)\Vert_{\mathcal{B}^{R}(X)}\\
        \leq  &~\frac{M}{|w|}\cdot L\vert q_{2}-q_{1}\vert^{\alpha}
        =:\frac{C_{1}}{|w|}\vert q_{2}-q_{1}\vert^{\alpha},
    \end{align*}
    and hence \eqref{kexi-ic7} holds.

    Next, we will show \eqref{kexi-ic8}.
    Let $ s_{1}, s_{2}\in\Sigma_{\theta}\cap\mathbb{C}_{I}$ and $u_{I,0}\in X$,
    since
    \begin{align*}
        &~\big(A(\xi)\odot^{R}\mathcal{S}_{\xi}(s_{2})\big)u_{I,0}-
        \big(A(\xi)\odot^{R}\mathcal{S}_{\xi}(s_{1})\big)u_{I,0}\\
        =&~ -\int_{\gamma_{[s_{1},s_{2},I]}}A(\xi)\odot^{R}\big(A(\xi)
        \odot^{R}\mathcal{S}_{\xi}(r)\big)u_{I,0}~\mathrm{d}r\\
        = &~-\int_{\gamma_{[s_{1},s_{2},I]}}\big(A(\xi)\odot^{R}\mathcal{S}_{\xi}(r/2)\big)^{2}x~\mathrm{d}r,
    \end{align*}
    from \eqref{e6.14}, it follows that
    there exists a constant $\widetilde{C}_{1}$ such that
    \begin{align*}
        &~\big\Vert \big(A(\xi)\odot^{R}\mathcal{S}_{\xi}(s_{2})\big) u_{I,0}-
        \big(A(\xi)\odot^{R}\mathcal{S}_{\xi}(s_{1}) \big)
        x\big\Vert_{X}\\
        \leq &~  4\widetilde{C}_{1}\Vert u_{I,0}\Vert_{X}\int_{\gamma_{[s_{1},s_{2},I]}}
        \frac{1}{|r|^{2}} ~|\mathrm{d}r|
        =  4\widetilde{C}^{2}_{1}\Vert u_{I,0}\Vert_{X}
        \frac{\big||s_{1}|-|s_{2}|\big|}{|s_{1}|\cdot|s_{2}|}\\
        \leq&~ 4\widetilde{C}^{2}_{1}\Vert u_{I,0}\Vert_{X}
        \frac{\big|s_{1}-s_{2}\big|}{|s_{1}|\cdot|s_{2}|}.
    \end{align*}
    Since $(\mathscr{P}_{3})$ is satisfied, then there exists a constant
    $\widetilde{C}_{2}>0$ such that
    \[
        \Vert A(q)\odot^{R}A(\xi)^{-1}\Vert_{\mathcal{B}^{R}(X)}\leq \widetilde{C}_{2}.
    \]
   Therefore,
    \begin{align*}
        &~\big\Vert A(q)\odot^{R}\big(\mathcal{S}_{\xi}(s_{1})-\mathcal{S}_{\xi}(s_{2})\big)
        \big\Vert_{\mathcal{B}^{R}(X)}\\
        \leq  &~\big\Vert A(q)\odot^{R}A(\xi)^{-1}\big\Vert_{\mathcal{B}^{R}(X)}
        \cdot\big\Vert A(\xi)\odot^{R}
        \big(\mathcal{S}_{\xi}(s_{1})-\mathcal{S}_{\xi}(s_{2})\big)\big\Vert_{\mathcal{B}^{R}(X)}\\
        \leq  &~ 4\widetilde{C}_{1}^{2}\widetilde{C}_{2}
        \frac{\big||s_{1}|-|s_{2}|\big|}{|s_{1}|\cdot|s_{2}|}
        =:C\frac{\big||s_{1}|-|s_{2}|\big|}{|s_{1}|\cdot|s_{2}|},
    \end{align*}
    i.e., \eqref{kexi-ic8} holds. This completes the proof.
\end{proof}

From Lemma \ref{ta0}, we can obtain an estimate inequality for $\mathscr{R}_{1}\big(q,\phi_{I}(z_{0})\big)$.
\begin{corollary}\label{la3}
    Let $\{A(q)\}_{q\in\Sigma_{\theta}} $
    be a family of closed linear quaternionic operators which satisfy
    the conditions $(\mathscr{P}_{1})$-$(\mathscr{P}_{3})$. Then the operator
    $\mathscr{R} _{1}\big(q,\phi_{I}(z_{0})\big)=\big(A(\phi_{I}(z_{0}))-A(q)\big)\odot^{R}\mathcal{S}_{\phi_{I}(z_{0})}\big(q-\phi_{I}(z_{0})\big)$
    $\big( q\in\mathscr{X}_{\theta}(\phi_{I}(z_{0}))\big)$ is
    right hyperholomorphic for every $I\in\mathbb{S}$ and
    there exists a constant $C>0$ such that
    \begin{equation}\label{am2}
        \big\Vert \mathscr{R} _{1}(q,\phi_{I}(z_{0}))\big\Vert_{\mathcal{B}^{R}(X)}  \leq   C
        \big|q-\phi_{I}(z_{0})\big|^{\alpha-1}\quad
        (0<\alpha \leq  1 ).
    \end{equation}
\end{corollary}

\begin{proof}
    The right hyperholomorphy of $\mathscr{R} _{1}\big(q,\phi_{I}(z_{0})\big)$ follows from the
    right slice regularity of
    $A(q)\odot^{R}\mathcal{S}_{\xi}\big(\phi_{I}(z_{0})\big)$
    in $\mathcal{B}^{R}(X)$. Furthermore, by $(\mathscr{P}_{3})$ and \eqref{e6.14},
    it follows that for $q\in\mathscr{X}_{\theta}\big(\phi_{I}(z_{0})\big)$, there exists a constant $C_{1}>0$ such that
    \begin{align*}
        &~\big\Vert \mathscr{R} _{1}\big(q,\phi_{I}(z_{0})\big)\big\Vert_{\mathcal{B}^{R}(X)}\\
         \leq &~\big\Vert\big(A(q)-A(\phi_{I}(z_{0}))\big)\odot^{R}
        A(\phi_{I}(z_{0}))^{-1}\big\Vert_{\mathcal{B}^{R}(X)}\cdot\big\Vert A(\phi_{I}(z_{0}))\odot^{R}
        \mathcal{S}_{\phi_{I}(z_{0})}(q-\phi_{I}(z_{0}))\big\Vert_{\mathcal{B}^{R}(X)}\\
         \leq &~LC_{1}\big|q-\phi_{I}(z_{0})\big|^{\alpha}|q-\phi_{I}(z_{0})|^{-1}=:C\big|q-\phi_{I}(z_{0})\big|^{\alpha-1}.
    \end{align*}
    This completes the proof.
\end{proof}

In what follows, we will establish our main result of this subsection.
\begin{theorem}\label{ht4}
    Under the assumptions $(\mathscr{P}_{1})$-$(\mathscr{P}_{3})$,
    there is a unique quaternionic evolution
    system $\mathcal{U}(q,w)$
    satisfying
    \begin{description}
        \item[$(\mathfrak{E}_{1})'$] $\Vert\mathcal{U}(q,w)\Vert_{\mathcal{B}^{R}(X)} \leq  C$
        for $q\in\mathscr{X}_{\theta}(w)$.

        \item[$(\mathfrak{E}_{2})^{+}$] For
        $q\in\mathscr{X}_{\theta}(\phi_{I}\big(z_{0})\big)$,
        $\mathcal{U}\big(q,\phi_{I}(z_{0})\big):X\rightarrow\mathfrak{D}$ and
        $q\mapsto\mathcal{U}\big(q,\phi_{I}(z_{0})\big)$ is right hyperholomorphic in $X$ for every $I\in\mathbb{S}$.
        Moreover,
        \begin{align}
            &\frac{\partial_{S}^{R}}{\partial q}\mathcal{U}\big(q,\phi_{I}(z_{0})\big)
            +A(q)\odot^{R}\mathcal{U}\big(q,\phi_{I}(z_{0})\big)=0 ,
            \label{fdq0}
            \\
            &\left\Vert\frac{\partial_{S}^{R}}{\partial q}\mathcal{U}\big(q,\phi_{I}(z_{0})\big)
            \right\Vert_{\mathcal{B}^{R}(X)}=\Vert A(q)\odot^{R}\mathcal{U}\big(q,\phi_{I}(z_{0})\big)\Vert_{\mathcal{B}^{R}(X)}
             \leq C\big|q-\phi_{I}(z_{0})\big|^{-1}\label{ka5}
        \end{align}
        and
        \begin{equation}\label{oa5}
            \big\Vert A(q)\odot^{R}\mathcal{U}\big(q,\phi_{I}(z_{0})\big)\odot^{R}
            A\big(\phi_{I}(z_{0})\big)^{-1}\big\Vert_{\mathcal{B}^{R}(X)} \leq  C  .
        \end{equation}

        \item[$(\mathfrak{E}_{3})^{+}$] For every $u_{I,0}\in\mathfrak{D}$,
        $q\in\mathscr{X}_{\theta}\big(\phi_{I}(z_{0})\big)$,
        $\mathcal{U}(q,w) u_{I,0}$ $\big(w=\phi_{I}(z_{0})\big)$
        admits right slice derivative with respect to
        $w$ and
        \begin{equation}\label{oa7}
            \frac{\partial_{S}^{R}}{\partial w}\mathcal{U}(q,w) u_{I,0}
            =\big(\mathcal{U}(q,w)\odot^{R}A(w)\big)
            u_{I,0}.
        \end{equation}
    \end{description}
\end{theorem}

In the sequel, we need to
divide the proof of Theorem \ref{ht4} into the following steps in which some
useful lemmas will be established.

\begin{lemma}\label{fd9}
    Let $\{A(q)\}_{q\in\Sigma_{\theta}}$ be a family of closed right linear
    quaternionic operators
    satisfying the assumptions $(\mathscr{P}_{1})$-$(\mathscr{P}_{3})$.
    Then there exists a quaternionic operator
    $\mathcal{U}(q,w)$ for
    $(q,w)\in \Sigma_{\theta} \times \Sigma_{\theta} $
    $(I\in\mathbb{S})$,
    satisfying  $(\mathfrak{E}_{1})'$.
\end{lemma}

\begin{proof}
    For $q\in\mathscr{X}_{\theta}\big(\phi_{I}(z_{0})\big)$,
    note that the integral equation \eqref{ib3} for $\mathscr{R}\big(q,\phi_{I}(z_{0})\big)$ can be solved
    by successive approximations. For $m\geqslant 1$, we define inductively
    \begin{equation}\label{ta7}
        \mathscr{R} _{m+1}\big(q,\phi_{I}(z_{0})\big)
        =\int_{\gamma_{[\phi_{I}(z_{0}),q, I ]}}\mathscr{R} _{1}(q,\xi)\odot^{R}
        \mathscr{R} _{m}\big(\xi,\phi_{I}(z_{0})\big)~\mathrm{d}\xi.
    \end{equation}
    By induction, the right hyperholomorphy of the quaternionic operator
    $\mathscr{R} _{m}\big(q,\phi_{I}(z_{0})\big)$ for every $I\in\mathbb{S}$
    is immediate.
    Since
    \begin{equation}\label{ya0}
        \int_{\gamma_{[\phi_{I}(z_{0}),q, I ]}}(q-\xi)^{\alpha-1}
        \big(\xi-\phi_{I}(z_{0})\big)^{\beta-1}~\mathrm{d}\xi
        =\big(q-\phi_{I}(z_{0})\big)^{\alpha+\beta-1} \Gamma(\alpha)\Gamma(\beta) \Gamma(\alpha+\beta)^{-1}
    \end{equation}
    holds for every $\alpha,\beta>0$, then it follows that
    \begin{equation}\label{ta6}
        \big\Vert \mathscr{R} _{m}(q,\phi_{I}(z_{0}))\big\Vert_{\mathcal{B}^{R}(X)}\leq
         \big(\widetilde{C}\Gamma(\alpha)\big)^{m} \Gamma(m\alpha)^{-1} \big|q-\phi_{I}(z_{0})\big|^{m\alpha-1}.
    \end{equation}
    Note that the integral given by \eqref{ta7} is an improper integral,
    and its existence can be directly inferred from \eqref{ta6}. Moreover,
    since $\mathscr{R} _{m}\big(q,\phi_{I}(z_{0})\big)$ and $\mathscr{R} _{1}\big(q,\phi_{I}(z_{0})\big)$ are right hyperholomorphic for every $I\in\mathbb{S}$, then
    according to \eqref{ta6},
    $\mathscr{R} _{m+1}\big(q,\phi_{I}(z_{0})\big)$ is also right hyperholomorphic
    for every $I\in\mathbb{S}$. On the other hand,
    by \eqref{ta6}, it follows that
    \[
        \mathscr{R}\big(q,\phi_{I}(z_{0})\big)
        =\sum_{m=1}^{\infty}\mathscr{R}_{m}\big(q,\phi_{I}(z_{0})\big)
    \]
    converges uniformly in the uniform operator topology for
    $q\in  \mathscr{X}_{\theta}\big(\phi_{I}(z_{0})\big)$.
    Therefore, $\mathscr{R}\big(q,\phi_{I}(z_{0})\big)$ is right hyperholomorphic
    in $\mathcal{B}^{R}(X)$, where $I\in\mathbb{S}$.
    Then we have
    \begin{equation}\label{ta8}
        \begin{aligned}
        \mathscr{R}\big(q,\phi_{I}(z_{0})\big)
        =&~\sum_{m=1}^{\infty}\mathscr{R} _{m}\big(q,\phi_{I}(z_{0})\big)\\
        =&~\mathscr{R}_{1}\big(q,\phi_{I}(z_{0})\big)
        +\sum_{m=1}^{\infty}
        \int_{\gamma_{[\phi_{I}(z_{0}),q,I]}}\mathscr{R}_{1}(q,\xi)\odot^{R}
        \mathscr{R}_{m}\big(\xi,\phi_{I}(z_{0})\big)~\mathrm{d}\xi.
        \end{aligned}
    \end{equation}
    Since $\mathscr{R}_{m}\big(q,\phi_{I}(z_{0})\big)$
    $(m\geqslant 1)$ is right hyperholomorphic for every $I\in\mathbb{S}$, then
    \eqref{am2} and \eqref{ta6} justify the interchange of summation and integration in
    \eqref{ta8}. This demonstrates that $\mathscr{R}\big(q,\phi_{I}(z_{0})\big)$ is a solution to the integral equation \eqref{ib3}. Furthermore, invoking Stirling's formula, we obtain that
    \begin{align}\label{ta9}
        \begin{split}
            \big\Vert \mathscr{R}\big(q,\phi_{I}(z_{0})\big)\big\Vert_{\mathcal{B}^{R}(X)}\leq   &~\sum_{m=1}^{\infty}
            \Gamma(m\alpha)^{-1}\big[\widetilde{C}\Gamma(\alpha)\big]^{m}|q-\phi_{I}(z_{0})|^{m\alpha-1}\\
            \leq   &~\left(\sum_{m=1}^{\infty}\Gamma(m\alpha)^{-1}\big[\widetilde{C}\Gamma(\alpha)\big]^{m}\tau^{\alpha(m-1)}\right)|q-\phi_{I}(z_{0})|^{\alpha-1}\\
            \leq   &~\widetilde{C}_{1}|q-\phi_{I}(z_{0})|^{\alpha-1}.
        \end{split}
    \end{align}
    By \eqref{e6.13}, there exists a constant $C_{1}$ such that
    $\big\Vert\mathcal{S}_{q}(\phi_{I}(z_{0}))\big\Vert_{\mathcal{B}^{R}(X)}\leq   C_{1}$.
    Constructing $\mathcal{U}\big(q,\phi_{I}(z_{0})\big)$ by \eqref{ib4},
    then from the right slice regularity of $\mathcal{S}_{q}\big(\phi_{I}(z_{0})\big)$ and
    \eqref{ta9}, it follows that $\mathcal{U}\big(q,\phi_{I}(z_{0})\big)$ is right hyperholomorphic for
    $q\in \mathscr{X}_{\theta}\big(\phi_{I}(z_{0})\big)$ $(I\in\mathbb{S})$
    and that
    \begin{align*}
        \big\Vert\mathcal{U}(q,\phi_{I}(z_{0}))\big\Vert_{\mathcal{B}^{R}(X)}
        \leq  &~\big\Vert\mathcal{S}_{\phi_{I}(z_{0})}(q-\phi_{I}(z_{0}))\big\Vert_{\mathcal{B}^{R}(X)}\\
        &~+
        \int_{\gamma_{[\phi_{I}(z_{0}),q,I]}}\big\Vert\mathcal{S}_{\xi}(q-\xi)\big\Vert_{\mathcal{B}^{R}(X)}\cdot
        \big\Vert \mathscr{R}(\xi,\phi_{I}(z_{0})) \big\Vert_{\mathcal{B}^{R}(X)}~|\mathrm{d}\xi|\\
        \leq  &~C_{1}+C_{1}\widetilde{C}_{1}\int_{\gamma_{[\phi_{I}(z_{0}),q,I]}}|\xi-\phi_{I}(z_{0})|^{\alpha-1}~|\mathrm{d}\xi|\leq   C,
    \end{align*}
    which implies that $(\mathfrak{E}_{1})'$ is satisfied.  This completes the proof.
\end{proof}

To demonstrate that the quaternionic operator
$\mathcal{U}(q,w)$
given by \eqref{ib4} is
a quaternionic evolution system, we still need to show that
$\mathcal{U}\big(q, \phi_{I}(z_{0})\big) = \mathcal{U}(q, r)\odot^{R}\mathcal{U}\big(r, \phi_{I}(z_{0})\big)$ for
$r\in\mathscr{X}_{\theta}\big(\phi_{I}(z_{0})\big)$,
$q\in\mathscr{X}_{\theta}(r)$.
The answer to this question will be derived from the uniqueness of the solution
of the quaternionic Cauchy problem \eqref{ib0}
which will be proved below (see Theorem \ref{ht3}), and the fact
that by $(\mathfrak{E}_{2})^{+}$ the solution of \eqref{ib0} restricted to
$\big\{p\in\Sigma_{\theta}:p\in\mathscr{X}_{\theta}(\phi_{I}(z_{0}))\big\}$ is
$\mathcal{U}\big(\phi_{I}(z),\phi_{I}(z_{0})\big) u_{I,0}$.
Firstly, let's show that $\mathcal{U}(q,w)$ satisfies $(\mathfrak{E}_{2})^{+}$,
for this, we need to establish some lemmas.
\begin{lemma}\label{la0}
    For each $\beta$ satisfies $0<\beta \leq \alpha \leq  1$, there is a constant $C_{\beta}$ such that
    \begin{equation}\label{ya1}
        \big\Vert\mathscr{R}_{1}(q,w)-\mathscr{R}_{1}(\xi,w)\big\Vert_{\mathcal{B}^{R}(X)}
         \leq  C_{\beta}|q-\xi|^{\beta} |\xi-w|^{\alpha-\beta-1} ,
    \end{equation}
    where $q\in\mathscr{X}_{\theta}(\xi)$, $\xi\in\mathscr{X}_{\theta}(w)$.
\end{lemma}

\begin{proof}
    For $q\in\mathscr{X}_{\theta}(\xi)$, $\xi\in\mathscr{X}_{\theta}(w)$,
    since
    \begin{align*}
        \mathscr{R}_{1}(q,w)-\mathscr{R}_{1}(\xi,w)
        =&~\big(A(\xi)-A(q)\big)\odot^{R}\mathcal{S}_{w}(q-w)\\ &~
        +
        \big(A(w)-A(\xi)\big)\odot^{R}\big(\mathcal{S}_{w}(q-w)-\mathcal{S}_{w}(\xi-w)\big),
    \end{align*}
    by \eqref{kexi-ic7}, we obtain that
    \begin{align*}
        &~\big\Vert \big(A(\xi)-A(q)\big)\odot^{R}
        \mathcal{S}_{w}(q-w)\big\Vert_{\mathcal{B}^{R}(X)}\\
        \leq&~
        \widetilde{C}
        |t-\xi|^{\alpha}|q-w|^{-1}\leq   \widetilde{C} |q-\xi|^{\alpha}|\xi-w|^{-1}.
    \end{align*}
    Moreover, by \eqref{ta3} and Lemma \ref{ta0}, we have
    \begin{align*}
        &~\big\Vert \big(A(w)-A(\xi)\big)\odot^{R}
        \big(\mathcal{S}_{w}(q-w)-\mathcal{S}_{w}(\xi-w)\big)
        \big\Vert_{\mathcal{B}^{R}(X)}\\
        \leq   &~\big\Vert \big(A(w)-A(\xi)\big) \odot^{R}A(w)^{-1}\big\Vert_{\mathcal{B}^{R}(X)}
        \cdot\big\Vert A(w) \odot^{R}\big(\mathcal{S}_{w}(q-w)-\mathcal{S}_{w}(\xi-w)\big)
        \big\Vert_{\mathcal{B}^{R}(X)}\\
        \leq  &~L|\xi-w|^{\alpha}\cdot \widetilde{C}_{1}|q-w|^{-1}|\xi-w|^{-1}|q-\xi|\\
        \leq   &~\widetilde{C}_{2}|\xi-w|^{\alpha-2}|q-\xi|.
    \end{align*}
    On the other hand, from \eqref{ta3}, it follows that
    \[
        \big\Vert\big(A(w)-A(\xi)\big)\odot^{R}
        \big(\mathcal{S}_{w}(q-w)
        -\mathcal{S}_{w}(\xi-w)\big)\big\Vert_{\mathcal{B}^{R}(X)}
        \leq   \widetilde{C}_{3}|\xi-w|^{\alpha-1}.
    \]
    Thus,
    \begin{align*}
        &~\big\Vert\big(A(w)-A(\xi)\big)\odot^{R}\big(\mathcal{S}_{w}(q-w)-\mathcal{S}_{w}(\xi-w)\big)\big\Vert_{\mathcal{B}^{R}(X)}\\
        \leq  &~ C_{1}\big(|\xi-w|^{\alpha-2}|q-\xi|\big)^{\alpha}\big(|\xi-w|^{\alpha-1}
        \big)^{1-\alpha}
        =\widetilde{C}_{4}|q-\xi|^{\alpha}|\xi-w|^{-1},
    \end{align*}
    which implies that
    \begin{equation}\label{ta4}
       \big\Vert \mathscr{R}_{1}(q,w)-\mathscr{R}_{1}(\xi,w)\big\Vert_{\mathcal{B}^{R}(X)}\leq
       C_{2}|q-\xi|^{\alpha}
       |\xi-w|^{-1}.
    \end{equation}
    Moreover, by \eqref{am2}, we obtain that
    \begin{align}\label{ta5}
        \begin{split}
        \big\Vert \mathscr{R}_{1}(q,w)-\mathscr{R}_{1}(\xi,w)\Vert_{\mathcal{B}^{R}(X)}\leq   &~\big\Vert\mathscr{R}_{1}(q,w)
        \big\Vert_{\mathcal{B}^{R}(X)}+\big\Vert\mathscr{R}_{1}(\xi,w)\big\Vert_{\mathcal{B}^{R}(X)}\\
        \leq &~ C_{3}\big(|q-w|^{\alpha-1}+|\xi-w|^{\alpha-1}\big)
        \leq C_{4}|\xi-w|^{\alpha-1}.
        \end{split}
    \end{align}
    Interpolating \eqref{ta4} and \eqref{ta5} for
    $\big\Vert\mathscr{R}_{1}(q,w)-\mathscr{R}_{1}(\xi,w)\big\Vert_{\mathcal{B}^{R}(X)}$ we have
    \begin{align*}
        \big\Vert \mathscr{R}_{1}(q,w)-\mathscr{R}_{1}(\xi,w)\big\Vert_{\mathcal{B}^{R}(X)}\leq
        &~C_{\beta,1}\big(|q-\xi|^{\alpha}|\xi-w|^{-1}\big)^{\beta/\alpha}\big(|\xi-w|^{\alpha-1}\big)^{1-\beta/\alpha} \\
        \leq  &~ C_{\beta}|q-\xi|^{\beta}|\xi-s|^{\alpha-\beta-1}.
    \end{align*}
    This completes the proof.
\end{proof}

Naturally, the following estimate inequality of $\mathscr{R}(q,w)-\mathscr{R}(\xi,w) $
can be obtained based on Lemma \ref{la0}.
\begin{corollary}\label{ka4}
    For every $\beta$ with $0<\beta<\alpha \leq  1$, there is a constant $C_{\beta}>0$ such that
    \begin{equation}\label{kexi-ya2}
        \big\Vert \mathscr{R}(q,w)-\mathscr{R}(\xi,w) \big\Vert_{\mathcal{B}^{R}(X)} \leq
        C_{\beta}|q-\xi|^{\beta}|\xi-w|^{\alpha-\beta-1},
    \end{equation}
    where
    $\xi\in\mathscr{X}_{\theta}(w)$, $q\in\mathscr{X}_{\theta}(\xi)$.
\end{corollary}

\begin{proof}
    For $\xi\in\mathscr{X}_{\theta}(w)$, $q\in\mathscr{X}_{\theta}(\xi)$,
    it follows from \eqref{ib3} that
    \begin{align*}
        \mathscr{R}(q,w) -\mathscr{R}(\xi,w)
        =&~\mathscr{R}_{1}(q,w)-\mathscr{R}_{1}(\xi,w)
        +\int_{\gamma_{[\xi,q, I ]}}\mathscr{R}_{1}(q,r)\odot^{R}
        \mathscr{R}(r,w)~\mathrm{d}r\\
        &~+\int_{\gamma_{[w,\xi, I ]}}\big(\mathscr{R}_{1}(q,r)-\mathscr{R}_{1}(\xi,r)\big)
        \odot^{R}
        \mathscr{R}(r,w)~\mathrm{d}r.
    \end{align*}
    According to \eqref{am2} and \eqref{ta9}, we have
    \begin{align}\label{ya3}
        \begin{split}
        &~\bigg\Vert\int_{\gamma_{[\xi,q, I ]}}\mathscr{R}_{1}(q,r)\odot^{R}
         \mathscr{R}(r,w)~\mathrm{d}r
        \bigg\Vert_{\mathcal{B}^{R}(X)}\\
        \leq   &~ C_{1}\int_{\gamma_{[\xi,q, I ]}}|q-r|^{\alpha-1}|r-w|^{\alpha-1}~|\mathrm{d}r |\\
        \leq  &~
        C_{1}|\xi-w|^{\alpha-1} \int_{\gamma_{[\xi,q, I ]}}|q-r|^{\alpha-1}~|\mathrm{d}r|\\
        = &~ \frac{C_{1}}{\alpha}|\xi-w|^{\alpha-1}|q-\xi|^{\alpha}
        = \frac{C_{1}}{\alpha}|\xi-w|^{\alpha-\beta-1}|\xi-w|^{\beta}|q-\xi|^{\alpha}\\
        =&~  \frac{C_{1}}{\alpha}|\xi-w|^{\alpha-\beta-1}|q-\xi|^{\beta}
        |q-\xi|^{\alpha-\beta}|\xi-w|^{\beta}\\
        \leq  &~\frac{\zeta ^{\alpha}C_{1}}{\alpha}|\xi-w|^{\alpha-\beta-1}|q-\xi|^{\beta},
        \end{split}
    \end{align}
    where $\zeta=|q-w|$.
    Therefore, it follows from \eqref{ya0}, \eqref{ta9} and \eqref{ya1} that
    \begin{align}\label{kexi-ya4}
        \begin{split}
        &~\bigg\Vert \int_{\gamma_{[w,\xi, I ]}}
        \big(\mathscr{R}_{1}(q,r)-\mathscr{R}_{1}(\xi,r)\big)\odot^{R}
        \mathscr{R}(r,w)~|\mathrm{d}r|\bigg\Vert_{\mathcal{B}^{R}(X)}\\
        \leq  &~\widetilde{C}_{\beta}\widetilde{C}_{1}(q-\xi)^{\beta}\int_{\gamma_{[w,\xi, I ]}}
        |\xi-r|^{\alpha-\beta-1}|r-w|^{\alpha-1}~|\mathrm{d}r|\\
        \leq   & ~\frac{\widetilde{C}_{\beta}\widetilde{C}_{1}
        \Gamma(\alpha-\beta)\Gamma(\alpha)}{\Gamma(2\alpha-\beta)}
        |q-\xi|^{\beta}|\xi-w|^{2\alpha-\beta-1}\\
        \leq   &~\frac{\tau^{\alpha}\widetilde{C}_{\beta}\widetilde{C}_{1}\Gamma(\alpha-\beta)\Gamma(\alpha)}{\Gamma(2\alpha-\beta)}|q-\xi|^{\beta}|\xi-w|^{\alpha-\beta-1}.
        \end{split}
    \end{align}
    By virtue of \eqref{ya1}, \eqref{ya3} and \eqref{kexi-ya4},  then \eqref{kexi-ya2} follows. This completes the proof.
\end{proof}

Now we will show that  $(\mathfrak{E}_{2})^{+}$ is satisfied for $\mathcal{U}(q,w)$.
\begin{lemma}\label{ao4}
    Let $\mathcal{U}(q,w)$ be the quaternionic operator
    given by \eqref{ib4}.
    Then $\mathcal{U}(q,w)$ satisfies $(\mathfrak{E}_{2})^{+}$.
\end{lemma}

\begin{proof}
    First, for $q\in\mathscr{X}_{\theta}\big(\phi_{I}(z_{0})\big)$, we will show the right hyperholomorphy of $\mathcal{U}\big(q,\phi_{I}(z_{0})\big)$ for every $I\in\mathbb{S}$.
    Since $\mathcal{S}_{\phi_{I}(z_{0})}\big(q-\phi_{I}(z_{0})\big)$
    is right hyperholomorphic for $q\in\mathscr{X}_{\theta}\big(\phi_{I}(z_{0})\big)$, then
    \[
        \frac{\partial_{S}^{R}}{\partial q}\mathcal{S}_{\phi_{I}(z_{0})}\big(q-\phi_{I}(z_{0})\big)
        =-A\big(\phi_{I}(z_{0})\big) \odot^{R}
        \mathcal{S}_{\phi_{I}(z_{0})}\big(q-\phi_{I}\big(z_{0})\big)
    \]
    is a bounded right linear closed quaternionic operator which is right hyperholomorphic in
    $\mathcal{B}^{R}(X)$, it is sufficient to show that
    $\mathcal{W}\big(q,\phi_{I}(z_{0})\big)$ $\big(q\in\mathscr{X}_{\theta}(\phi_{I}(z_{0}))\big)$ is right hyperholomorphic
    for every $I\in\mathbb{S}$, where
    $\mathcal{W}\big(q,\phi_{I}(z_{0})\big)$ is given by \eqref{ib4}. Now let
    $\gamma :[a,b]\rightarrow \mathbb{C}_{ I }$ be a $C^{1}$-path
    from $\phi_{I}(z_{0})$ to $q$ and let
    \begin{equation}\label{oa2}
        \mathcal{W}_{\varepsilon}\big(q,\phi_{I}(z_{0})\big)
        =\int_{\gamma [\phi_{I}(z_{0}),\gamma(b-\varepsilon),I]}\mathcal{S}_{\xi}(q-\xi)\odot^{R}
        \mathscr{R}\big(\xi,\phi_{I}(z_{0})\big) ~\mathrm{d}\xi \quad (0<\varepsilon <b-a).
    \end{equation}
    Then it follows that
    $\mathcal{W}_{\varepsilon}\big(q,\phi_{I}(z_{0})\big)\rightarrow\mathcal{W}\big(q,\phi_{I}(z_{0})\big)$ as $\varepsilon\rightarrow 0$ and
    $\mathcal{W}_{\varepsilon}\big(q,\phi_{I}(z_{0})\big)$ is
    right hyperholomorphic in $q$ and
    \begin{equation}\label{oa3}\begin{split}
        \frac{\partial_{S}^{R}}{\partial q}\mathcal{W}_{\varepsilon}\big(q,\phi_{I}(z_{0})\big)
        = &~\mathcal{S}_{\gamma(b-\varepsilon)}\big(q-\gamma(b-\varepsilon)\big)\odot^{R}
       \mathscr{R}\big(\gamma(b-\varepsilon),\phi_{I}(z_{0})\big) \\
       &~
        -\int_{\gamma_{[\phi_{I}(z_{0}),\gamma(b-\varepsilon), I ]}}A(\xi) \odot^{R}
        \mathcal{S}_{\xi}(q-\xi) \odot^{R}
        \mathscr{R}\big(\xi,\phi_{I}(z_{0})\big)~\mathrm{d}\xi.
    \end{split}
    \end{equation}
    Since
    $A(q)\odot^{R}\mathcal{S}_{q}(q-\xi)=(\partial_{S}^{R}/\partial \xi)\mathcal{S}_{q}(q-\xi)$,
    we obtain that
    \begin{align}\label{fd8}
        \begin{split}
         &~\frac{\partial_{S}^{R}}{\partial q}\mathcal{W}_{\varepsilon}(q,\phi_{I}\big(z_{0})\big) \\
        =&~\mathcal{S}_{\gamma(b-\varepsilon)}\big(q-\gamma(b-\varepsilon)\big) \odot^{R}
        \mathscr{R}
        \big(\gamma(b-\varepsilon),\phi_{I}(z_{0})\big)\\
        &~+
        \int_{\gamma[\phi_{I}(z_{0}),\gamma(b-\varepsilon), I ]}\big(A(q)\odot^{R}
        \mathcal{S}_{q}(q-\xi)-A(\xi)\odot^{R}
        \mathcal{S}_{\xi}(q-\xi)\big)\odot^{R}\mathscr{R}\big(\xi,\phi_{I}(z_{0})\big)~\mathrm{d}\xi\\
        &~+\int_{\gamma[\phi_{I}(z_{0}),\gamma(b-\varepsilon),I]}A(q)\odot^{R}
        \mathcal{S}_{q}(q-\xi)\odot^{R}
        \big(\mathscr{R}\big(q,\phi_{I}(z_{0})\big)-\mathscr{R}\big(\xi,\phi_{I}(z_{0})\big)\big)~\mathrm{d}\xi\\
        &~+
        \big(\mathcal{S}_{q}\big(q-\phi_{I}(z_{0})\big)-\mathcal{S}_{q}(q-\gamma(b-\varepsilon))\big)\odot^{R}
        \mathscr{R}\big(q,\phi_{I}(z_{0})\big).
        \end{split}
    \end{align}
    From \eqref{e6.13},
    it follows that there exists a constant $\widetilde{C}>0$ such that
    $\Vert \mathcal{S}_{q}(w)\Vert_{\mathcal{B}^{R}(X)}\leq   \widetilde{C}$,
    by \eqref{ta9}, we have the first and the last terms on the
    right-hand side of \eqref{fd8} are bounded in norm by
    $2\widetilde{C}\widetilde{C}_{1}|\phi_{ I }(b-\varepsilon)-\phi_{I}(z_{0})|^{\alpha-1}$. Moreover, by
    \eqref{kexi-ic7} and \eqref{ic9} we obtain
    that there exist constants $C_{1}, C_{2}>0$ such that
    \[
        \Vert A(q)\odot^{R}
        \mathcal{S}_{q}(q-\xi)-A(\xi)\odot^{R}
        \mathcal{S}_{\xi}(q-\xi)\Vert_{\mathcal{B}^{R}(X)}\leq
        (C_{1}+C_{2})|q-\xi|^{\alpha-1},
    \]
    and by \eqref{ya0} we have
    \begin{align*}
        &~\bigg\Vert\int_{\gamma_{[\phi_{I}(z_{0}),\gamma(b-\varepsilon), I ]}}\big(A(q)\odot^{R}
        \mathcal{S}_{q}(q-\xi)-A(\xi)\odot^{R}
        \mathcal{S}_{\xi}
        (-\xi)\big) \odot^{R}
       \mathscr{R}\big(\xi,\phi_{I}(z_{0})\big)~\mathrm{d}\xi\bigg\Vert_{\mathcal{B}^{R}(X)}\\
        \leq   &~(C_{1}+C_{2}) \widetilde{C}_{1}\int_{\gamma_{[\phi_{I}(z_{0}),\gamma(b-\varepsilon), I ]}}
        |q-\xi|^{\alpha-1}\big|\xi-\phi_{I}(z_{0})\big|^{\alpha-1}~|\mathrm{d}\xi|\\
        \leq   &~ \frac{(C_{1}+C_{2})\widetilde{C}_{1}\Gamma(\alpha)^{2}}{\Gamma(2\alpha)}\big|q-\phi_{I}(z_{0})\big|^{2\alpha-1}\\
        \leq   &~\frac{(C_{1}+C_{2})\widetilde{C}_{1}\Gamma(\alpha)^{2}}{\Gamma(2\alpha)}\big|q-\phi_{I}(z_{0})\big|^{\alpha}\big|\gamma(b-\varepsilon)-\phi_{I}(z_{0})\big|^{\alpha-1}\\
        \leq  &~\frac{(C_{1}+C_{2})\widetilde{C}_{1}\tau^{\alpha}\Gamma(\alpha)^{2}}{\Gamma(2\alpha)}\big|\gamma(b-\varepsilon)-\phi_{I}(z_{0})\big|^{\alpha-1}.
    \end{align*}
    It follows from \eqref{e6.14} and \eqref{kexi-ya2} that there exists constant $C_{3}>0$ such that
    \begin{align*}
        &~\bigg\Vert\int_{\gamma_{[\phi_{I}(z_{0}),\gamma(b-\varepsilon),I]}}A(q)\odot^{R}
        \mathcal{S}_{q}(q-\xi)\odot^{R}
         \big(\mathscr{R}\big(q,\phi_{I}(z_{0})\big)-
        \mathscr{R}\big(\xi,\phi_{I}(z_{0})\big)\big)
        ~ \mathrm{d}\xi \bigg\Vert_{\mathcal{B}^{R}(X)}\\
        \leq  &~C_{3}C_{\beta}\int_{\gamma_{[\phi_{I}(z_{0}),\gamma(b-\varepsilon), I ]}} |q-\xi|^{\beta-1}\big|\xi-\phi_{I}(z_{0})\big|^{\alpha-\beta-1}~|\mathrm{d}\xi|\\
        \leq   &~ \frac{C_{3}C_{\beta}\Gamma(\beta)\Gamma(\alpha-\beta)}{\Gamma(\alpha)}\big|q-\phi_{I}(z_{0})\big|^{\alpha-1}
        \leq   \frac{C_{3}C_{\beta}\Gamma(\beta)\Gamma(\alpha-\beta)}{\Gamma(\alpha)}
        \big|\gamma(b-\varepsilon)-\phi_{I}(z_{0})\big|^{\alpha-1}.
    \end{align*}
    Finally, we obtain that there exists a constant $C_{4}>0$ such that
    \begin{equation}\label{oa0}
        \bigg\Vert\frac{\partial_{S}^{R}}{\partial q}\mathcal{W}_{\varepsilon}\big(q,\phi_{I}(z_{0})\big)
        \bigg\Vert_{\mathcal{B}^{R}(X)}
        \leq   \frac{C_{4}}{\big|\gamma(b-\varepsilon)-\phi_{I}(z_{0})\big|^{1-\alpha}},
    \end{equation}
    where $C_{4}$ is independent of $\varepsilon>0$.
    Letting $\varepsilon\rightarrow 0$ on
    the right-hand side of \eqref{fd8}, from Definition \ref{kexi-lian},
    it follows that
    $(\partial_{S}^{R}/\partial q)\mathcal{W}_{\varepsilon}\big(q,\phi_{I}(z_{0})\big)$
    converges strongly as $\varepsilon\rightarrow 0$. Denoting its limit by
    $\widetilde{\mathcal{W}}\big(q,\phi_{I}(z_{0})\big)$, then
    \begin{align}\label{oa8}
        \begin{split}
        \widetilde{\mathcal{W}}\big(q,\phi_{I}(z_{0})\big)=&~\mathcal{S}_{q}\big(q-\phi_{I}(z_{0})\big)\odot^{R}
        \mathscr{R}\big(q,\phi_{I}(z_{0})\big)\\
        &~+\int_{\gamma_{[\phi_{I}(z_{0}),q,I]}}
        \big(A(q)\odot^{R}
        \mathcal{S}_{q}(q-\xi)
        -A(\xi)\odot^{R}
        \mathcal{S}_{\xi}(q-\xi)\big)\odot^{R}
        \mathscr{R}\big(\xi,\phi_{I}(z_{0})\big)~\mathrm{d}\xi\\
        &~+\int_{\gamma_{[\phi_{I}(z_{0}),q,I]}}A(q)\odot^{R}
        \mathcal{S}_{q}(q-\xi)\odot^{R}
        \big(\mathscr{R}\big(q,\phi_{I}(z_{0})\big)-\mathscr{R}\big(\xi,\phi_{I}(z_{0})\big)\big)~\mathrm{d}\xi,
        \end{split}
    \end{align}
    it follows that $\widetilde{\mathcal{W}}\big(q,\phi_{I}(z_{0})\big)$ is right hyperholomorphic for
    every
    $I\in\mathbb{S}$, where $q\in\mathscr{X}_{\theta}\big(\phi_{I}(z_{0})\big)$.
    Letting $\varepsilon\rightarrow 0$ in \eqref{oa0} yields
    \begin{equation}\label{oa1}
        \big\Vert\widetilde{\mathcal{W}} \big(q,\phi_{I}(z_{0})\big)\Vert_{\mathcal{B}^{R}(X)}\leq
        C_{4}\big|q-\phi_{I}(z_{0})\big|^{\alpha-1}.
    \end{equation}
    Now, letting $\varepsilon\rightarrow 0$ in
    \[
        \mathcal{W}_{\varepsilon}\big(q_{2},\phi_{I}(z_{0})\big)-\mathcal{W}_{\varepsilon}
        \big(q_{1},\phi_{I}(z_{0})\big)
        =\int_{\gamma_{[q_{1},q_{2}, I ]}}
        \frac{\partial_{S}^{R}}{\partial\xi}\mathcal{W}_{\varepsilon}\big(\xi,\phi_{I}(z_{0})\big)~\mathrm{d}\xi,
    \]
    where $q_{1},q_{2}\in \gamma_{[\phi_{I}(z_{0}),q, I ]}([a,b])$,
    we obtain that
    \[
    \mathcal{W}\big(q_{2},\phi_{I}(z_{0})\big)-\mathcal{W}\big(q_{1},\phi_{I}(z_{0})\big)=\int_{\gamma_{[q_{1},q_{2}, I ]}}
    \widetilde{\mathcal{W}}\big(\xi,\phi_{I}(z_{0})\big)~\mathrm{d}\xi.
    \]
    Since $\widetilde{\mathcal{W}} \big(q,\phi_{I}(z_{0})\big)$ is
    right hyperholomorphic for every $I\in\mathbb{S}$,
    where $q\in\mathscr{X}_{\theta}\big(\phi_{I}(z_{0})\big)$,
    then $\mathcal{W}\big(q,\phi_{I}(z_{0})\big)$ is right hyperholomorphic with respect
    to $q$ for every $I\in\mathbb{S}$ and
    \[
        \frac{\partial_{S}^{R}}{\partial q}\mathcal{W}
        \big(q,\phi_{I}(z_{0})\big)=\widetilde{\mathcal{W} }
        \big(q,\phi_{I}(z_{0})\big).
    \]
    Thus, $\mathcal{U}\big(q,\phi_{I}(z_{0})\big)$ is right hyperholomorphic and satisfies
    \[
        \frac{\partial_{S}^{R}}{\partial q}\mathcal{U}\big(q,\phi_{I}(z_{0})\big)
        =-A\big(\phi_{I}(z_{0})\big)\odot^{R}
        \mathcal{S}_{\phi_{I}(z_{0})}\big(q-\phi_{I}(z_{0})\big)
        +\frac{\partial_{S}^{R}}{\partial q}
        \mathcal{W}\big(q,\phi_{I}(z_{0})\big),
    \]
    then it follows from \eqref{e6.14} and \eqref{oa1} that there exists a constant
    $C>0$ such that
    \[
        \bigg\Vert\frac{\partial_{S}^{R}}{\partial q}\mathcal{U}\big(q,\phi_{I}(z_{0})\big)
        \bigg\Vert_{\mathcal{B}^{R}(X)}\leq \frac{C}{\big|q-\phi_{I}(z_{0})\big|},
    \]
    which implies that \eqref{ka5} holds.

    Now we will show that $\mathcal{U}\big(q,\phi_{I}(z_{0})\big)$ satisfies \eqref{fdq0}. Setting
    \[
        \mathcal{U}_{\varepsilon}\big(q,\phi_{I}(z_{0})\big)=\mathcal{S}_{s}\big(q-\phi_{I}(z_{0})\big)+\mathcal{W}_{\varepsilon}\big(q,\phi_{I}(z_{0})\big) \quad \big(q\in \mathscr{X}_{\theta}(\phi_{I}(z_{0}))\big),
    \]
    where $\varepsilon>0$.
    Since $\mathcal{U}_{\varepsilon}\big(q,\phi_{I}(z_{0})\big):X\rightarrow\mathfrak{D}$,
    from \eqref{oa2} and \eqref{oa3}, we have
    \begin{align}\label{oa4}
        \begin{split}
        &~\frac{\partial_{S}^{R}}{\partial q}\mathcal{U}_{\varepsilon}\big(q,\phi_{I}(z_{0})\big)+A(q)\odot^{R}
        \mathcal{U}_{\varepsilon}(q,\phi_{I}\big(z_{0})\big)\\
        =&~\mathcal{S}_{\gamma(b-\varepsilon)}\big(q-\gamma(b-\varepsilon)\big)\odot^{R}
        \mathscr{R}\big(\gamma(b-\varepsilon),\phi_{I}(z_{0})\big)-\mathscr{R}_{1}\big(q,\phi_{I}(z_{0})\big) \\
        &~
         - \int_{\gamma_{[\phi_{I}(z_{0}),\gamma(b-\varepsilon),I]}}\mathscr{R}_{1}(q,\xi)\odot^{R}
        \mathscr{R}\big(\xi,\phi_{I}(z_{0})\big)~\mathrm{d}\xi.
        \end{split}
    \end{align}
    Letting $\varepsilon\rightarrow 0$, the right hand side of \eqref{oa4}
    tends strongly to zero. Since
    $(\partial_{S}^{R}/\partial q)\mathcal{U}_{\varepsilon}\big(q,\phi_{I}(z_{0})\big)$ strongly converges to
    $(\partial_{S}^{R}/\partial q)\mathcal{U}\big(q,\phi_{I}(z_{0})\big)$,
    it follows from \eqref{oa4}
    that $A(q)\odot^{R}\mathcal{U}_{\varepsilon}\big(q,\phi_{I}(z_{0})\big)$ converges strongly as $\varepsilon\rightarrow 0$.
    Let $u_{0}\in X$, since
    $\mathcal{U}_{\varepsilon}\big(q,\phi_{I}(z_{0})\big) u_{I,0}$ converges to
    $\mathcal{U}\big(q,\phi_{I}(z_{0})\big)u_{I,0}$, by the closedness of $A(q)$,
    it follows that
    $\mathcal{U}\big(q,\phi_{I}(z_{0})\big) u_{I,0}\in\mathfrak{D}$ and that
    \[
        \big(A(q)\odot^{R}\mathcal{U}_{\varepsilon}(q,\phi_{I}(z_{0}))\big) u_{I,0}\rightarrow
       \big( A(q)\odot^{R}\mathcal{U}(q,\phi_{I}(z_{0}))\big) u_{I,0}.
    \]
    Therefore, taking strong limit as $\varepsilon\rightarrow 0$ in
    \eqref{oa4} yields
    \[
        \frac{\partial_{S}^{R}}{\partial q}\mathcal{U}\big(q,\phi_{I}(z_{0})\big)+A(q)\odot^{R}\mathcal{U}\big(q,\phi_{I}(z_{0})\big)=0,
        \quad \big(q\in\mathscr{X}_{\theta}(\phi_{I}(z_{0}))\big)   ,
    \]
    which implies that $\mathcal{U}\big(q,\phi_{I}(z_{0})\big)$ satisfies \eqref{fdq0}.

    Next, we will show that $\mathcal{U}\big(q,\phi_{I}(z_{0})\big)$ satisfies
    \eqref{oa5}. Let $u_{0}\in X$ and consider the function
    \[
        \psi\big(\phi_{I}(z_{0})\big)
        =\big(\mathcal{S}_{q}(q-\phi_{I}(z_{0}))\odot^{R}
        \mathcal{U}\big(\phi_{I}(z_{0}),\xi)\odot^{R}A(\xi)^{-1}\big)u_{I,0} \quad
        \big(\phi_{I}(z_{0})\in\mathscr{X}_{\theta}(\xi)\big).
    \]
    It follows that $\psi$ admits right slice derivative with respect to
    $w$ $\big(w=\phi_{I}(z_{0})\big)$ and that
    \[
         \partial_{S}^{R} \psi(w)
        =\big(\mathcal{S}_{q}(q-w)\odot^{R} (A(q)-A(w))\odot^{R}\mathcal{U}(w,\xi)\odot^{R}
        A(\xi)^{-1}\big)u_{I,0}.
    \]
    Integrating $(\partial_{S}^{R}/\partial_{S}^{R} w)\psi(w)$ along
    $\gamma_{[\xi,q, I ]}$ and applying $A(q)$ to the result we obtain that
    \begin{equation}\label{e6.39}
        \begin{aligned}
        \mathscr{Z}(q,\xi) u_{I,0}
        =&~\big(A(q)\odot^{R}\mathcal{S}_{q}(q-\xi)\odot^{R}A(\xi)^{-1}\big) u_{I,0}\\
        &~
        +\int_{\gamma_{[\xi,q, I ]}}
        \big(\mathscr{Y}(q,w)\odot^{R}\mathscr{Z}(w,\xi)\big) u_{I,0}~\mathrm{d}w,
        \end{aligned}
    \end{equation}
    where
    \[
        \mathscr{Y}(q,w)=A(q)\odot^{R}\mathcal{S}_{q}(q-w)\odot^{R}
        \big(A(q)-A(w)\big)\odot^{R}A(w)^{-1},
        \]
    and
    \[
        \mathscr{Z}(q,\xi)=A(q)\odot^{R}\mathcal{U}(q,\xi)\odot^{R}A(\xi)^{-1}.
    \]
    By \eqref{ta3} and \eqref{e6.13},
    it follows that there exists a constant $\widetilde{M}$ such that
    \begin{align*}
        \Vert A(q)\odot^{R}\mathcal{S}_{q}(q-\xi)\odot^{R}A(\xi)^{-1}\Vert_{\mathcal{B}^{R}(X)}
        \leq   &~\widetilde{M},
    \end{align*}
    and from \eqref{ta3} and \eqref{e6.14}, it follows
    that there exists a constant $\widetilde{M}_{1}$ such that
    \begin{align*}
        \Vert \mathscr{Y}(q,w)\Vert_{\mathcal{B}^{R}(X)}
        \leq &~  \Vert A(q)\odot^{R}\mathcal{S}_{q}(q-w)\Vert_{\mathcal{B}^{R}(X)}\cdot
        \big\Vert(A(q)-A(w))\odot^{R}A(w)^{-1}\big\Vert_{\mathcal{B}^{R}(X)}\\
        \leq &~  L\widetilde{M}_{1}|q-w|^{\alpha-1}.
    \end{align*}
    Therefore, by \eqref{e6.39}, we obtain that
    \[
        \Vert \mathscr{Z}(q,\xi) u_{I,0}\Vert_{X}\leq
        \widetilde{M}M\tau^{-1}\Vert u_{I,0}\Vert_{X}+L\widetilde{M}_{1}\int_{\gamma_{[\xi,q, I ]}}|q-w|^{\alpha-1}
        \Vert \mathscr{Z}(w,\xi)u_{I,0}\Vert_{X}~|\mathrm{d}w|.
    \]
    It follows from \eqref{ya0} that there exists a constant $C>0$ such that
    $\Vert \mathscr{Z}(q,\xi) u_{I,0}\Vert_{X}\leq   C\Vert u_{I,0}\Vert_{X}$, and hence
    \[
        \Vert \mathscr{Z}(q,\xi)\Vert_{\mathcal{B}^{R}(X)}=\Vert A(q) \odot^{R}
        \mathcal{U}(q,\xi) \odot^{R}
        A(\xi)^{-1}\Vert_{\mathcal{B}^{R}(X)}\leq    C,
    \]
    which implies \eqref{oa5} holds. This completes the proof.
\end{proof}

Next, we will prove the uniqueness of $\mathcal{U}(q,w)$.

\begin{lemma}\label{ao5}
    Let $\mathcal{U}(q,w)$ be the quaternionic operator given by \eqref{ib4}.
    Then $\mathcal{U}(q,w)$ is unique and
    satisfies the condition $(\mathfrak{E}_{3})^{+}$ for
    $q\in\mathscr{X}_{\theta}\big(\phi_{I}(z_{0})\big)$.
\end{lemma}

\begin{proof}
    By $(\mathscr{P}_{2})$, there exists a constant $C_{1}$ such that
    \begin{equation}\label{na1}
        \left\Vert \frac{\partial_{S}^{R}}{\partial q} S_{L}^{-1}\big(\lambda,A(\phi_{I}(z_{0}))\big)
        \right\Vert_{\mathcal{B}^{R}(X)}
        \leq   \frac{C_{1}}{\vert\lambda\vert+1}\quad
        \text{for } \lambda\in\rho_{S}(A(\phi_{I}(z_{0}))).
    \end{equation}
    Moreover, since $-A(q)$ is spherical $\theta$-sectorial
    operator with vertex $\omega=0$, it follows from Theorem 6.9 in \cite{kexi-srs0} that
    \begin{equation*}
        \begin{aligned}
        \mathcal{S}_{\phi_{I}(z_{0})}\big(q-\phi_{I}(z_{0})\big)
        =\frac{1}{2\pi}\int_{ \Gamma(I;r;\theta)}S_{L}^{-1}\big(\lambda,A\big(\phi_{I}(z_{0})\big)\big)  I ^{-1}
        ~\mathrm{d}\lambda ~ \exp^{-\lambda}\big(q-\phi_{I}(z_{0})\big)  ,
        \end{aligned}
    \end{equation*}
    where $r\in (0,\infty)$ and $\Gamma(I;r;\theta)=\gamma_{1}\cup\gamma_{2}\cup\gamma_{3}$
    and
        \begin{align*}
         \gamma_{1} =&~\{-t e^{- I (\theta+\pi/2-\varepsilon)}:t\in(-\infty,-r]\},\\
        \gamma_{2} =&~
        \{r e^{I\alpha }: -(\pi/2+\theta-\varepsilon)\leq \alpha \leq (\pi/2+\theta-\varepsilon)\},\\
         \gamma_{3} =&~\{t e^{ I (\pi/2+\theta-\varepsilon)}:r\leq t <\infty\} .
    \end{align*}
    On the other hand, since $q\in\mathscr{X}_{\theta}\big(\phi_{I}(z_{0})\big)$,
    and $\mathcal{S}_{w}(q-w)$ is right hyperholomorphic in $w$ $\big(w=\phi_{I}(z_{0})\big)$
    for every $I\in\mathbb{S}$, we have
    \begin{align*}
        \frac{\partial_{S}^{R}}{\partial w}\mathcal{S}_{w}(q-w)
        =&~\frac{1}{2\pi}
        \int_{ \Gamma(I;r;\theta)}
         S_{L}^{-1}\big(\lambda,A(w)\big)\lambda  I ^{-1}~\mathrm{d}\lambda~
         \exp^{-\lambda}(q-w) \\
        &~+\frac{1}{2\pi} \int_{\Gamma(I;r;\theta)}
         \frac{\partial_{S}^{R}}{\partial q} S_{L}^{-1}\big(\lambda,A(w)\big)I^{-1}_{w}~
        \mathrm{d}\lambda~\exp^{-\lambda}(q-w) \\
        =&~-\frac{\partial_{S}^{R}}{\partial q}\mathcal{S}_{w}(q-w)
        +\frac{1}{2\pi}
        \int_{ \Gamma(I;r;\theta)}
         \frac{\partial_{S}^{R}}{\partial q} S_{L}^{-1}\big(\lambda,A(w)\big)I^{-1}_{w}
        ~\mathrm{d}\lambda~\exp^{-\lambda}(q-w).
    \end{align*}
    To show $\mathsf{(\mathfrak{E}_{3})^{+}}$ is satisfied, let $\mathfrak{V}(q,w)$
    $\big(w=\phi_{I}(z_{0})\big)$
    be a quaternionic operator-valued function satisfying
    \begin{equation}\label{na3}
        \left\{\begin{aligned}
             &\frac{\partial_{S}^{R}}{\partial w}\mathfrak{V}(q,w) u_{I,0}
             =\big(\mathfrak{V}(q,w)\odot^{R}A(w)\big) u_{I,0},\quad
               q\in  \mathscr{X}_{\theta}(w),\\
             &\mathfrak{V}(q,q)=\mathcal{I},
        \end{aligned}\right.
    \end{equation}
    and construct it as the same as the construction of $\mathcal{U}(q,w)$ given
    by \eqref{ib4}.
    Let
    \begin{align}\label{xc1}\begin{split}
        \mathscr{Q}_{1}(q,w)=&~\left(\frac{\partial_{S}^{R}}{\partial q}
        +\frac{\partial_{S}^{R}}{\partial w}\right)\mathcal{S}_{w}
        (q-w)\\
        =&~\frac{1}{2\pi}\int_{\Gamma(I;r;\theta)}
        \frac{\partial_{S}^{R}}{\partial w}S_{L}^{-1}
        \big(\lambda,A(w)\big) I ^{-1}~\mathrm{d}\lambda~
        \exp^{-\alpha}(q-w).
    \end{split}
    \end{align}
    It follows from \eqref{na1} that there exists a constant $\widetilde{C}_{1}>0$ such that
\begin{equation*}
\begin{aligned}
    &~\bigg\Vert \frac{1}{2\pi}\int_{\gamma_{3}}
        \frac{\partial_{S}^{R}}{\partial q}
        S_{L}^{-1}\big(\lambda,A(w)\big)
        I^{-1} \mathrm{d}\lambda ~\exp^{-\lambda}(q-w) \bigg\Vert_{\mathcal{B}^{R}(X)} \\
    \leq &~\frac{C_{1}}{2\pi}
    \Bigg\Vert \int_{\gamma_3} e^{-\lambda(q-w)}
    \frac{1}{|\lambda|+1}\, \mathrm{d}\lambda \Bigg\Vert
    \leq  \frac{C_{1}}{2\pi}\int_r^\infty
    \frac{\bigl| e^{-t e^{i(\pi/2+\theta-\varepsilon)}(q-w)} \bigr|}{t+1}\, \mathrm{d}t \\
    \leq&~ \frac{C_{1}}{2\pi}
    \int_r^\infty \frac{e^{-t\, |q-w|\, \cos(\pi/2+\theta-\varepsilon+\arg(q-w))}}{t}\, \mathrm{d}t \\
    \leq&~\frac{C_{1}}{2\pi}
    \int_r^\infty \frac{e^{-t\, c_\varepsilon}}{t}\, dt \quad \text{for some } c_\varepsilon>0 \\
    =&~ \frac{C_{1}}{2\pi}\int_{r\, c_\varepsilon}^{\infty} \frac{e^{-\xi}}{\xi}\,
    \mathrm{d}\xi \le \widetilde{C}_{1}.
    \end{aligned}
\end{equation*}
    The integral on $\gamma_{1} $ is estimated similarly as
    $\gamma_{3} $. For $\gamma_{2} $, we have there exists a constant $C_{2}>0$ such that
    \[
\begin{aligned}
\Bigg\Vert \int_{\gamma_2} e^{-\lambda(q-w)} \frac{1}{|\lambda|+1}\, \mathrm{d}\lambda \Bigg\Vert
&\leq \frac{C_{1}}{2\pi}\int_{-(\pi/2+\theta-\varepsilon)}^{\pi/2+\theta-\varepsilon}
   \frac{\bigl| e^{-r e^{i\alpha} (q-w)} \bigr|}{r+1} \, r\, \mathrm{d}\alpha \\
&= \frac{C_{1}}{2\pi}\int_{-(\pi/2+\theta-\varepsilon)}^{\pi/2+\theta-\varepsilon}
   \frac{e^{-r |q-w| \cos(\alpha+\arg(q-w))}}{r+1} \, r\, \mathrm{d}\alpha \\
&\leq \frac{C_{1}}{2\pi}\int_{-(\pi/2+\theta-\varepsilon)}^{\pi/2+\theta-\varepsilon}
   e^{-r c_\varepsilon} \, \mathrm{d}\alpha \quad \text{for some } c_\varepsilon>0 \\
&= \frac{C_{1}}{2\pi} 2(\pi/2+\theta-\varepsilon) \, e^{-r c_\varepsilon} \le C_{2}.
\end{aligned}
\]
Thus, there is a constant $C>0$ such that
    $\Vert \mathfrak{Q}_{1}(q,w)\Vert_{\mathcal{B}^{R} (X)}\leq   C$.

    Now, set
    \begin{equation}\label{na2}
        \mathfrak{Q}\big(q,\phi_{I}(z_{0})\big)=\mathfrak{Q}_{1}\big(q,\phi_{I}(z_{0})\big)+\int_{\gamma_{[\phi_{I}(z_{0}),q, I ]}}\mathfrak{Q}(q,\xi)
        \odot^{R}\mathfrak{Q}_{1}\big(\xi,\phi_{I}(z_{0})\big)~\mathrm{d}\xi\quad
        \big(q\in\mathscr{X}_{\theta}(\phi_{I}(z_{0}))\big).
    \end{equation}
    Note that the solution $\mathfrak{Q}\big(q,\phi_{I}(z_{0})\big)$ of \eqref{na2} can be obtained by the same way as the solution of the integral equation \eqref{ib3}
    in the proof of Lemma \ref{fd9}. On the other hand,
    since in this case $\mathfrak{Q}_{1}\big(q,\phi_{I}(z_{0})\big)$ is uniformly bounded,
    the solution $\mathfrak{Q}\big(q,\phi_{I}(z_{0})\big)$ of  \eqref{na2}
    will satisfy
    \[
        \Vert \mathfrak{Q}\big(q,\phi_{I}(z_{0})\big)\Vert_{\mathcal{B}^{R}(X)}\leq   C.
    \]
    Let
    \[
        \mathfrak{V}(q,w)=\mathcal{S}_{w}(q-w)+\int_{\gamma_{[w,q,I]}}\mathfrak{Q}(q,\xi)\odot^{R}
        \mathcal{S}_{w}(\xi-w)~\mathrm{d}\xi \quad \big(w=\phi_{I}(z_{0})\big).
    \]
    Then $\Vert \mathfrak{V}(q,w)\Vert_{\mathcal{B}^{R}(X)}\leq   C$
    and $\mathfrak{V}(q,w) u_{I,0}$ admits right slice derivative
    at $w$ with $u_{I,0}\in\mathfrak{D}$.
    Calculating the right slice derivative of $\mathfrak{V}(q,w) u_{I,0}$
    with respect to $w$ yields
    \begin{align*}
        &~\frac{\partial_{S}^{R}}{\partial w}\mathfrak{V}(q,w) u_{I,0}-\big(
            \mathfrak{V}(q,w)  \odot^{R}A(w)\big) u_{I,0}\\
        =&~\mathfrak{Q}_{1}(q,w) u_{I,0}+\int_{\gamma_{[w,q, I ]}}
        \big( \mathfrak{Q}(q,\xi) \odot^{R}\mathfrak{Q}_{1}(\xi,w)\big)
         u_{I,0}~\mathrm{d}\xi-\mathfrak{Q}(q,w)  u_{I,0}=0.
    \end{align*}
    From the definition of $\mathfrak{V}(q, w)$, it follows  that $\mathfrak{V}(q, q) = \mathcal{I}$ and
    hence $\mathfrak{V}(q, w)$ is a
    solution of \eqref{na3}.

    For given $u_{0}\in X$, $q\in \mathscr{X}_{\theta}(r)$, $r\in\mathscr{X}_{\theta}\big(\phi_{I}(z_{0})\big)$,
     the function
    $r\mapsto \mathfrak{V}(q,r)\odot^{R}\mathcal{U}\big(r,\phi_{I}(z_{0})\big) u_{I,0}$ admits right slice
    derivative in $r$
    and
    \begin{align*}
        &~\Big(\frac{\partial_{S}^{R}}{\partial r}\mathfrak{V}(q,r)
        \odot^{R}\mathcal{U}(r,\phi_{I}(z_{0}))\Big) u_{I,0}\\
        =&~\big(\mathfrak{V}(q,r)\odot^{R}A(r)\odot^{R}\mathcal{U}(r,\phi_{I}(z_{0}))\big) u_{I,0}
        -\big(\mathfrak{V}(q,r)\odot^{R}A(r)\odot^{R}\mathcal{U}
        \big(r,\phi_{I}(z_{0})\big)\big) u_{I,0}=0,
    \end{align*}
    which implies that $\big(\mathfrak{V}(q,r)\odot^{R}\mathcal{U}(r,\phi_{I}(z_{0}))\big) u_{I,0}$ is independent of $r$ for $q\in \mathscr{X}_{\theta}(r)$, $r\in\mathscr{X}_{\theta}\big(\phi_{I}(z_{0})\big)$. Letting
    $r\rightarrow \phi_{I}(z_{0})$ and $r\rightarrow q$ respectively,
    we obtain that $\mathfrak{V}\big(q,\phi_{I}(z_{0})\big) u_{I,0}=\mathcal{U}\big(q,\phi_{I}(z_{0})\big) u_{I,0}$ for every
    $u_{I,0}\in X$. Thus,
    $\mathcal{U}\big(q,\phi_{I}(z_{0})\big)=\mathfrak{V}\big(q,\phi_{I}(z_{0})\big)$
    and $\mathcal{U}\big(q,\phi_{I}(z_{0})\big)$ satisfies
    \begin{equation*}
        \frac{\partial_{S}^{R}}{\partial w}\mathcal{U}(q,w) u_{I,0}
        =\big(\mathcal{U}(q,w) \odot^{R}A(w)\big) u_{I,0},\quad w=\phi_{I}(z_{0}),
    \end{equation*}
    i.e., \eqref{oa7} holds. This completes the proof.
\end{proof}

In what follows, we are ready to show that the quaternionic operator constructed by
\eqref{ib4} is the unique quaternionic evolution system for the quaternionic Cauchy problem.
\begin{theorem}\label{ht3}
    Let $\{A(q)\}_{q\in\Sigma_{\theta}}$ satisfy the conditions
    $(\mathscr{P}_{1})$-$(\mathscr{P}_{3})$.
    For each $u_{0}\in X$,
    the quaternionic Cauchy problem
    \begin{equation}\label{ya7}
        \left\{\begin{aligned}
            & \partial_{S}^{R}  u(q)+A(q)u(q)=0, \\
            &u(s_{0})=u_{0},
        \end{aligned}\right.\qquad s_{0},q\in\Sigma_{\theta}
    \end{equation}
    has a unique solution $u$ and its restriction to
    $\big\{p\in\Sigma_{\theta}:p\in\mathscr{X}_{\theta}(\phi_{I}(z_{0}))\big\}$ is given by
    $u(q)=\mathcal{U}\big(q,\phi_{I}(z_{0})\big) u_{I,0}$ for
    $q\in\mathscr{X}_{\theta}\big(\phi_{I}(z_{0})\big)$, where $\mathcal{U}(q,w)$ is
    given by \eqref{ib4}.
\end{theorem}

\begin{proof}
    It follows from \eqref{fdq0} that
    $u(q)=\mathcal{U}\big(q,\phi_{I}(z_{0})\big) u_{I,0}$ is a solution of
    \eqref{ya7} restricted to
    $\big\{p\in\Sigma_{\theta}:p\in\mathscr{X}_{\theta}(\phi_{I}(z_{0}))\big\}$, where
    $q\in\mathscr{X}_{\theta}(s_{0})$. Suppose that $v(q)$ is another solution of \eqref{ya7}.
    Since for
    every $r\in \mathscr{X}_{\theta}\big(\phi_{I}(z_{0})\big)$
    such that $q\in\mathscr{X}_{\theta}(r)$, we have
    $v(q)\in\mathfrak{D}$, then by \eqref{oa7} and \eqref{ya7}, it follows that the function
    $r\mapsto\mathcal{U}(q,r) v(r)$ admits right slice
    derivative and
    \begin{align*}
         \frac{\partial_{S}^{R}}{\partial r}\mathcal{U}(q,r) v(r)
        = \big(\mathcal{U}(q,r)\odot^{R}A(r)\big) v(r)-\big(\mathcal{U}(q,r)
        \odot^{R}A(r)\big) v(r)=0,
    \end{align*}
    and hence $\mathcal{U}(q,r) v(r)$ is constant
    for $r\in\mathscr{X}_{\theta}\big(\phi_{I}(z_{0})\big)$
    $\big(q\in\mathscr{X}_{\theta}(r)\big)$. Since it is right hyperholomorphic for every
    $I\in\mathbb{S}$,
    we can let $r\rightarrow q$ and $r\rightarrow \phi_{I}(z_{0})$
    to obtain $\mathcal{U}\big(q,\phi_{I}(z_{0})\big) u_{I,0}=v(q)$ and the
    uniqueness of the solution of \eqref{ya7} follows. This completes the proof.
\end{proof}

\noindent
\textbf{The Proof of Theorem \ref{ht4}.}
\begin{proof}
    By virtue of Lemmas \ref{fd9}, \ref{ao4} and \ref{ao5}, there is a
     quaternionic operator $\mathcal{U}(q,w)$
    satisfying $(\mathfrak{E}_{1})'$, $(\mathfrak{E}_{2})^{+}$ and $(\mathfrak{E}_{3})^{+}$. Moreover,
    by Theorem \ref{ht3}, we obtain that for $u_{0}\in X$,
    \[
        \mathcal{U}\big(q,\phi_{I}(z_{0})\big) u_{I,0}
        =\big(\mathcal{U}(q,r)\odot^{R}\mathcal{U}\big(r,\phi_{I}(z_{0})\big)\big) u_{I,0}
        \quad \big(r\in\mathscr{X}_{\theta}(\phi_{I}(z_{0})),q\in\mathscr{X}_{\theta}(r) \big),
    \]
    and hence $\mathcal{U}(q,w)$ is a quaternionic evolution system satisfying $(\mathfrak{E}_{1})'$
    $(\mathfrak{E}_{2})^{+}$ and
    $(\mathfrak{E}_{3})^{+}$. If $V(q,w)$ is
    another quaternionic evolution system satisfying $(\mathfrak{E}_{1})'$
    and $(\mathfrak{E}_{2})^{+}$, then
    $V\big(q,\phi_{I}(z_{0})\big) \nu$ is a solution of \eqref{ya7} and it follows from Theorem \ref{ht3}  that
    $V\big(q,\phi_{I}(z_{0})\big) \nu=\mathcal{U}\big(q,\phi_{I}(z_{0})\big) \nu$, thus $V(q,w)=\mathcal{U}(q,w)$,
    i.e., $\mathcal{U}(q,w)$
    is the unique quaternionic evolution
    system satisfying
    $(\mathfrak{E}_{1})'$, $(\mathfrak{E}_{2})^{+}$ and
    $(\mathfrak{E}_{3})^{+}$. This completes the proof.
\end{proof}

\subsubsection{The Inhomogeneous Evolution Equation under the Parabolic-Like Conditions}

Let $\{A(q)\}_{q\in\Sigma_{\theta}}$
satisfy the conditions $(\mathscr{P}_{1})$-$(\mathscr{P}_{3})$
and consider the quaternionic inhomogeneous Cauchy problem as follows
\begin{equation}\label{ua0}
    \left\{\begin{aligned}
        & \partial_{S}^{R}u (q)+Au(q)=f(q), \\
        &u(s_{0})=u_{0},
    \end{aligned}\right.\qquad s_{0},q\in\Sigma_{\theta}
\end{equation}
in the Banach bilateral $\mathbb{H}$-module $X$.

In the following, a theorem for the existence and uniqueness of solutions to \eqref{ua0} will be
established.

\begin{theorem}\label{th0}
    Let $\{A(q)\}_{q\in\Omega_{\mathscr{D}}}$ satisfy the conditions $(\mathscr{P}_{1})$-$(\mathscr{P}_{3})$
    and let $\mathcal{U}(q,w)$ be the quaternionic
    evolution system provided by Theorem \ref{ht4}. If $f$ is right slice regular
     on $\Omega_{\mathscr{D}}$, then the quaternionic Cauchy problem \eqref{ua0} has
    a unique solution $u$ and its restriction to
    $\big\{p\in\Omega_{\mathscr{D}}:p\in\mathscr{X}_{\theta}(\phi_{I}(z_{0}))\big\}$
    is given by
    \begin{equation*}
        u(q)=
            \mathcal{U}\big(q,\phi_{I}(z_{0})\big) u_{I,0}
            +\int_{\gamma_{[\phi_{I}(z_{0}),q, I ]}}
            \mathcal{U}(q,r) f(r)~\mathrm{d}r,
    \end{equation*}
    where $q\in\mathscr{X}_{\theta}\big(\phi_{I}(z_{0})\big)$.
\end{theorem}

\begin{proof}
    It follows from Theorem \ref{ht3} that
    $\mathcal{U}(\phi_{I}(z),\phi_{I}(z_{0})) u_{I,0}$  $(z\in\mathscr{D},\phi_{I}(z)=q)$
    is the unique solution of
    the homogeneous quaternionic Cauchy problem
    \begin{equation}\label{ua7}
         \partial_{S}^{R} u(q)+A(q)u(q)=0,\quad u(s_{0})=u_{0}.
    \end{equation}
    To show the existence of a solution of the Cauchy problem \eqref{ua0},
    it is sufficient to show that
    \[
        v\big(\phi_{I}(z)\big)=\int_{\gamma_{[\phi_{I}(z_{0}),\phi_{I}(z),I]}}
        \mathcal{U}(q,r) f(r)~\mathrm{d}r\quad (z\in\mathscr{D},I\in\mathbb{S})
    \]
    is a solution of the Cauchy problem
    \begin{equation}\label{ua8}
        \left\{
            \begin{aligned}
                & \partial_{S}^{R} v(q)+A(q)v(q)=f(q)\\
                & v(s_{0})=0.
            \end{aligned}\right.
    \end{equation}
    Thus, set
    \begin{equation}\label{ua5}
        v_{\varepsilon}(q)
        =\int_{\gamma_{[\phi_{I}(z_{0}),\gamma(b-\varepsilon),I]}}
        \mathcal{U}(q,r) f(r)~\mathrm{d}r.
    \end{equation}
    According to Theorem \ref{ht4}, we obtain that $v_{\varepsilon}(q)$
    admits right slice derivative at $q=\phi_{I}(z)$ and
    \begin{align*}
        \frac{\partial}{\partial q}v_{\varepsilon}(q)
        =&~\frac{\partial_{S}^{R}}{\partial q}\int_{\gamma_{[\phi_{I}(z_{0}),\gamma(b-\varepsilon), I ]}}
        \mathcal{U}(q,r) f(r)~\mathrm{d}r\\
        =&~\mathcal{U}\big(q,\gamma(b-\varepsilon)\big) f(\gamma(b -\varepsilon))
         -\int_{\gamma_{[\phi_{I}(z_{0}),\gamma(b-\varepsilon), I ]}}\big(A(q)\odot^{R}
         \mathcal{U}(q,r)\big) f(r)~\mathrm{d}r,
    \end{align*}
    which implies that
    \begin{equation}\label{ua2}
         \partial_{S}^{R} v_{\varepsilon}(q)+ A(q)
         v_{\varepsilon}(q)=\mathcal{U}\big(q,\gamma(b-\varepsilon)\big)
         f\big(\gamma(b-\varepsilon)\big).
    \end{equation}
    We conclude the proof by letting $\varepsilon\rightarrow0$ in \eqref{ua2}. To substantiate this passage to the limit, we will demonstrate that
    $ \partial_{S}^{R} v_{\varepsilon}(q)$ has a
    limit as $\varepsilon \to 0$.

   By \eqref{ib4} and Theorem \ref{ht4}, it follows that
    \begin{equation*}
        \mathcal{U}\big(q,\phi_{I}(z_{0})\big)
        =\mathcal{S}_{\phi_{I}(z_{0})}\big(q-\phi_{I}(z_{0})\big)
        +\mathcal{W}\big(q,\phi_{I}(z_{0})\big)\quad
        \big(q\in \mathscr{X}_{\theta}(\phi_{I}(z_{0}))\big),
    \end{equation*}
    where $\{\mathcal{S}_{q}(w)\}_{w\in\Sigma_{\theta}'}$ is
    the quaternionic right slice regular semigroup generated by $-A(q)$ and
    $\mathcal{W}\big(q,\phi_{I}(z_{0})\big)$ admits right slice derivative
    at $q$ and
    \begin{align*}
        &~ \frac{\partial_{S}^{R}}{\partial q}\int_{\gamma_{[\phi_{I}(z_{0}),\gamma(b-\varepsilon),I]}}\mathcal{W}(q,r) f(r)~\mathrm{d}r\\
        = &~\mathcal{W}\big(q,\gamma(b-\varepsilon)\big) f\big(\gamma(b-\varepsilon)\big)
        +\int_{\gamma_{[\phi_{I}(z_{0}),\gamma(b-\varepsilon),I]}}
        \frac{\partial_{S}^{R}}{\partial q} \mathcal{W}(q,r) f(r)~\mathrm{d}r.
    \end{align*}
    On the other hand,
    by \eqref{oa1}, we obtain that
    $\big\Vert(\partial_{S}^{R}/\partial q)\mathcal{W}\big(q,\phi_{I}(z_{0})\big)\big\Vert_{\mathcal{B}^{R}(X)}\leq   C\big|q-\phi_{I}(z_{0})\big|^{\alpha-1}$, and
    since
    $\mathcal{W}(q,q)=0$, then
    \[
        \lim\limits_{\varepsilon\rightarrow 0}\frac{\partial_{S}^{R}}{\partial q}
        \int_{\gamma_{[\phi_{I}(z_{0}),\gamma(b-\varepsilon),I]}}
        \mathcal{W}(q,r) f(r)~\mathrm{d}r
        =\int_{\gamma_{[\phi_{I}(z_{0}), q,I]}}
        \frac{\partial_{S}^{R}}{\partial q}\mathcal{W} (q,r) f(r)~\mathrm{d}r.
    \]
    Moreover,
    \begin{align*}
        \begin{split}
            &~\frac{\partial_{S}^{R}}{\partial q}\int_{\gamma_{[\phi_{I}(z_{0}),\gamma(b-\varepsilon), I ]}}\mathcal{S}_{r}(q-r) f(r)~\mathrm{d}r\\
            =&~\mathcal{S}_{\gamma(b-\varepsilon)}\big(q-\gamma(b-\varepsilon)\big)
            f\big(\gamma(b-\varepsilon)\big)
            +\int_{\gamma_{[\phi_{I}(z_{0}),\gamma(b-\varepsilon), I ]}}\frac{\partial_{S}^{R}}{\partial q}
            \mathcal{S}_{r}(q-r) f(r)~\mathrm{d}r\\
            =&~\mathcal{S}_{\gamma(b-\varepsilon)} \big(q-\gamma(b-\varepsilon)\big)
            f(\gamma(b-\varepsilon))\\
            &~+\int_{\gamma_{[\phi_{I}(z_{0}),\gamma(b-\varepsilon), I ]}}
            \big[A(q)\odot^{R}\mathcal{S}_{q}(q-r)-A(r)\odot^{R}
            \mathcal{S}_{r}(q-r)\big] f(r)~\mathrm{d}r
            \\ &~-\int_{\gamma_{[\phi_{I}(z_{0}),\gamma(b-\varepsilon), I ]}}
            \big(A(q)\odot^{R} \mathcal{S}_{q}(q-r)\big) f(r)~\mathrm{d}r\\
            =:&~\mathscr{I}_{1}+\mathscr{I}_{2}-\mathscr{I}_{3}.
        \end{split}
    \end{align*}
    It follows from Definition \ref{kexi-lian} that
    \[
        \lim\limits_{\varepsilon\rightarrow 0 }\mathscr{I}_{1}=
        \lim\limits_{\varepsilon\rightarrow 0 }\mathcal{S}_{\gamma(b-\varepsilon)}
        \big(q-\gamma(b-\varepsilon)\big)
        f\big(\gamma(b-\varepsilon)\big)=f(q).
    \]
    Since the assumptions $(\mathscr{P}_{1})$-$(\mathscr{P}_{3})$ are satisfied,
    then
    \[
        \big\Vert A(q)\odot^{R}\mathcal{S}_{q}(q-r)-A(r)\odot^{R}\mathcal{S}_{r}
        (q-r)\big\Vert_{\mathcal{B}^{R}(X)}
        \leq   C_{1}|q-r|^{\alpha-1},
    \]
    which implies that
    \begin{align*}
        \lim\limits_{\varepsilon \rightarrow 0 }\mathscr{I}_{2}
        =&~\lim\limits_{\varepsilon\rightarrow 0 }
        \int_{\gamma_{[\phi_{I}(z_{0}),\gamma(b-\varepsilon), I ]}}\big[A(q)\odot^{R}
        \mathcal{S}_{q}(q-r)-A(r)\odot^{R}
        \mathcal{S}_{r}(q-r)\big] f(r)~\mathrm{d}r\\
        =&~\int_{\gamma_{[\phi_{I}(z_{0}),q, I ]}}\big[A(q)\odot^{R}\mathcal{S}_{q}(q-r)-A(r) \odot^{R}
        \mathcal{S}_{r}(q-r)\big]
        f(r)~\mathrm{d}r.
    \end{align*}
    As for $\mathscr{I}_{3}$, we have
    \begin{align}\label{ua4}
        \begin{split}
            &~\int_{\gamma_{[\phi_{I}(z_{0}),\gamma(b-\varepsilon), I ]}}\big(A(q) \odot^{R}
            \mathcal{S}_{q}(q-r) \big)f(r)~\mathrm{d}r\\
            =&~\int_{\gamma_{[\phi_{I}(z_{0}),\gamma(b-\varepsilon), I ]}}\big(A(q) \odot^{R}
            \mathcal{S}_{q}(q-r)\big)
            \big(f(r)-f(q)\big)~\mathrm{d}r\\
            &~+
            \int_{\gamma_{[\phi_{I}(z_{0}),\gamma(b-\varepsilon), I ]}}
            \big(A(q)\odot^{R}\mathcal{S}_{q}(q-r)\big) f(q)~\mathrm{d}r\\
            =&~\int_{\gamma_{[\phi_{I}(z_{0}),\gamma(b-\varepsilon), I ]}}\big(A(q)\odot^{R}
            \mathcal{S}_{q}(q-r) \big)\big(f(r)
             -f(q)\big)~\mathrm{d}r\\
             &~-\big[\mathcal{S}_{q}(q-\phi_{I}(z_{0}))
             -\mathcal{S}_{q}(q-\gamma(b-\varepsilon))\big] f(q).
        \end{split}
    \end{align}
        Since $f$ is right slice regular and
    $\big\Vert A(q)\odot^{R}\mathcal{S}_{q}(q-r)\big\Vert_{\mathcal{B}^{R}(X)}\leq  C_{3}|q-r|^{-1}$, then
    the integrand on the first term on the right-hand side of \eqref{ua4} is integrable and hence
    \begin{align*}
        &~\lim\limits_{\varepsilon\rightarrow 0 }
        \int_{\gamma_{[\phi_{I}(z_{0}),\gamma(b-\varepsilon), I ]}}\big(A(q)\odot^{R}
        \mathcal{S}_{q}(q-r) \big)f(r)~\mathrm{d}r\\
        =&~\int_{\gamma_{[\phi_{I}(z_{0}),q, I ]}}\big(A(q)\odot^{R}\mathcal{S}_{q}(q-r) \big)
        \big(f(r)-f(q)\big)~\mathrm{d}r
         -\mathcal{S}_{q}\big(q-\phi_{I}(z_{0})\big) f(q)+f(q).
    \end{align*}
    Therefore,
    \[
     \frac{\partial_{S}^{R}}{\partial q}\int_{\gamma_{[\phi_{I}(z_{0}),\gamma(b-\varepsilon), I ]}}\mathcal{S}_{r}(q-r) f(r)~\mathrm{d}r
    \]
    converges as $\varepsilon\rightarrow 0$ and noticing that
    \[
     \frac{\partial^{R}_{S}}{\partial q}\int_{\gamma_{[\phi_{I}(z_{0}),\gamma(b-\varepsilon),I]}}\mathcal{W}(q,r) f(r)~\mathrm{d}r
    \]
    converges as $\varepsilon\rightarrow 0$, we have that
    $\lim_{\varepsilon\rightarrow0} \partial_{S}^{R} v_{\varepsilon}(q)$ exists and
    this limit is uniform on $\gamma_{[\gamma(a+\varepsilon),q, I ]}$ for every
    $\varepsilon>0$.

    Moreover, it follows from \eqref{ua2} that $A(q) v_{\varepsilon}(q)$ converges as
    $\varepsilon\rightarrow 0$,
    and by \eqref{ua5}, we have $v_{\varepsilon}(q)\rightarrow v(q)$ as $\varepsilon\rightarrow 0 $.
    Since $A(q)$ is closed,
    then $v(q)\in\mathcal{D}\big(A(q)\big)=\mathfrak{D}$ and
    $A(q) v_{\varepsilon}(q)\rightarrow A(q) v(q)$ for
    $q\in \mathscr{X}_{\theta}\big(\phi_{I}(z_{0})\big)$.
    Integrating \eqref{ua2} along $\gamma_{[\gamma(b-h),q,I]}$ we obtain that
    \begin{align*}
        v_{\varepsilon}(q)-v_{\varepsilon}(\gamma(b-h))=&~
        -\int_{\gamma_{[\gamma(b-h),q,I]}}A(\xi) v_{\varepsilon}(\xi)~\mathrm{d}\xi\\
        &~+
        \int_{\gamma_{[\gamma(b-h),q,I]}}
        \mathcal{U}\big(\xi,\gamma(\gamma^{-1}(\xi)-\varepsilon)\big)
        f\big(\gamma(\gamma^{-1}(\xi)-\varepsilon)\big)~\mathrm{d}\xi.
    \end{align*}
    Passing to the limit as $\varepsilon\rightarrow 0 $ yields
    \begin{equation}\label{ua6}
        v (q)-v \big(\gamma(b-h)\big)
        =-\int_{\gamma_{[\gamma(b-h),q,I]}}A(\xi) v(\xi)~\mathrm{d}\xi
        +\int_{\gamma_{[\gamma(b-h),q,I]}}f(\xi)~\mathrm{d}\xi.
    \end{equation}
    Lastly, multiplying right by $[q-\gamma(b-h)]^{-1}$  on the both side of \eqref{ua6}
    and letting $h\rightarrow 0 $, it yields that
    $v(q)$ satisfies \eqref{ua8}. On the other hand,
    the uniqueness of the solution to the quaternionic Cauchy problem \eqref{ua0}
    follows directly from the uniqueness of the solution to the quaternionic
    homogeneous Cauchy problem \eqref{ua7}, as demonstrated in Theorem \ref{ht3}.
    This completes the proof.
\end{proof}

\section{Cauchy Problems in Banach Bilateral $\mathbb{A}$-Modules}\label{sec5}

This section extends the Cauchy problems previously discussed for the Banach bilateral $\mathbb{H}$-modules
to the Banach $\mathbb{A}$-bimodules.
Since Banach bilateral $\mathbb{A}$-modules can be regarded as a general algebraic extension of
Banach bilateral $\mathbb{H}$-modules,  all our results can be generalized to Banach bilateral
modules having as a set of scalar any real associative $^*$-algebra, Clifford algebras $\mathbb{R}_{n}$ included.
In the following, we will lift the results to the Banach bilateral $\mathbb{A}$-modules and $\mathbb{A}$-algebras. Since the proofs of all the results in this section are analogous to those presented in Section \ref{sec4}, we will not repeat these redundant proof details.

\subsection{$\sigma$-Analyticity in Banach Bilateral $\mathbb{A}$-Modules}

For this purpose, we first introduce some basic definitions.
For a finite dimensional associative nontrivial real algebra $\mathbb{A}$
with unit element $1_{\mathbb{A}}$,
we define the $\ast$-involution as follows:
\begin{align*}
    (q^{c})^{c}=q,\quad (pq)^{c}=q^{c}p^{c},\quad r^{c}=r\quad
    (\forall q,p\in \mathbb{A},r\in\mathbb{R}).
\end{align*}
The scalar components of an element $q \in \mathbb{A}$, namely its real part $\mathrm{Re}^{\mathbb{A}}(q)$ and imaginary part $\mathrm{Im}^{\mathbb{A}}(q)$, are given by the following relations:
\[
    \mathrm{Re}^{\mathbb{A}}(q) := \frac{1}{2}(q + q^c) \quad \text{and} \quad \mathrm{Im}^{\mathbb{A}}(q)
    := \frac{1}{2}(q - q^c).
\]
In addition, the imaginary sphere and slice subspaces of $\mathbb{A}$ are defined by:
\[
    \mathbb{S}_{\mathbb{A}}:=\{q\in \mathbb{A}:q^{c}=-q,q^{2}=-1\},
\]
and
\[
    \mathbb{C}_{\mathfrak{V}}:=\{a+b\mathfrak{V}\in \mathbb{A}:a,b\in\mathbb{R}\},
    \quad \mathfrak{V}\in \mathbb{S}_{\mathbb{A}}.
\]
Furthermore, we define the quadratic cone $Q_{\mathbb{A}}$
by setting:
\[
    Q_{\mathbb{A}}:=\bigcup_{\mathfrak{V}\in\mathbb{S}_{\mathbb{A}}}\mathbb{C}_{\mathfrak{V}}
    \text{ if }\mathbb{S}_{\mathbb{A}} \neq \emptyset,\quad \text{and}
    \quad Q_{\mathbb{A}}:=\mathbb{R}\text{ otherwise}.
\]
For more information on the algebra $\mathbb{A}$, please refer to \cite{kexi-srs0}.

The isomorphism
$\phi_{I}$ and slice axially symmetric domain $\Omega_{\mathscr{D}}$
(see Definition \ref{kexi-def3.3})
was introduced to construct slice regular
functions and slice product.
It is easily observed that this isomorphism is also applicable to the
algebra $\mathbb{A}$. In this case,
the isomorphism $\phi_{\mathfrak{V}}$ means that
\[
    \phi_{\mathfrak{V}}(a + b i ) := a + b\mathfrak{V},
    \quad (a, b \in \mathbb{R},\mathfrak{V}\in\mathbb{S_{A}}).
\]
Let $\mathscr{D}$ be a subset of $\mathbb{C}$,
    symmetric with respect to complex conjugation. Then,
    the subset $\Theta_{\mathscr{D}}$ of $Q_{\mathbb{A}}$ is defined as follows:
    \[
        \Theta_{\mathscr{D}} := \bigcup_{\mathfrak{V}\in \mathbb{S_{A}}}
        \phi_{\mathfrak{V}}(\mathscr{D})
        = \{a + b\mathfrak{V} \in Q_{\mathbb{A}} : a, b \in \mathbb{R}, a + b i  \in \mathscr{D},
        \mathfrak{V} \in \mathbb{S_{A}}\}.
    \]

In \cite{kexi-srs0}, Ghiloni and Recupero considered $\mathbb{A}$
to be a finite dimensional, associative, nontrivial real Banach algebra with
a unit element $1_{\mathbb{A}}$.
In this case, an abelian group $(X, +)$ is defined to
have the structure of a right $\mathbb{A}$-module if, when equipped with a scalar
multiplication $X \times \mathbb{A} \rightarrow X$, it satisfies the distributive
and associative laws, specifically ensuring that $(xp)q = x(pq)$ and $x1 = x$ for all
$x \in X$ and $p, q \in \mathbb{A}$.
The definition of a left $\mathbb{A}$-module is given in an entirely analogous and
symmetric manner.
An abelian group $(X, +)$ is then called a bilateral (or two-side)
$\mathbb{A}$-module if it is simultaneously equipped with both left and
right scalar multiplications that satisfy $p(xq) = (px)q$ and
$rx = xr$ for all $r \in \mathbb{R}$ and $x \in X$.

Building upon this algebraic foundation,
$X$ is defined as a (associative) bilateral
$\mathbb{A}$-algebra if it is endowed with an associative product
$X\times X\rightarrow X:(x,y)\mapsto xy$ such that
\begin{align*}
    &x(y+z)=xy+xz,\quad (x+y)z=xz+yz\quad \forall x,y,z\in X,\\
    &q(xy) = (qx)y,\quad (xy)q=x(yq)\quad \forall x,y\in X,\forall q\in\mathbb{A}.
\end{align*}
For such an algebra, if the norm $|\cdot|$ of $\mathbb{A}$ has the property
$|pq|=|p||q|$ on specific subspaces $\mathbb{C}_{\mathfrak{V}}$ $(\mathfrak{V}\in\mathbb{S_{A}})$, the structure
may be enriched into a normed bilateral $\mathbb{A}$-algebra provided its
$\mathbb{A}$-norm $\|\cdot\|$ is submultiplicative. In cases where $X$ is complete,
it is called a Banach bilateral $\mathbb{A}$-algebra, which may further
be required to possess a unit element $1_X$ of unit norm.

\begin{remark}
In the following discussion, we will assume that
$\mathbb{A}$ is
    a finite dimensional associative nontrivial real Banach algebra
    with unit element $1_{\mathbb{A}}$ and $\mathbb{S_{A}}\neq \emptyset$
    such that
    the norm $|\cdot|$ of $\mathbb{A}$
    obeys the property $|pq|=|p||q|$ whenever $p,q\in\mathbb{C}_{\mathfrak{V}}$
    for $\mathfrak{V}\in\mathbb{S_{A}}$. In addition,
$X$ is a Banach bilateral
$\mathbb{A}$-module with $\mathbb{A}$-norm $\|\cdot\|$.
Furthermore,
$\mathscr{L}^{r}(X)$ denote the set of right linear bounded operators $A$ with $\mathcal{D}(A)=X$.
\end{remark}

Next, we discuss the concept of analyticity over Banach bilateral $\mathbb{A}$-algebras.
Following the approach in \cite{kexi-jishu}, which incorporates the multiplicative
structure from \cite{kexi-jishu-07} to establish the definition of
quaternionic power series and the concept of $\sigma$-analyticity,
we define a corresponding multiplication within the framework of
Banach bilateral $\mathbb{A}$-modules.
To begin with, we introduce a specific multiplicative structure that will
be utilized in the series representation.

\begin{definition}\label{kexi-jishu-09}
    Suppose that
    $f:U\subset\mathbb{A}\rightarrow X$ and $g:U\subset\mathbb{A}\rightarrow X$ are two functions
    satisfying:
    \begin{align*}
        f(v)=\sum_{n\in\mathbb{N}}a_{n}v^{n}\quad \text{and} \quad
        g(v)=\sum_{n\in\mathbb{N}}b_{n}v^{n}\quad (v,a_{n},b_{n}\in \mathbb{A}),
    \end{align*}
    then define
    \begin{align*}
        (f\ast_{r}g)(v)=\sum_{\ell=0}^{n}a_{\ell}b_{n-\ell}\sum_{n\in\mathbb{N}}v^{n}.
    \end{align*}
    Furthermore, for a given $p\in\mathbb{A}$ and a sequence
    $\{c_{n}\}_{n\in\mathbb{N}}\subset X$.
    The power series centered at $p$ with coefficients $\{c_{n}\}_{n\in\mathbb{N}}$ is
    defined as
    \begin{equation}\label{kexi-jishu-08}
        \sum_{n\in\mathbb{N}}c_{n}(v-p)^{\ast_{r}n}\quad (v\in\mathbb{A}),
    \end{equation}
    where $(v-p)^{\ast_{r}n}=(v-p)\ast_{r}\cdots \ast_{r}(v-p)$ denotes the $n$-th power of
    the binomial with respect to $\ast_{r}$-multiplication.
\end{definition}

The multiplicative structure introduced in Definition
\ref{kexi-jishu-09} is crucial for both the representation of power series
and the verification of slice regularity,
serving as a fundamental approach to investigating noncommutative algebraic structures.

In \cite{kexi-jishu}, Gentili and Stoppato
introduced the notion of strong analyticity by means of quaternionic power series.
Subsequently, they defined a distance function $\sigma$ on $\mathbb{H}$
that is not topologically equivalent to the Euclidean distance.
The topology induced by $\sigma$ is finer than the Euclidean one.
Furthermore, they established that the convergence domain of a quaternionic power series
is a $\sigma$-ball. Analogously,
the function $\sigma$ is introduced over the $\mathbb{A}$-module structure.

\begin{definition}
    For any $p,v\in Q_{\mathbb{A}}$, define
    \[
        \omega(v,p)=\sqrt{\big(\mathrm{Re}^{\mathbb{A}}(v)-\mathrm{Re}^{\mathbb{A}}(p)\big)^{2}+\big(|\mathrm{Im}^{\mathbb{A}}(v)|+|\mathrm{Im}^{\mathbb{A}}(p)|\big)^{2}}
    \]
    and
    \begin{equation*}
        \sigma(v,p)=\left\{
            \begin{aligned}
                &|v-p|\quad &&\text{if }p,v \in  \mathbb{C}_{\mathfrak{V}}\text{ for some }
                \mathfrak{V}\in\mathbb{S_{A}},\\
                & \omega(v,p)&& \text{otherwise}.
            \end{aligned}
        \right.
    \end{equation*}
    In addition, define the $\sigma$-ball $\Sigma_{\mathbb{A}}(p,R)$ as follows:
    \[
        \Sigma_{\mathbb{A}}(p,R):=\{v\in Q_{\mathbb{A}}:\sigma(v,p)<R,R\in\mathbb{R}^{+}\}.
    \]
\end{definition}

\begin{remark}
    Note that the definition of the function $\sigma$ is restricted to $Q_{\mathbb{A}}$.
    This restriction is determined by the construction of $\Theta_{\mathscr{D}}$ and the
    assumption that $\mathbb{S_{A}} \neq \emptyset$. Such a setting ensures that
    $\sigma$ serves as a well-defined metric on
    $Q_{\mathbb{A}}$.
\end{remark}
Similar to Theorem 8 of \cite{kexi-jishu},
the convergence region of the series \eqref{kexi-jishu-08}
is found to be a $\sigma$-ball $\Sigma_{\mathbb{A}}(p,R)$ with center $p$,
where $1/R=\lim\sup_{n\rightarrow\infty}\|c_{n}\|^{1/n}$.
Moreover, for the power series $f(v)=\sum_{n\in\mathbb{N}}c_{n}v^{n}$, if we set
the derivative of $f$ by
\[
    f'(v)=\sum_{n\in\mathbb{N}}(n+1)c_{n+1}v^{n},
\]
then we have the following result directly.

\begin{theorem}
    Let $f(v)=\sum_{n\in\mathbb{N}}c_{n}v^{n}$ be a
    power series with radius of convergence $R>0$. If
    $p\in \{v\in Q_{\mathbb{A}}:|v|<R\}$. Then
    \begin{equation*}
        f(v)=\sum_{n\in\mathbb{N}}\frac{1}{n!}f^{(n)}(p)(v-p)^{\ast_{r}n},
    \end{equation*}
    where $v\in \Sigma_{\mathbb{A}}(p,R-|p|)$.
\end{theorem}

Now, based on the preceding discussion,
we are in a position to introduce the concept of $\sigma$-analyticity
in Banach bilateral $\mathbb{A}$-modules.

\begin{definition}
   Let $U$ be a $\sigma$-open ball of $Q_{\mathbb{A}}$ and $X$ be a
   Banach bilateral $\mathbb{A}$-module.
   A function $f:U\rightarrow X$ is called right
    $\sigma$-analytic at $p\in U$ if there exists power series
    $\sum_{n\in\mathbb{N}}c_{n}(v-p)^{{\ast_{r}} n}$ $(c_{n}\in X,v,p\in Q_{\mathbb{A}})$
    converging in a $\sigma$-neighborhood
    $W$ of $p$ in $U$ such that
    $f(v)=\sum_{n\in\mathbb{N}}c_{n}(v-p)^{{\ast_{r}} n}$ for all $v\in W$. Furthermore,
    $f$ is called right $\sigma$-analytic if it is right
    $\sigma$-analytic at all $p\in U$.
\end{definition}

Based on the concept of $X\otimes_{\mathbb{R}}\mathbb{C}$ (see page 5006 in \cite{kexi-srs0})
and stem function (see Definition 4.2 in \cite{kexi-srs0}),
Ghiloni and Recupero introduced the definition of slice regular functions
on Banach bilateral $\mathbb{A}$-modules;
specifically, when $\mathbb{A} = \mathbb{H}$,
it reduces to the definition originally proposed by Alpay et al. in
\cite{kexi-srs1} (see Remark 4.5 in \cite{kexi-srs0}).
Furthermore, the relationship between $\sigma$-analytic functions and
slice regular functions is established by the following theorem.
The proof is completely analogous to that of Corollary 2 in
\cite{kexi-jishu} and is therefore omitted here.

\begin{theorem}
    Let $X$ be a Banach bilateral $\mathbb{A}$-module. Then a function
    $f:U\subset Q_{\mathbb{A}}\rightarrow X$ is right $\sigma$-analytic if
    and only if it is right slice regular in the same domain.
\end{theorem}

\subsection{Evolution System in Banach Bilateral $\mathbb{A}$-Algebra and Cauchy Problem}

Let us define
$\Sigma^{\mathbb{A}}_{\varphi}$ as follows:
\[
    \Sigma^{\mathbb{A}}_{\varphi}:=\big\{s\in Q_{\mathbb{A}}\backslash\{0\}:\mathrm{arg}(s)<\varphi\big\}.
\]
Moreover, we will set
$\widetilde{\Sigma}^{\mathbb{A}}_{\varphi}:=\Sigma^{\mathbb{A}}_{\varphi}\cup\{0\}$.
In this subsection, we investigate the following equation
on a Banach bilateral $\mathbb{A}$-module $X$:
\begin{equation}\label{kexi-e4.1-5}
    \left\{\begin{aligned}
     &\partial_{S}^{R} u (q)
    =A(q)u(q)+f(q), \\
     &u(s_{0})=u_{0},
    \end{aligned}\right.\qquad s_{0},q\in \Sigma^{\mathbb{A}}_{\varphi},
\end{equation}
where $A(q)\in \mathscr{L}^{r}(X)$, $f$ is a given function.

Following the proof technique of Theorem \ref{kexi-t4.1},
we first establish the existence and uniqueness of the solution to
\begin{equation}\label{kexi-ftf-0h}
    \left\{\begin{aligned}
     &\partial_{S}^{R}u (q)=A(q)u(q),  \\
        &u(s_{0})=u_{0},
    \end{aligned}\right. \quad s_{0},q\in\Sigma^{\mathbb{A}}_{\theta},
\end{equation}
on the complex plane $\mathbb{C}_{\mathfrak{V}}$
$(\mathfrak{V}\in\mathbb{S_{A}})$ containing the initial value.
Subsequently, by virtue of the representation formula,
we extend this result to any other complex plane
$\mathbb{C}_{\mathfrak{U}}$ for $\mathfrak{U} \in \mathbb{S_{A}}$.
Consequently, we obtain the following local existence and uniqueness theorem.

\begin{theorem}\label{kexi-ftf4.1}
    Let $X$ be a Banach bilateral $\mathbb{A}$-module and
   $\theta\in(0,\pi/2]$.
    Assume that the operator-valued function
    $A: \Sigma^{\mathbb{A}}_{\theta}\to \mathscr{L}^{r}(X)$
    is right slice regular and $A(q)$ is uniformly bounded closed right linear
    operator for each $q\in\Sigma^{\mathbb{A}}_{\theta}$, then
    \eqref{kexi-ftf-0h}
    has a unique right $\sigma$-analytic solution $u$
    in $\overline{\Sigma_{\mathbb{A}}(s_{0},R)}\subset \Sigma^{\mathbb{A}}_{\theta}$.
\end{theorem}

\begin{remark}
    Note that the set $\Sigma_{\theta}^{\mathbb{A}}$ is a slice axially symmetric domain
    in the sense of $\mathbb{A}$-algebras;
    hence, in what follows, we shall use $\Theta_{\mathscr{D}}$ to
    replace $\Sigma_{\theta}^{\mathbb{A}}$.
\end{remark}

Based on the types of $\sigma$-balls,
we can extend the solution in three distinct cases until it is defined over the entire
domain $\Sigma_{\theta}^{\mathbb{A}}$.
Consequently, we establish the following extension theorem for the solution.

\begin{theorem}
Let $X$ be a
Banach bilateral $\mathbb{A}$-module.
Assume that $A(q)\in\mathscr{L}^{r}(X)$
satisfies the hypotheses of Theorem \ref{kexi-ftf4.1},
then the right $\sigma$-analytic solution of \eqref{kexi-ftf-0h} on
$\overline{\Sigma_{\mathbb{A}}(s_{0},R)}$
can be extended to a global solution on $\Theta_{\mathscr{D}}$.
\end{theorem}

Having established the existence and uniqueness of the solution to the homogeneous equation
\eqref{kexi-ftf-0h}, we now provide a unified expression for the solution via
the isomorphism $\phi_{\mathfrak{V}}$.

\begin{proposition}
    Let $X$ be a Banach bilateral $\mathbb{A}$-module.
Assume that $A(q)$
satisfies the hypotheses of Theorem \ref{kexi-ftf4.1} and
$s_{0}=a+b\mathfrak{U}=\phi_{\mathfrak{U}}(z_{0})$
$(\mathfrak{U}\in\mathbb{S_{A}})$. Then the classical solution
of \eqref{kexi-ftf4.1} admits the following unified form:
\begin{equation*}
    u\big(\phi_{\mathfrak{V}}(z)\big)
    =u_{\mathfrak{V},0}+\int_{\gamma_{[\phi_{\mathfrak{V}}(z_{0}),\phi_{\mathfrak{V}}(z),\mathfrak{V}]}}
    A(\xi)u(\xi)~\mathrm{d}\xi,
\end{equation*}
where $\phi_{\mathfrak{V}}(z_{0})=a+b\mathfrak{V}$, $z\in\mathscr{D}$ and
\[
    u_{\mathfrak{V},0}=u\big(\phi_{\mathfrak{V}}(z_{0})\big)=
    \left\{\begin{aligned}
    &u_{0}, &&\mathfrak{V}=\mathfrak{U},\\
    &\frac{1}{2}u(\overline{s_{0}})(1+\mathfrak{U}\mathfrak{V})+\frac{1}{2}u(s_{0})(1-\mathfrak{U}\mathfrak{V}),\quad &&\mathfrak{V}\neq \mathfrak{U}.
    \end{aligned}\right.
\]
\end{proposition}

\begin{remark}
    It is worth noting that in the cross-slice representation of solutions in both
    Theorem \ref{kexi-t4.1} and Theorem \ref{kexi-ftf4.1},
    we initiate the process from a $\sigma$-ball
    at the initial time $s_{0}$, utilizing a point and its conjugate.
    However, such a ball does not necessarily contain both a point and its
    conjugate simultaneously. To address this, we can first extend the solution within the
    complex plane containing the initial point, and then perform the cross-slice
    representation. This approach making our presentation clear and consistent.
\end{remark}

To establish the notion of evolution systems in the sense of Banach
bilateral $\mathbb{A}$-modules,
we require the following definition.

\begin{definition}\label{kexi-tianjia-0}
    For $\theta\in(0,\pi/2]$ and $p\in \mathbb{C}_{\mathfrak{V}}$ with $\mathfrak{V}\in\mathbb{S_{A}}$, define the set 
    \[
        \mathscr{X}^{\mathbb{A}}_{\theta}(p)
        :=\big\{q\in  \mathbb{C}_{\mathfrak{V}}: 
        \mathrm{Re}^{\mathbb{A}}(q)\geq\mathrm{Re}^{\mathbb{A}}(p) \mbox{ and }
        |\mathrm{Im}^{\mathbb{A}}(q)-\mathrm{Im}^{\mathbb{A}}(p)|\leq 
        \big(\mathrm{Re}^{\mathbb{A}}(q)-\mathrm{Re}^{\mathbb{A}}(p)\big)\tan \theta
        \big\}.
    \]
\end{definition}

\begin{remark}
    Actually, for an arbitrary  $w\in \mathbb{A}$, 
    it is not true in general that 
    $\mathrm{Re}^{\mathbb{A}}(w) \in \mathbb{R}$. 
    However, in Definition \ref{kexi-tianjia-0}, 
    since the chosen $p$ and $q$ both lie in some  $\mathbb{C}_{\mathfrak{V}}$,
    it follows from Proposition 2.3 in \cite{kexi-srs0} that 
    $\mathrm{Re}^{\mathbb{A}}(p), \mathrm{Re}^{\mathbb{A}}(q) \in \mathbb{R} $. 
    Therefore, Definition \ref{kexi-tianjia-0} is well-defined.
\end{remark}

To further represent the solution, we now introduce the definition
of evolution systems of
\eqref{kexi-ftf-0h} in the sense of Banach bilateral $\mathbb{A}$-modules.

\begin{definition}\label{kexi-def4-5}
    Let $u(q)$ be the solution of \eqref{kexi-ftf-0h} in $\Theta_{\mathscr{D}}$
    and $s_{0}=\phi_{\mathfrak{U}}(z_{0})$ $(z_{0}\in\mathscr{D},\mathfrak{U}\in\mathbb{S_{A}})$.
    A two-parameter operator-valued function
    $\mathcal{O}:\Theta_{\mathscr{D}}\times \Theta_{\mathscr{D}}\rightarrow \mathscr{L}^{r}(X)$
    is called the evolution system of \eqref{kexi-ftf-0h} if
    $(\mu,\varrho)\mapsto\mathcal{O}(\mu,\varrho)$ is right holomorphic with respect to both parameters
    $\mu$ and $\varrho$ as $\mu,\varrho\in\mathbb{C}_{\mathfrak{L}}$ for every $\mathfrak{L}\in\mathbb{S_{A}}$
    and the following relations hold
    \[
        \mathcal{O} \big(\phi_{\mathfrak{V}}(z),\phi_{\mathfrak{V}}(z_{0})\big)
        u\big(\phi_{\mathfrak{V}}(z_{0})\big)
        =u\big(\phi_{\mathfrak{V}}(z)\big) \quad
        \big(\mathfrak{V}\in\mathbb{S_{A}},z\in  \mathscr{D},\phi_{\mathfrak{V}}(z)\in
        \mathscr{X}^{\mathbb{A}}_{\theta}(\phi_{\mathfrak{V}}(z_{0}))\big)
    \]
    and
    \[
        \mathcal{O}\big(\phi_{\mathfrak{V}}(z_{1}),\phi_{\mathfrak{V}}(z_{1})\big)=\mathcal{I}, \quad
       \mathcal{O}\big(\phi_{\mathfrak{V}}(z_{2}),\phi_{\mathfrak{V}}(z_{1})\big)\circledcirc^{R}
       \mathcal{O}\big(\phi_{\mathfrak{V}}(z_{1}),\phi_{\mathfrak{V}}(z_{0})\big)
       =\mathcal{O}\big(\phi_{\mathfrak{V}}(z_{2}),\phi_{\mathfrak{V}}(z_{0})\big),
    \]
    where $z_{1},z_{2}\in\mathscr{D}$,
    $\phi_{\mathfrak{V}}(z_{1})\in \mathscr{X}^{\mathbb{A}}_{\theta}\big(\phi_{\mathfrak{V}}(z_{0})\big)$,
    $\phi_{\mathfrak{V}}(z_{2})\in \mathscr{X}^{\mathbb{A}}_{\theta}\big(\phi_{\mathfrak{V}}(z_{1})\big)$ and $\circledcirc^{R}$ denote the right slice product.
\end{definition}

Thus, based on the definition of evolution systems,
we can formulate the expression for the solution to the
inhomogeneous equation \eqref{kexi-e4.1-5} over the restricted set.

\begin{proposition}
    Let $X$ be a Banach bilateral $\mathbb{A}$-module and
    $f:\Theta_{\mathscr{D}}\to X$ be right slice regular function.
    Assume that the uniformly bounded operator-valued function
    $A: \Theta_{\mathscr{D}} \to \mathscr{L}^{r}(X)$
    is right slice regular and $s_{0}=\phi_{\mathfrak{U}}(z_{0})$
    $(z_{0}\in\mathscr{D},\mathfrak{U}\in\mathbb{S_{A}})$.
    If there is an
    evolution system $\mathcal{O}(\mu,\varrho)$ of the
    Cauchy problem \eqref{kexi-ftf-0h} such that 
    the inclusion
    $u_{\mathfrak{V},0}\in\mathcal{D}\big(A(\phi_{\mathfrak{V}}(z_{0}))\big)$ 
    implies that
    $\mathcal{O}\big(q,\phi_{\mathfrak{V}}(z_{0})\big) u_{\mathfrak{V},0}\in\mathcal{D}\big(A(q)\big)$ for all $q\in\mathscr{X}^{\mathbb{A}}_{\theta}\big(\phi_{\mathfrak{V}}(z_{0})\big)$ 
    and
    $\mathcal{O}\big(q,\phi_{\mathfrak{V}}(z_{0})\big) u_{\mathfrak{V},0}$
    admits right slice derivative both in $q$
    and $\varrho$ $\big(\varrho=\phi_{\mathfrak{V}}(z_{0})\big)$ satisfying
    \begin{align*}
         \frac{\partial_{S}^{R}}{\partial q}
         \mathcal{O}\big(q,\phi_{\mathfrak{V}}(z_{0})\big) u_{\mathfrak{V},0}
        = &~\big( A(q) \circledcirc^{R} \mathcal{O}(q,\phi_{\mathfrak{V}}(z_{0})) \big) u_{\mathfrak{V},0},\\
        \frac{\partial_{S}^{R}}{\partial \varrho}\mathcal{O}(q,\varrho)
        u_{\mathfrak{V},0}
        = &~\big(-\mathcal{O}(q,\varrho)\circledcirc^{R} A(\varrho) \big) u_{\mathfrak{V},0},
    \end{align*}
    then every classical solution $u$ of \eqref{kexi-e4.1-5}
    with $u_{\mathfrak{V},0}\in\mathcal{D}\big(A(\phi_{I}(z_{0}))\big)$,
    restricted to
    \[
    \big\{p\in\Theta_{\mathscr{D}}:p\in \mathscr{X}^{\mathbb{A}}_{\theta}(\phi_{\mathfrak{V}}(z_{0})),\mathfrak{V}\in\mathbb{S_{A}}\big\},
    \]
    is given by
    \begin{equation*}
        u(q)=
        \mathcal{O}\big(q,\phi_{\mathfrak{V}}(z_{0})\big) u_{\mathfrak{V},0}
        +\int_{\gamma_{[\phi_{\mathfrak{V}}(z_{0}),q,\mathfrak{V}]}}
        \mathcal{O}(q,\xi)  f(\xi)~ \mathrm{d}\xi,\quad
        q \in  \mathscr{X}^{\mathbb{A}}_{\theta}\big(\phi_{\mathfrak{V}}(z_{0})\big).
    \end{equation*}
\end{proposition}

Finally, by generalizing conditions $(\mathfrak{H}_1)$-$(\mathfrak{H}_3)$ from Subsection
\ref{subsec5}
and $(\mathscr{P}_1)$-$(\mathscr{P}_3)$ from Subsection
\ref{subsec6}, we establish the existence and
uniqueness of evolution systems under these two distinct frameworks, respectively.
Based on these results, we further derive the existence and uniqueness of the solutions
to the corresponding equations.
First, regarding the generalized conditions
$(\mathfrak{H}_1)^{\mathbb{A}}$-$(\mathfrak{H}_3)^{\mathbb{A}}$, the following results hold.

\begin{theorem}\label{ht1-5}
    Let the family $\{A(q)\}_{q\in \Theta_{\mathscr{D}}}$
    satisfy the stipulated hypotheses
    $(\mathfrak{H}_{1})^{\mathbb{A}}$-$(\mathfrak{H}_{3})^{\mathbb{A}}$. Then
    there exists a unique evolution
    system $ \mathcal{O}(\mu, \varrho)$,
    $(\mu,\varrho)\in \Theta_{\mathscr{D}} \times \Theta_{\mathscr{D}}$
    of \eqref{kexi-ftf-0h},
    such that
    \begin{description}
        \item[$(\mathfrak{E}_{1})^{\mathbb{A}}$]
        $\Vert\mathcal{O}(\mu,\varrho)\Vert_{\mathscr{L}^{r}(X)} \leq  M e^{\omega\cdot\mathrm{Re}^{\mathbb{A}}(\mu-\varrho)}$ for $\varrho=\phi_{\mathfrak{V}}(z_{0})$,
        $\mu\in \mathscr{X}^{\mathbb{A}}_{\theta}\big(\phi_{\mathfrak{V}}(z_{0})\big)$ and
        $\mathfrak{V}\in\mathbb{S_{A}}$,

        \item[$(\mathfrak{E}_{2})^{\mathbb{A}}$]
        $\frac{\partial_{S}^{R}}{\partial \mu}\mathcal{O}\big(\mu,\phi_{\mathfrak{V}}(z_{0})\big) u_{\mathfrak{V},0} =\big(A(\mu)\circledcirc^{R}\mathcal{O}(\mu,\phi_{\mathfrak{V}}(z_{0}))\big)u_{\mathfrak{V},0}$
        for $u_{\mathfrak{V},0}\in \mathfrak{D}^{\mathbb{A}}$,
        $\mu \in \mathscr{X}^{\mathbb{A}}_{\theta}\big(\phi_{\mathfrak{V}}(z_{0})\big)$ and $\mathfrak{V}\in\mathbb{S_{A}}$,

        \item[$(\mathfrak{E}_{3})^{\mathbb{A}}$]
        $\frac{\partial_{S}^{R}}{\partial \varrho}\mathcal{O}(\mu,\varrho) u_{\mathfrak{V},0}=\big(-\mathcal{O}(\mu,\varrho)\circledcirc^{R} A(\varrho)\big) u_{\mathfrak{V},0}$ for $u_{\mathfrak{V},0}\in \mathfrak{D}^{\mathbb{A}}$, $\varrho=\phi_{\mathfrak{V}}(z_{0})$,
        $\mu\in \mathscr{X}^{\mathbb{A}}_{\theta}\big(\phi_{\mathfrak{V}}(z_{0})\big)$ and $\mathfrak{V}\in\mathbb{S_{A}}$,
    \end{description}
where $\mathfrak{D}^{\mathbb{A}}$ is a subspace of
    Banach bilateral $\mathbb{A}$-module $X$ with a dense and continuous embedding in $X$.
\end{theorem}

\begin{corollary}
    Let $\{A(q)\}_{q\in\Theta_{\mathscr{D}}}$ satisfy the conditions
    $(\mathfrak{H}_{1})^{\mathbb{A}}$-$(\mathfrak{H}_{3})^{\mathbb{A}}$.
    Then the Cauchy problem
    \eqref{kexi-ftf-0h} associated with each $u_{0}\in  \mathfrak{D}^{\mathbb{A}}$
    has a unique $\mathfrak{D}^{\mathbb{A}}$-valued solution on $\Theta_{\mathscr{D}}$.
\end{corollary}

Lastly, under the generalized conditions
$(\mathscr{P}_1)^{\mathbb{A}}$-$(\mathscr{P}_3)^{\mathbb{A}}$, we obtain the following results.

\begin{theorem}\label{ht4-5}
    If the assumptions $(\mathscr{P}_{1})^{\mathbb{A}}$-$(\mathscr{P}_{3})^{\mathbb{A}}$ hold, then
    there is a unique evolution
    system $\mathcal{O}(\mu,\varrho)$ such that 
    \begin{description}
        \item[$(\mathfrak{E}_{1})^{+}_{\mathbb{A}}$] $\Vert\mathcal{O}(\mu,\varrho)\Vert_{\mathscr{L}^{r}(X)} \leq  C$
        for $\mu\in \mathscr{X}^{\mathbb{A}}_{\theta}(\varrho)$,

        \item[$(\mathfrak{E}_{2})^{+}_{\mathbb{A}}$] 
        the map
        $\mu\mapsto\mathcal{O}\big(\mu,\phi_{\mathfrak{V}}(z_{0})\big)$ 
        with 
        $\mu\in \mathscr{X}^{\mathbb{A}}_{\theta}(\phi_{\mathfrak{V}}\big(z_{0})\big)$ and
        $\mathcal{O}\big(\mu,\phi_{\mathfrak{V}}(z_{0})\big):X\rightarrow  \mathfrak{D}^{\mathbb{A}}$  
        is right hyperholomorphic in $X$ for every $\mathfrak{V}\in\mathbb{S_{A}}$
        and, moreover,
        \begin{align*}
            &\frac{\partial_{S}^{R}}{\partial \mu}\mathcal{O}\big(\mu,\phi_{\mathfrak{V}}(z_{0})\big)
            +A(\mu)\circledcirc^{R}\mathcal{O}\big(\mu,\phi_{\mathfrak{V}}(z_{0})\big)=0 ,
            \\
            &\left\Vert\frac{\partial_{S}^{R}}{\partial \mu}\mathcal{O}\big(\mu,\phi_{\mathfrak{V}}(z_{0})\big)
            \right\Vert_{\mathscr{L}^{r}(X)}=\Vert A(\mu)\circledcirc^{R}\mathcal{O}\big(\mu,\phi_{\mathfrak{V}}(z_{0})\big)\Vert_{\mathscr{L}^{r}(X)}
             \leq C\big|\mu-\phi_{\mathfrak{V}}(z_{0})\big|^{-1}
             \\
            &\big\Vert A(\mu)\circledcirc^{R}\mathcal{O}\big(\mu,\phi_{\mathfrak{V}}(z_{0})\big)\circledcirc^{R}
            A\big(\phi_{\mathfrak{V}}(z_{0})\big)^{-1}\big\Vert_{\mathscr{L}^{r}(X)} \leq  C,
        \end{align*}
        and  
        \item[$(\mathfrak{E}_{3})^{+}_{\mathbb{A}}$] For every $u_{\mathfrak{V},0}\in\mathfrak{D}^{\mathbb{A}}$,
        $\mu\in\mathscr{X}^{\mathbb{A}}_{\theta}\big(\phi_{\mathfrak{V}}(z_{0})\big)$,
        $\mathcal{O}(\mu,\varrho) u_{\mathfrak{V},0}$ $\big(\varrho=\phi_{\mathfrak{V}}(z_{0})\big)$
        admits right slice derivative with respect to
        $\varrho$ and
        \begin{equation*}
            \frac{\partial_{S}^{R}}{\partial \varrho}\mathcal{O}
            (\mu,\varrho) u_{\mathfrak{V},0}
            =\big(\mathcal{O}(\mu,\varrho)\circledcirc^{R}A(\varrho)\big)
            u_{\mathfrak{V},0}.
        \end{equation*}
    \end{description}
\end{theorem}

\begin{theorem}
    Let $\{A(q)\}_{q\in\Theta_{\mathscr{D}}}$ satisfy the conditions
    $(\mathscr{P}_{1})^{\mathbb{A}}$-$(\mathscr{P}_{3})^{\mathbb{A}}$.
    Then for each $u_{0}\in X$,
    the Cauchy problem
    \begin{equation*}
        \left\{\begin{aligned}
            & \partial_{S}^{R}  u(q)+A(q)u(q)=0, \\
            &u(s_{0})=u_{0},
        \end{aligned}\right.\qquad s_{0},q\in\Theta_{\mathscr{D}}
    \end{equation*}
    has a unique solution $u$ and its restriction to
    \[
    \big\{p\in\Theta_{\mathscr{D}}:p\in \mathscr{X}^{\mathbb{A}}_{\theta}(\phi_{\mathfrak{V}}(z_{0}))\big\}
    \]
    is given by
    $u(q)=\mathcal{O}\big(q,\phi_{\mathfrak{V}}(z_{0})\big) u_{\mathfrak{V},0}$ for
    $q\in \mathscr{X}^{\mathbb{A}}_{\theta}\big(\phi_{\mathfrak{V}}(z_{0})\big)$.
\end{theorem}

\section{Conclusion}\label{sec6}

In this paper, we have answered the key questions raised in the introduction part.
We have established precise sufficient conditions that guarantee the existence and uniqueness of a local $\sigma$-analytic solution to the evolution equation \eqref{99p2}. These conditions consist of the combination of one-parameter regularity along a fixed slice and $\sigma$-analyticity across distinct slices. Moreover, the standard notion of a quaternionic evolution system is defined via the integral representation obtained from the quaternionic sliced geometric path integral. This definition ensures both existence and uniqueness of the evolution system under the aforementioned regularity and analyticity conditions.

Building on these results, we have proved that the local $\sigma$-analytic solution of the quaternionic Cauchy problem \eqref{99p2} is uniquely extendable from a closed $\sigma$-ball to the entire spherical sectorial region. The extension is achieved by iterative $\sigma$-analytic continuation, and uniqueness follows from the identity principle for slice regular functions. Furthermore, hyperbolic-like and parabolic-like constraints independently guarantee global existence and uniqueness, both constraints lead to the global well-posedness result.

With global well-posedness established, we further demonstrate that after global extension the $\sigma$-analyticity of the solution to the inhomogeneous evolution equation (represented via the standard quaternionic evolution system) is preserved. Nontrivial evolution does not destroy the cross-slice analytic structure; the solution remains $\sigma$-analytic on the whole spherical sectorial region without being forced to a single isolated slice. Finally, the entire well-posedness theory has been successfully lifted from $\mathbb{H}$-algebras to the more general setting of $\mathbb{A}$-algebras. This lifting reveals that any noncommutative Banach algebra equipped with a suitable slice regularity condition inherently admits a universal well-posedness framework for evolution equations, thereby establishing the universality of the theory.

\end{document}